\documentclass[11pt,a4paper]{article}
\usepackage[utf8]{inputenc}
\usepackage{amsmath,bm}
\usepackage{enumitem}
\usepackage{amsmath,nccmath}
\usepackage{amsthm}
\usepackage{amsfonts}
\usepackage[square, comma, sort&compress,numbers]{natbib}
\usepackage{amssymb}
\usepackage{lineno}
\usepackage{setspace}
\usepackage{graphicx}
\usepackage{epstopdf}
\usepackage{color}
\usepackage[linkcolor=blue,colorlinks=true,linktocpage=true]{hyperref}
\usepackage[marginal]{footmisc}
\usepackage[section]{placeins}
\usepackage{indentfirst}
\usepackage[left=2cm,right=2cm,top=2cm,bottom=2cm]{geometry}
\usepackage[title]{appendix}

\author{}
\usepackage{float}
\theoremstyle{plain}
\newtheorem{thm}{Theorem}
\newtheorem*{lem*}{Lemma}
\newtheorem*{thm*}{Theorem}
\newtheorem{lem}{Lemma}[section]
\newtheorem{prop}[lem]{Proposition}
\newtheorem*{prop*}{Proposition}
\newtheorem{cor}{Corollary}

\newtheorem{thma}{Theorem}

\newtheorem{cora}{Corollary}

\newtheorem{propa}{Proposition}[section]

\theoremstyle{definition}
\newtheorem{defi}[lem]{Definition}
\newtheorem*{defi*}{Definition}

\newtheorem*{ques*}{Question}
\newtheorem*{conj*}{Conjecture}
\newtheorem{defia}[propa]{Definition}

\newtheorem{rem}[lem]{Remark}
\newtheorem{rema}[propa]{Remark}

\allowdisplaybreaks
\usepackage{mathtools}

\usepackage{authblk}
\usepackage{titlesec} 
\titleformat{\subsubsection}[runin]
  {\normalfont\normalsize\bfseries}{\thesubsubsection}{1em}{}

\numberwithin{equation}{section}

\definecolor{darkred}{rgb}{.56,0,0}

\makeatletter
\newcommand*{\rom}[1]{\expandafter\@slowromancap\romannumeral #1@}
\makeatother

\title{{\bf Persistence of Heterodimensional Cycles}}
\author{Dongchen Li}
\author{Dmitry Turaev}
\affil{Department of Mathematics, Imperial College London}

\usepackage{tocloft}
\begin{document}
\def\D{\mathrm{D}}
\def\d{\mathrm{d}}
\def\tr{\mathrm{tr}\,}
\def\diff{Di\!f\!f}
\def\eps{\varepsilon}
\bibliographystyle{plainnat}

\renewcommand{\geq}{\geqslant}
\renewcommand{\leq}{\leqslant}
\maketitle

\par{}
\noindent{\bf Abstract.} A heterodimensional cycle is an invariant set of a dynamical system consisting of two hyperbolic periodic orbits with different
dimensions of their unstable manifolds and a pair of orbits that connect them. For systems which are at least $C^2$, we show that bifurcations of a coindex-1 heterodimensional cycle within a generic 2-parameter family create robust heterodimensional dynamics, i.e., a pair of non-trivial hyperbolic basic sets with different numbers of positive Lyapunov exponents, such that the unstable manifold of each of the sets intersects the stable manifold of the second set and these intersections persist for an open set of parameter values. We also give a solution to the so-called local stabilization problem of coindex-1 heterodimensional cycles in any regularity class $r=2,\ldots,\infty,\omega$. The results are based on the observation that arithmetic properties of moduli of topological conjugacy of systems with heterodimensional cycles determine the emergence of Bonatti-D\'iaz blenders.
\tableofcontents

\section{Introduction}\label{sec:intro}

In this paper, we solve the $C^r$-persistence problem for heterodimensional cycles of coindex 1 in any regularity class $r=2,\ldots,\infty,\omega$. The result gives a heterodimensional counterpart to the Newhouse theorem on the $C^r$-persistence of non-transverse equidimensional cycles (homoclinic tangencies) \citep{Ne70}. It implies the ubiquity of heterodimensional dynamics, which is, in our opinion, one of the most basic properties of non-hyperbolic multi-dimensional dynamical systems with chaotic behavior.

For uniformly hyperbolic systems, all orbits within the same chain-recurrent class have the same number of positive Lyapunov exponents and the same number of negative ones. However, chaotic dynamics are often not hyperbolic, and then one can expect that orbits with different numbers of positive Lyapunov exponents coexist and are, in a sense, inseparable from each other. The first example of such sort was given by Abraham and Smale in \cite{AS70}. As an example of the non-density of hyperbolicity in the space of dynamical systems, they described an open region in the space of $C^1$-diffeomorphisms where each diffeomorphism has hyperbolic periodic orbits with different dimensions of unstable manifolds within the same transitive set. More examples followed, see e.g. \cite{Sh71,Si72,Ma78,DR92,Di95a,Di95b,GI}, with a general construction for robust heterodimensionality developed by Bonatti and D\'iaz in \citep{BD96,BD08}.

We study both the discrete-time and continuous-time cases (our approach allows for a simultaneous consideration of both cases, as explained in Section \ref{sec:intro2}). Unless otherwise stated, by a dynamical system, we always mean a diffeomorphism of a manifold of dimension 3 or higher\footnote{Heterodimensional cycles can also be defined for endomorphisms on two-dimensional manifolds, see e.g. \citep{BCP21}.}, or a flow on a manifold of dimension 4 or higher. We use

\begin{defi}[Heterodimensional dynamics]\label{defi:hdd}
Let a dynamical system $f$ have two compact,  transitive, and uniformly hyperbolic invariant sets $\Lambda_1$ and $\Lambda_2$. Let $ind(\Lambda)$ (the index of a transitive hyperbolic set $\Lambda$) denote the dimension of the unstable manifold of any of its orbits\footnote{When a hyperbolic set is transitive, i.e., when some of its orbits is dense in it, all of its orbits have unstable manifolds of the same dimension.}. We say that $f$ has {\em heterodimensional dynamics involving $\Lambda_1$ and $\Lambda_2$} if
\begin{itemize}[nolistsep]
\item $ind(\Lambda_1)\neq ind(\Lambda_2)$; and
\item the unstable sets $W^u(\Lambda_1)$ and $W^u(\Lambda_2)$ intersect the stable sets $W^s(\Lambda_2)$ and, respectively, $W^s(\Lambda_1)$.
\end{itemize}
The difference $|ind(\Lambda_1)-ind(\Lambda_2)|$ is called the coindex of the heterodimensional dynamics.
\end{defi}

Often, the term ``heterodimensional cycle'' is used for what we call the heterodimensional dynamics \citep{BD08,BDK12}. We, however, reserve the term for the most basic case, where $\Lambda_1$ and $\Lambda_2$ are trivial.

\begin{defi}[Heterodimensional cycles]\label{defi:hdc}
A {\em heterodimensional cycle} is a closed invariant set consisting of four orbits: two hyperbolic periodic orbits $L_1$ and $L_2$, with $ind(L_1)\neq ind(L_2)$, and two heteroclinic orbits, one from  $W^s(L_1)\cap W^u(L_2)$ and the other from $W^u(L_1)\cap W^s(L_2)$.
\end{defi}

\begin{figure}[!h]
\begin{center}
\includegraphics[width=0.65\textwidth]{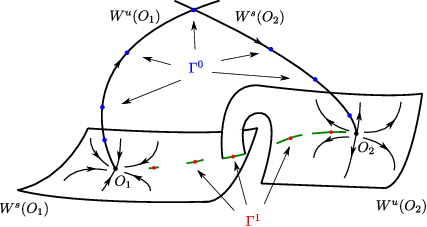}
\end{center}
\caption{A heterodimensional cycle involving two hyperbolic fixed points $O_1$ and $O_2$ (black dots) of a three-dimensional diffeomorphism. The cycle consists of the two fixed points, a fragile heteroclinic orbit $\Gamma^0$ (blue dots) in the non-transverse intersection of the one-dimensional invariant manifolds, and a robust heteroclinic orbit $\Gamma^1$ (red dots) in the transverse intersection (green curves) of the two-dimensional manifolds. See Section \ref{sec:intro2} for details.}
\label{fig:hdc0}
\end{figure}

See Figure \ref{fig:hdc0} for an illustration. Most of this paper is focused on bifurcations in this particular case. This does not lead to a loss in generality, because whenever we have heterodimensional dynamics, a heterodimensional cycle with periodic orbits can be created by a $C^r$-small perturbation (see discussion above Corollary \ref{cor:persis_hdd}).

Due to the difference between the dimensions of the unstable manifolds, some heterodimensional intersections can be fragile. Indeed, consider a diffeomorphism of a $d$-dimensional manifold with a heterodimensional cycle involving periodic orbits $L_1$ and $L_2$. Let ${\rm dim}\, W^u(L_1)=d_1$ and ${\rm dim}\, W^s(L_2)=d-d_2$, with $d_1<d_2$. The space spanned by the tangents to $W^u(L_1)$ and $W^s(L_2)$ at any of their intersection points has dimension less than the dimension $d$ of the full space. This means that every particular heteroclinic intersection of $W^u(L_1)$ and $W^s(L_2)$ is non-transverse and, by Kupka-Smale theorem, can be removed by an arbitrarily small perturbation.

However, the situation changes when we have heterodimensional dynamics involving two non-trivial hyperbolic sets $\Lambda_1$ and $\Lambda_2$. For example, when $ind(\Lambda_2)>ind(\Lambda_1)$, it may happen that when a non-transverse intersection of $W^u(\Lambda_1)$ with $W^s(\Lambda_2)$ at the points of some orbits is destroyed, a new one arises. In this case the heterodimensional dynamics are called {\em robust}.

Recall that basic (i.e., compact, transitive, and locally maximal) hyperbolic invariant sets continue uniquely when the dynamical system varies continuously (in the $C^1$ topology).

\begin{defi}[Robust heterodimensional dynamics]\label{defi:hdd_robust}
We say that a system exhibits $C^1$-{\em robust heterodimensional dynamics} if it has heterodimensional dynamics involving two  hyperbolic basic sets $\Lambda_1$ and $\Lambda_2$, where at least one of them is non-trivial, and there exists a $C^1$-neighborhood $\mathcal{U}$ of the original dynamical system such that every system from $\mathcal{U}$ has heterodimensional dynamics involving the hyperbolic continuations of $\Lambda_1$ and $\Lambda_2$.
\end{defi}

This was the case in the original Abraham-Smale example and in the other examples we mentioned. Moreover, Bonatti and D\'iaz proved in \citep{BD08} that any diffeomorphism with a heterodimensional cycle of coindex 1 can be arbitrarily well approximated, {\em in the $C^1$ topology}, by a diffeomorphism with $C^1$-robust heterodimensional dynamics. The result gave a nice topological characterization of the set of systems with heterodimensional dynamics of coindex 1: this set is the $C^1$-closure of its $C^1$-interior. However, the perturbation techniques used in \cite{BD08} (in an essential way) can only be valid in the $C^1$ topology. As a result, the $C^1$-small perturbations proposed in \cite{BD08} are {\em large in} $C^r$ for any $r>1$. 

This leads to the question whether the Bonatti-D\'iaz result survives in higher regularity\footnote{That the $C^1$ bifurcation theory of heterodimensional cycles cannot be straightforwardly translated to the $C^r$-case is illustrated by the result of \cite{AST17,AST18} which shows that the dynamics created by $C^2$-small or $C^3$-small perturbations of partially-hyperbolic heterodimensional cycles of coindex 1 can be very much different from what can be achieved by $C^1$-small perturbations. The reason is that the dynamics in the central direction are restricted by the signs of the second derivative and the Schwarzian derivative of the one-dimensional transition map; this signature can be changed by $C^1$-small perturbations but not by $C^3$-small perturbations.}. In this paper, we close the question with a positive answer.

\begin{thma}\label{thm:persis_hdd}
Any dynamical system of class $C^r$ ($r=1,\ldots,\infty,\omega$) having a heterodimensional cycle of coindex 1 can be $C^r$-approximated by a system which has $C^1$-robust heterodimensional dynamics.
\end{thma}

\begin{rem}
A partial case of Theorem \ref{thm:persis_hdd} can be derived from the result in \citep{DP19} about renormalization near heterodimensional cycles of three-dimensional diffeomorphisms with two saddle-foci (see also \citep{DP23}). For two-dimensional endomorphisms under the partial hyperbolicity condition, the result of Theorem \ref{thm:persis_hdd} is Theorem B in \citep{BCP21}.
\end{rem}

It is well-known that every point in a transitive, uniformly hyperbolic set $\Lambda$ is a limit point of hyperbolic periodic points (with the same dimension of the unstable manifolds equal to $ind(\Lambda)$), and the $C^r$-closure of invariant manifolds of these periodic points contains the stable and unstable sets of $\Lambda$, see e.g. \cite[Theorem 6.4.15]{KH95}. Hence, whenever we have heterodimensional dynamics involving two such sets, we can always obtain, by an arbitrarily small perturbation, a heterodimensional cycle associated with some hyperbolic periodic orbits. Thus, Theorem \ref{thm:persis_hdd} implies

\begin{cora}\label{cor:persis_hdd}
Any dynamical system of class $C^r$ ($r=1,\ldots,\infty,\omega$) having heterodimensional dynamics of coindex 1 can be $C^r$-approximated by a system which has $C^1$-robust heterodimensional dynamics.
\end{cora}

We stress that the results hold true, in particular, in the real-analytic case ($r=\omega$): given a real-analytic dynamical system on a real-analytic manifold we consider any complex neighborhood $\mathcal{M}$ of this manifold such that the system admits a holomorphic extension on $\mathcal{M}$; then the $C^\omega$-topology in Theorem \ref{thm:persis_hdd} and Corollary \ref{cor:persis_hdd} is the topology of uniform convergence on compacta in $\mathcal{M}$.

In fact, we obtain Theorem \ref{thm:persis_hdd} from its ``constructive version'', Theorem \ref{thm:persis_hdd_2p} below. Namely, to obtain the result of Theorem \ref{thm:persis_hdd}, one proceeds as follows. First, bring a given heterodimensional cycle into a general position (so that it satisfies the non-degeneracy conditions defined in Sections \ref{sec:gcon} and \ref{sec:gc4}) --  this can be done by an arbitrarily small $C^r$-perturbation of any heterodimensional cycle. Then, we embed our system $f$ into a finite-parameter family of perturbations $f_\varepsilon$ with {\em at least 2 parameters}. We formulate certain explicit conditions in Section \ref{sec:gfamily} which define an open and dense set in the space of $C^r$-families $f_\varepsilon$ such that $f_0=f$. An arbitrary family from this set is called a {\em proper unfolding} of $f$.

\begin{thma}\label{thm:persis_hdd_2p}
Let $f$ be of class $C^r$ $(r=2,\dots,\infty,\omega)$ and have a non-degenerate heterodimensional cycle of coindex 1, and let $f_\varepsilon$ be a proper unfolding of $f$. Then, arbitrarily close to $\varepsilon=0$ in the space of parameters there exist open regions where $f_\varepsilon$ has $C^1$-robust heterodimensional dynamics.
\end{thma}

\begin{rem}
As the family $f_\varepsilon$ is of class $C^r$, small $\varepsilon$ correspond to $C^r$-small perturbations of $f$. Thus, this theorem implies Theorem \ref{thm:persis_hdd} indeed. It should be stressed, however, that a result similar to Theorem \ref{thm:persis_hdd_2p} does not necessarily hold for one-parameter families, see Theorem \ref{thm:rationaltheta}.
\end{rem}

\begin{rem}
The proof of Theorem \ref{thm:persis_hdd_2p} (and all remaining results in this paper) makes use of the existence of certain $C^2$ coordinate transformations (local partial linearization, see Section \ref{sec:gc4}). We do not know, therefore, if these results hold for $C^1$ systems (except for Theorem \ref{thm:stanew}, whose $C^1$-analogue is proven in \cite{BDK12}). Thus, when dealing with systems of class $C^1$, in order to derive Theorem \ref{thm:persis_hdd} and Corollary \ref{cor:persis_hdd} from Theorem \ref{thm:persis_hdd_2p}, one can first perturb the system to make it $C^2$ and only then apply Theorem \ref{thm:persis_hdd_2p}.
\end{rem}

The non-degeneracy/propriety conditions of Theorem \ref{thm:persis_hdd_2p} are explicit, and checking them requires only a finite amount of computations with a finite number of periodic and heteroclinic orbits. Thus, Theorem \ref{thm:persis_hdd_2p} provides a universal tool for detecting and demonstrating the robust coindex-1 heterodimensional dynamics in multi-dimensional systems. This theorem can be used for showing that robust heterodimensional dynamics exist in specific restrictive settings, for example for polynomial perturbations, perturbations which keep various sorts of symmetry, etc.., and can be directly applied to dynamical systems coming from scientific applications, which usually have a form of finite-parameter families of differential equations or maps.  

\subsection{Stabilization of heterodimensional cycles}

The hyperbolic basic sets involved in the robust heterodimensional dynamics described in Theorem \ref{thm:persis_hdd_2p} are not necessarily homoclinically related\footnote{Two hyperbolic basic sets of the same index are homoclinically related if their stable and unstable manifolds have transverse intersections.} to the continuations of the periodic orbits from the original heterodimensional cycle. This means that even though parameter values corresponding to the existence of heterodimensional cycles are dense in the open regions of robust heterodimensional dynamics given by Theorem \ref{thm:persis_hdd_2p}, it may happen that these heterodimensional cycles do not contain the continuations of the periodic orbits of the original cycle.

To address this question, we employ the following adaptation of a definition from \citep{BD12,BDK12}. Recall that a heterodimensional cycle in our terminology refers to an invariant set consisting of only four orbits (two periodic and two connecting ones), see Definition \ref{defi:hdc}.

\begin{defi}[Stabilization of heterodimensional cycles]\label{defi:stabilization}
Let a dynamical system $f$ of class $C^r$ have a  heterodimensional cycle involving two  hyperbolic periodic orbits $L_1$ and $L_2$. The cycle is called $C^r$-{\em stabilizable} if arbitrarily close, in $C^r$, to $f$ there exists a dynamical system $g$, which exhibits $C^1$-robust heterodimensional dynamics involving non-trivial hyperbolic basic sets $\Lambda_1$ and $\Lambda_2$ that contain the continuations of  $L_1$ and $L_2$, respectively.
\end{defi}

In particular, since the stable and unstable manifolds of $L_j$ are dense in the stable and, respectively, unstable sets of the basic set $\Lambda_j$, $j=1,2$, it follows that systems with heterodimensional cycles involving the continuations of $L_1$ and $L_2$ accumulate on $g$ in the $C^r$-topology (and, more generally, on any $C^r$-system which is sufficiently close to $g$ in $C^1$).

\begin{defi}[Local stabilization of heterodimensional cycles]
A heterodimensional cycle is  {\em locally $C^r$-stabilizable}, if, given any neighborhood $U$ of the cycle, the $C^r$-close to $f$ system $g$ of Definition \ref{defi:stabilization} can be chosen such that the sets the sets $\Lambda_{1,2}$ and the corresponding new heterodimensional cycles involving the continuations of $L_1$ and $L_2$ all lie in $U$.
\end{defi}

Bonatti and D\'iaz constructed in \citep{BD12} an example of diffeomorphisms with heterodimensional cycles which cannot be $C^1$-stabilized (hence they cannot be $C^r$-stabilized). This work motivated the paper \citep{BDK12} by Bonatti, D\'iaz and Kiriki, showing that all heterodimensional cycles except for the so-called twisted ones (this class contains the cycles from the example in \citep{BD12}) can be locally $C^1$-stabilized. In the above definitions, we require a higher regularity of the stabilizing perturbations. We solve the question of local $C^r$-stabilization in Theorem \ref{thm:stanew} below.

\begin{figure}[!h]
\begin{center}
\includegraphics[width=0.9\textwidth]{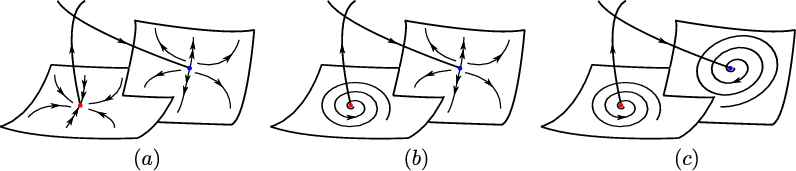}
\end{center}
\caption{Three cases of a heterodimensional cycle with two hyperbolic fixed points of a three-dimensional diffeomorphism. The central multipliers -- corresponding to the weakest contraction rate at the left fixed point and the weakest expansion rate at the right fixed point -- are both real in the saddle case (a), one real and one complex in the saddle-focus case (b), and both complex in the double-focus case (c).}
\label{fig:3cases}
\end{figure}

We distinguish three main cases of heterodimensional cycles, as depicted in Figure \ref{fig:3cases}: saddle, saddle-focus, and double-focus, depending on whether the central multipliers are real or not (see Section \ref{sec:local} for the precise definition). We show that in the saddle-focus and double-focus cases, the robust heterodimensional dynamics given by Theorem \ref{thm:persis_hdd_2p} are always associated with hyperbolic basic sets which are homoclinically related to the continuations of the periodic orbits in the original heterodimensional cycle, see Theorem \ref{thm:sad-foc_main}. However, in the saddle case, whether this homoclinic relation holds or not, this depends on whether the heterodimensional cycle is of type I or type II, as described in Section \ref{sec:types}. Thus, we have

\begin{thma}\label{thm:stanew}
Given any $r=2,\dots,\infty,\omega$, a non-degenerate heterodimensional cycle of coindex 1 in the saddle-focus and double-focus cases is locally $C^r$-stabilizable. In the saddle case, the cycle is locally $C^r$-stabilizable if and only if it is not of type I.
\end{thma}

Up to technical details, the type-I cycles correspond to twisted cycles from \citep{BDK12}. So, this theorem is the high regularity counterpart of the main result of \citep{BDK12}. Similarly to \citep{BDK12}, the type-I cycles become $C^r$-stabilizable (though not locally) when at least one of the periodic orbits in the cycle has a transverse homoclinic, see Corollary \ref{cor:stabilization_lambda}.

\subsection{Heterodimensionality vs. equidimensionality}

In the simplest setting, a {\em cycle} is an invariant set which consists of a cyclically ordered finite collection of periodic orbits and orbits that connect them (heteroclinic orbits), such that for each periodic orbit in the cycle 
there is exactly one heteroclinic orbit that tends to this periodic orbit in backward time and to the next periodic orbit in forward time. If all periodic orbits in the cycle have the same dimension of the unstable manifold and the same dimension of the stable manifold, the cycle is called {\em equidimensional}, and it is {\em heterodimensional} otherwise, see \cite{NP76}.

In this paper we consider the simplest case of heterodimensional cycles: coindex-1, and with 2 periodic orbits. The simplest case of an equidimensional cycle has just one periodic orbit and one homoclinic (the orbit of an 
intersection of the stable and unstable manifolds of the periodic orbit).

For a uniformly hyperbolic system, any cycle is equidimensional and the connecting orbits correspond to transverse intersections of the stable and unstable manifolds. Therefore, heterodimensional cycles and equidimensional cycles with non-transverse intersections ({\em homoclinic} or  {\em heteroclinic  tangencies}) are two obvious obstructions to the hyperbolicity.

A non-trivial fact is that both heterodimensional cycles and non-transverse equidimensional cycles can be robust. Thus, any hope that uniformly hyperbolic systems could be dense in the space of dynamical systems was destroyed with the above mentioned example by Abraham and Smale \cite{AS70} of a $C^1$-open region in the space of 4-dimensional diffeomorphisms where diffeomorphisms with heterodimensional cycles are dense and the example by Newhouse \citep{Ne70} of a $C^2$-open region in the space of 2-dimensional diffeomorphisms where diffeomorphisms with non-transverse equidimensional cycles are dense.

Newhouse built a theory of {\em thickness} of hyperbolic sets (for two-dimensional $C^2$-diffeomorphisms) and introduced a concept of a {\em wild} hyperbolic set: a non-trivial hyperbolic basic set whose stable and unstable sets have a tangency, for the system itself and for every $C^2$-close system. Essentially, if a hyperbolic set
is ``thick enough'' and its stable and unstable sets have a tangency, then the tangency is, typically, $C^2$-robust, i.e., the hyperbolic set is wild. Based on this theory, Newhouse 
proved in in \cite{Ne79} the $C^2$-{\em persistence of homoclinic tangencies}: for any generic one-parameter family of $C^r$ surface diffeomorphisms ($r\geq 2$) which unfolds a quadratic homoclinic tangency or an equidimensional cycle with a quadratic heteroclinic tangency,   
there exist open intervals of parameter values for which a wild
hyperbolic set exists and parameter values corresponding to quadratic homoclinic tangencies are dense in these intervals (a generalization to multi-dimensional systems was done in \citep{GTS93,PV94,Ro95}).  

Although many examples of robust heterodimensional cycles appeared after the Abraham-Smale construction, a heterodimensional analogue of the Newhouse theory was missing until the discovery of {\em blender} by Bonatti and D\'iaz \citep{BD96,BD08}. A blender is a hyperbolic basic set whose projection along strong-stable or strong-unstable directions contains an open set, see Section \ref{sec:blender} and Appendix for the precise definition. This ``openness in projection'' property make the heterodimensional dynamics that involve a blender robust. Thus, blenders play the same role in the creation of robust heterodimensional dynamics as Newhouse thick horseshoes do for persistent homoclinic tangencies. In particular, finding a blender near a perturbed heterodimensional cycle is the essence of the Bonatti-D\'iaz $C^1$-persistence result \cite{BD08} and our $C^r$-persistence results (Theorems \ref{thm:persis_hdd} and \ref{thm:persis_hdd_2p}).

Note that Newhouse ``thick horseshoe'' construction of the robust non-transverse intersections of stable and unstable manifolds is different from the Abraham-Smale construction and the later Bonatti-D\'iaz blender construction -- in particular, the homoclinic tangencies in the Newhouse construction are $C^2$-persistent, but not necessarily $C^1$-persistent\footnote{We have also learnt that Crovisier and Gourmelon obtained results on the $C^{1+\varepsilon}$-persistence of homoclinic tangencies, in a work under preparation.}. Indeed, Ures \citep{Ur95} discovered that for $C^1$ surface diffeomorphisms, the thickness of a horseshoe loses continuous dependence on the system, which is a crucial condition for the Newhouse construction; and later Moreira \citep{Mo11} built a more general theory and proved the non-existence of persistent homoclinic tangencies for $C^1$ surface diffeomorphisms. However, there are several examples of $C^1$-persistent homoclinic tangencies in higher dimensions, see \citep{Si72b,As08,BD12,LLST22}, which use blenders, or their variants, as supporting structures.

It is important to mention that there is a strong link between the homoclinic tangencies and the ``heterodimensionality''. The simplest manifestation of this is that bifurcations of homoclinic tangencies can lead to the birth of coexisting sinks and saddles \cite{GS73, Go83}. Moreover, Newhouse showed in \cite{Ne74,Ne79} that a wild hyperbolic set of a generic area-contracting surface diffeomorphism is in the closure of the set of sinks, i.e., the periodic orbits with different indices (here - saddles and sinks) are generically inseparable from each other. Without the contraction of areas, one has coexisting sets of sinks, saddles, and sources \cite{GTS97}, and, in the multi-dimensional case, coexisting saddles of different indices \cite{Ro95,GST08}. In \cite{LLST22}, we use the results of the present paper to give conditions under which the saddles of different indices that are born out of a homoclinic tangency get involved into the $C^1$-robust heterodimensional dynamics, as in Definition \ref{defi:hdd_robust}.

\subsection{Applications of Theorem \ref{thm:persis_hdd_2p}}

A commonly shared conjecture is that any $C^r$-diffeomorphism is either uniformly hyperbolic (i.e., every chain-recurrent class is uniformly hyperbolic) or it is arbitrarily close in $C^r$ to a diffeomorphism with wild hyperbolic sets (persistent homoclinic tangencies) or with robust heterodimensional dynamics, or both, see \cite{Pa00,Bo11} 
(for flows, one should also add a possibility of robust ``Lorenz-like'' dynamics \cite{ABS83,GW79,MPP04}). Irrespective of whether this conjecture is true or not, it is an empirical fact that homoclinic tangencies easily appear in many non-hyperbolic systems; we expect the same should be true for heterodimensional cycles.

In support of such claim, we have shown in a series of papers \citep{Li16,LT17,LT20,LLST22} that heterodimensional cycles emerge due to several types of homoclinic bifurcations. In fact, in the spirit of \cite{Tu96,Tu10,Go14}, one can conjecture that coindex-1 heterodimensional cycles can appear, with very few exceptions, in any homoclinic/heteroclinic bifurcation whose effective dimension allows it, i.e., when the dynamics of the map under consideration are not reduced to a two-dimensional invariant manifold and the map is not sectionally-dissipative (i.e., not area-contracting)\footnote{To have a heterodimensional cycle, we need saddles with different dimensions of unstable manifolds, and these conditions are obviously necessary. The conjecture is that they should also be sufficient for the birth of a heterodimensional cycle in most situations.}.

We believe this conjecture is true, so Theorem \ref{thm:persis_hdd_2p} allows for establishing the presence of robust heterodimensional dynamics whenever a non-hyperbolic chaotic behavior with more than one positive Lyapunov exponent (the ``hyperchaos'' in the terminology of \cite{Ro79}) is observed. In particular, it was also shown in \citep{Li16,LT17,LT20} that the coindex-1 heterodimensional cycles can be a part of a pseudohyperbolic chain-transitive attractor which appears in systems with Shilnikov loops \citep{Sh70,TS98,GKT21} or after a periodic perturbation of the Lorenz attractor \citep{TS08}. It, thus, follows from Theorems \ref{thm:persis_hdd} and \ref{thm:persis_hdd_2p} that the attractor in such systems remains heterodimensional for an open set of parameter values.

An important feature of robust heterodimensional dynamics is the robust presence of orbits with a zero Lyapunov exponent. In particular, the result of \cite{DG09} implies, for parametric families described by our Theorem \ref{thm:persis_hdd_2p}, the existence of open regions of parameter values where a generic system has an ergodic invariant measure with at least one zero Lyapunov exponent, i.e., the dynamics for such parameter values are manifestly non-hyperbolic.

For a dense set of parameter values from such regions the system has a non-hyperbolic periodic orbit. Bifurcations of such periodic orbits depend on the coefficients of the nonlinear terms of the Taylor expansion of the first-return map restricted to a center manifold. The degeneracy in the nonlinear terms increases complexity of the bifurcations. It follows from \cite{AST17,AST18} that once the so-called ``sign conditions'' are imposed on a heterodimensional cycle, the regions of robust heterodimensional dynamics given by Theorem \ref{thm:persis_hdd} contain a $C^r$-dense set of systems having infinitely degenerate (flat) non-hyperbolic periodic orbits. This fact also leads to the $C^\infty$-genericity (for systems from these regions) of a superexponential growth of the number of periodic orbits and the so-called $C^r$-universal dynamics, see \cite{AST17,AST18}.

\subsection{Blenders}\label{sec:blender}

As we mentioned, the main object responsible for the robustness of heterodimensional dynamics is a particular type of a hyperbolic set, a blender, introduced by Bonatti and D\'iaz in \cite{BD96}. It can be defined in many ways \citep{BDV,BCDW,BD12a,NP12, BKR14,BBD16}. For instance, the ``operational definition'' as in \cite[Definition 6.11]{BDV} can be rephrased as follows. Let $f$ be a dynamical system on a smooth manifold $\mathcal{M}$ with $\dim(\mathcal{M})\geqslant 3$ if $f$ is a diffeomorphism, or with $\dim(\mathcal{M})\geqslant 4$ if $f$ is a flow.

\begin{defi}[Blenders]\label{defi:blender_old}
A basic hyperbolic set $\Lambda$ of $f$ is called a center-stable (cs) blender if there exists a $C^1$-open set $\mathcal{D}^{ss}$ of  $d^{ss}$-dimensional discs (embedded copies of $\mathbb{R}^{d^{ss}}$) with $d^{ss}$ strictly less than the dimension of the stable manifolds of the orbits of $\Lambda)$, such that for every system $g$ which is $C^1$-close to $f$, for the hyperbolic continuation $\Lambda_g$ of the basic set $\Lambda$, the set $W^u(\Lambda_g)$ intersects every element from $\mathcal{D}^{ss}$; a center-unstable (cu) blender is a cs-blender for the dynamical system obtained from $f$ by the time-reversal.
\end{defi}

It is immediate from this definition that the existence of heterodimensional dynamics involving a blender is, essentially, a reformulation of the existence of the robust heterodimensional dynamics. Namely, for a system having two hyperbolic sets $\Lambda_{1,2}$ where $ind(\Lambda_2) > ind(\Lambda_1)$ and $\Lambda_1$ is a cs-blender, if there exists a transverse intersection of $W^u(\Lambda_2)$ with $W^s(\Lambda_1)$, and $W^s(\Lambda_2)$ contains a disc from $\mathcal{D}^{ss}$ defined as in the previous section (so there is an intersection of $W^u(\Lambda_1)$ with this piece of $W^s(\Lambda_2)$), then the system exhibits $C^1$-robust heterodimensional dynamics: the intersection of $W^u(\Lambda_2)$ with $W^s(\Lambda_1)$ survives small perturbations of the system because of the transversality, and the non-transverse intersection of $W^u(\Lambda_1)$ with $W^s(\Lambda_2)$ survives simply by the definition of the blender, as $W^s(\Lambda_2)$ varies continuously as the system varies.  

It should be noted that the blenders obtained in this paper have additional dynamical properties (partial hyperbolicity, the existence of a special Markov partition) which are not included in Definition \ref{defi:blender_old}. This additional structure is important for the actual construction of blenders and it is also crutial in applications, for example for proving the $C^1$-persistence of a certain type of homoclinic tangencies \cite{BD12a, LLST22}. We describe such blenders in
Definition \ref{defi:blender}, and call them {\em standard}. A standard blender is a version of the {\em blender-horseshoe} defined in \cite{BD12a}. The main difference  is that the latter has a Markov partition of exactly two elements, whereas we take the Markov partition with a large number of elements, like in \citep[Section 6.2]{BDV}. Since Definition \ref{defi:blender_old} suffices for discussing the persistence problem of heterodimensional cycles, the main goal of this paper, we do not make a digress to define standard blenders here. Instead, we detail the construction in the Appendix, see Definition \ref{defi:blender}. In Proposition  \ref{prop:stanble}, we prove that the blenders we construct in this paper are indeed standard blenders.

The central result (Theorem \ref{thm:blender} below) of the current paper is that we identify a class of heterodimensional cycles for which standard blenders exist in an arbitrarily small neighborhood of the cycle. In other words, any such cycle is {\em a limit of an infinite sequence of standard blenders}. We infer Theorem \ref{thm:persis_hdd_2p} from this result by showing that a proper unfolding of a non-degenerate heterodimensional cycle creates heterodimensional cycles of the ``blender-producing'' class described in Theorem \ref{thm:blender}, thus proving creation of blenders by a generic perturbation of an arbitrary heterodimensional cycle of coindex 1.

We always enumerate the periodic orbits $L_{1,2}$ in the heterodimensional cycle such that $ind(L_1)<ind(L_2)$, so the intersection $W^u(L_1)$ with $W^s(L_2)$ is fragile, whereas the intersection $W^u(L_2)\cap W^s(L_1)$ is transverse. By a multiplier of a periodic orbit, we mean an  eigenvalue of the derivative matrix of the return map at the corresponding fixed point, see Section \ref{sec:local}. The central multipliers $\lambda_{1,1}$ and $\gamma_{2,1}$ are the nearest to the unit circle among those multipliers of $L_1$ whose absolute value is smaller than $1$ and, respectively, the
nearest to the unit circle among those multipliers of $L_2$ whose absolute value is greater than $1$. Recall that we distinguish three cases of non-degenerate heterodimensional cycles of coindex 1: saddle, saddle-focus, and double-focus, depending on whether $\lambda_{1,1}$ and $\gamma_{2,1}$ are real or not (in the saddle-focus case, we assume that $\lambda_{1,1}$ is complex and $\gamma_{1,1}$ is real, with no loss of generality). In the saddle case (both $\lambda_{1,1}$ and $\gamma_{2,1}$ are real) we also have cycles of type I and type II (see Section \ref{sec:types}). It is well-known \cite{NPT83,St82} that the values of
$$\theta= - \frac{\ln |\lambda_{1,1}|}{\ln|\gamma_{2,1}|}$$
and of
$$\omega_1=Arg(\lambda_{1,1}), \qquad \omega_2=Arg(\gamma_{2,1})$$
(when $\lambda_{1,1}$ or $\gamma_{2,1}$ are complex) are {\em moduli}, i.e., invariants of topological equivalence, for systems with heterodimensional cycle.

\begin{thma}\label{thm:blender}
Let $f$ be of class $C^r$ $(r=2,\dots,\infty,\omega)$ and have a non-degenerate heterodimensional cycle of coindex-1, and let $U$ be any neighborhood of the cycle.
\begin{itemize}
\item In the saddle case, if the cycle is of type I and $\theta$ is irrational, then a standard blender exists in $U$.
\item In the saddle-focus case when $\theta$, $\frac{1}{2\pi}\omega_1$ and $1$ are rationally independent, and in the double-focus case when $\theta$, $\frac{1}{2\pi}\omega_1$, $ \frac{1}{2\pi}\theta\omega_2$ and $1$ are rationally independent, the system has simultaneously a standard cs-blender and a standard cu-blender in $U$. The blenders have different indices and their stable and unstable sets intersect $C^1$-robustly.\end{itemize}
\end{thma}

The result is the summary of Theorems \ref{thm:blender_in_type1} and \ref{thm:sad-foc_main} of Section \ref{sec:intro2}. It should be stressed that Theorem  \ref{thm:blender} is a non-perturbative result; this is the major difference with other works, where heterodimensional cycles are unfolded to obtain blenders, see  \citep{BD08,DP19,BCP21}. Note that by Theorem \ref{thm:rationaltheta} of Section \ref{sec:intro2}, no blenders exist in a sufficiently small neighborhood of the hetrodimensional cycle in the saddle case when $\theta$ is rational; in the case of complex central multipliers a similar result can also be derived \cite{Lipre}. Therefore, we conclude that when the values of the moduli change, {\em new blenders are ceaselessly produced by the heterodimensional cycle}. One can see here a parallel to Gonchenko's theory of a homoclinic tangency, which relates dynamics near a homoclinic tangency - the structure of hyperbolic sets, the existence of infinitely many sinks - to arithmetic properties of moduli of topological and $\Omega$-conjugacy \cite{GS86,GS90}.

The technique we use to establish the blender is based on the approximation of the first-return map near a heterodimensional cycle by an iterated function system (IFS) which is composed of a collection of affine maps of an interval, see formula (\ref{eq:blender:1}). This is similar to the approach used in many other works, see e.g. \citep{BDV,NP12, BD12a,BKR14,BBD16}. Note also that since the maps in our IFS are nearly affine, we can expect that the parablenders, introduced by Berger, can be implemented in our case too, cf. \citep{Be16,BCP15,BCP21}. \\ ~\\

Throughout the rest of the paper, by heterodimensional cycles/dynamics we always mean those of coindex 1, and the blenders we find/construct are always standard blenders as defined in the Appendix (see Proposition \ref{prop:stanble}). In Section \ref{sec:intro2}, we give precise definitions of the notions used in the paper, introduce non-degeneracy conditions for heterodimensional cycles, define the proper unfolding families, and give a complete formulation of the results. We start the proofs for the saddle case in Section \ref{sec:firstreturn}, where a computation for the first-return maps is carried out, and we prove the partial hyperbolicity of these maps. After that,  we find blenders near type-I cycles in Section \ref{sec:blendersintype1}. Next, in Section \ref{sec:sta_of_cycles}, we investigate the bifurcations of heterodymensional cycles in the saddle case and construct robust heterodimensional dynamics using the previously obtained blenders. Finally, we deal with the saddle-focus and double-focus cases in Section \ref{sec:complex_ev}.

\section{Robust heterodimensional dynamics in finite-parameter families}\label{sec:intro2}

In this section, we give a precise formulation of the results, which, in particular, imply Theorems \ref{thm:persis_hdd_2p}, \ref{thm:stanew} and \ref{thm:blender}. We consider the discrete and continuous-time cases. For both cases we define local maps and transition maps near the heterodimensional cycle (see Sections \ref{sec:local} and \ref{sec:gcon}). After that the proofs are solely based on the analysis of these maps and hence hold for both cases simultaneously.

We start with a more precise description of a heterodimensional cycle.
Let $f$ be a $C^r$-diffeomorphism of a $d$-dimensional manifold or a $C^r$-flow
of a $(d+1)$-dimensional manifold, where $d\geq 3$ and $r=2,\dots,\infty,\omega$. Let $f$ have a heterodimensional cycle $\Gamma$ of coindex $1$ associated with two hyperbolic periodic orbits $L_1$ and $L_2$ with
${\rm dim}\, W^u(L_2)={\rm dim}\, W^u(L_1)+1$. Along with the orbits $L_1$ and $L_2$, the heterodimensional cycle
$\Gamma$ 
consists of two heteroclinic orbits
$\Gamma^0\in W^u(L_1)\cap W^s(L_2)$ and $\Gamma^1\in W^u(L_2)\cap W^s(L_1)$.
Due to the difference in the dimensions of $W^u(L_1)$ and $W^u(L_2)$, the intersection $W^u(L_1)\cap W^s(L_2)$ is non-transverse and can be removed by 
a small perturbation. We call the orbit $\Gamma^0$ a {\em fragile heteroclinic orbit}. On the other hand, the intersection $W^s(L_1)\cap W^u(L_2)$ at the points of the orbit $\Gamma^1$ is assumed to be transverse and it gives a smooth one-parameter family of heteroclinic orbits. We call them {\em robust heteroclinic orbits}. See Figure \ref{fig:hdc0} for an illustration.

The two orbits $\Gamma^0$ and $\Gamma^1$ will be required to satisfy certain non-degeneracy conditions, introduced in Sections \ref{sec:gcon} and \ref{sec:gc4}. Our goal is to show how $C^1$-robust heterodimensional dynamics emerge in a small neighborhood of the cycle 
$\Gamma$ after $C^r$-small perturbations. The mechanisms for that depend on the type of the heterodimensional cycle, as described in detail below.

\subsection{Local maps near periodic orbits}\label{sec:local}
\begin{figure}[!h]
\begin{center}
\includegraphics[width=0.8\textwidth]{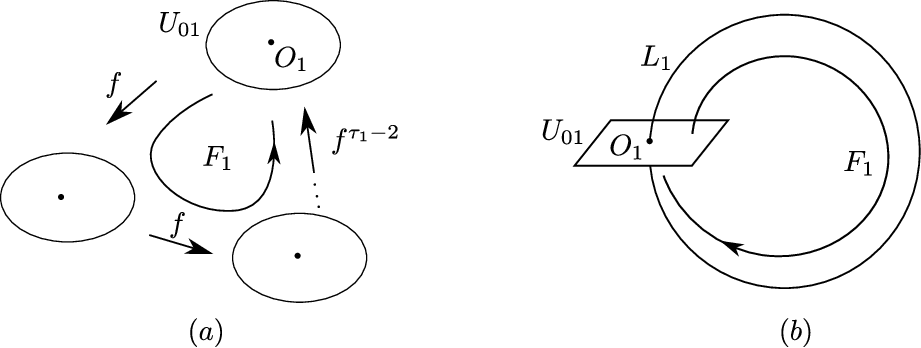}
\end{center}
\caption{The local maps near $O_1$ in the discrete-time case (a) and the continuous-time case (b).}
\label{fig:localmaps}
\end{figure}

In the discrete-time case, $f$ is a diffeomorphism. Let $O_1$ and $O_2$ be some points of the orbits $L_1$ and $L_2$. We take a small neighborhood $U_{0j}$ of the point $O_j$, $j=1,2$, and consider the first-return map $F_{j}$ in this 
neighborhood: $F_{j}=f^{\tau_j}$ where $\tau_j$ is the period of $O_j$ (see Figure \ref{fig:localmaps}). 

In the continuous-time case, the system $f$ is a flow generated by some smooth vector field. In this case we take some points $O_1\in L_1$ and $O_2\in L_2$ and let $U_{0j}$ ($j=1,2$) be small $d$-dimensional (i.e., of codimension 1) 
cross-sections transverse at $O_j$ to the vector field of $f$. Let $F_{j}$ be the first-return map (the Poincar\'e map) to the cross-section $U_{0j}$ (see Figure \ref{fig:localmaps}). 

In both cases $O_j$ is a hyperbolic fixed point of $F_{j}$: $F_{j} (O_j)= O_j$. The {\em multipliers} of $O_j$ are defined as the eigenvalues of the derivative of $F_{j}$ at $O_j$. The hyperbolicity means that no multipliers are equal to 1 in the absolute value; we assume that $d_j<d$ multipliers of $O_j$ lie outside the unit circle and $(d-d_j)$ multipliers lie inside. By our ``coindex-1 assumption'' ${\rm dim}\, W^u(L_2)={\rm dim}\, W^u(L_1)+1$, we have
\begin{equation*}
d_2=d_1+1.
\end{equation*}

We denote the multipliers of $O_j$, $j=1,2$, as $\lambda_{j,d-d_j}, \dots, \lambda_{j,1}, \gamma_{j,1}, \dots, \gamma_{j,d_j}$ and order them as follows:
\begin{equation}\label{eq:coindex1}
|\lambda_{j,{d-d_j}}|\leq \dots \leq |\lambda_{j,2}|\leq|\lambda_{j,1}|
<1<
|\gamma_{j,1}|\leq|\gamma_{j,2}|\leq\dots \leq |\gamma_{j,d_j}|.
\end{equation}
We call the largest in the absolute value multipliers inside the unit circle {\em center-stable} multipliers and those nearest to the unit
circle from the outside are called the {\em center-unstable} multipliers. The rest of the multipliers $\lambda$ and 
$\gamma$ are called strong-stable and, respectively, strong-unstable.

It is important whether the center-stable multipliers of $O_1$ and center-unstable multipliers
of $O_2$ are real or complex\footnote{
It is well-known that the dynamics are quite different in these cases: for real multipliers, one has partially hyperbolic dynamics with one-dimensional center typically, while nonreal multipliers mean higher-dimensional dynamics, e.g. can lead to the emergence of homoclinic tangencies, see e.g. \citep{DR92}.}. Note that by adding an arbitrarily small perturbation, if necessary, we can always bring the multipliers into a general position. In our situation,
this means that we can assume that $O_1$ has only one center-stable multiplier $\lambda_{1,1}$ which is real and simple, or a pair of simple complex conjugate center-stable multipliers
$\lambda_{1,1}=\lambda_{1,2}^*$; we also can assume that $O_2$ has either only one center-unstable multiplier 
$\gamma_{2,1}$ which is real and simple, or a pair of simple complex conjugate center-unstable multipliers
$\gamma_{2,1}=\gamma_{2,2}^*$.

Accordingly, we distinguish three main cases.
\begin{itemize}
\item Saddle case: here $\lambda_{1,1}$ and $\gamma_{2,1}$ are real and simple, i.e., we have 
$|\lambda_{1,2}|<|\lambda_{1,1}|$ and 
$|\gamma_{2,1}|<|\gamma_{2,2} |$.
\item Saddle-focus case: here either 
$$\lambda_{1,1}=\lambda_{1,2}^*=\lambda e^{i\omega}, \;\omega\in(0,\pi), \quad\mbox{and}\quad \gamma_{2,1}\quad \mbox{is real},$$
where $\lambda>|\lambda_{1,3}|$ and $|\gamma_{2,1}|<|\gamma_{2,2} |$, or
$$\gamma_{2,1}=\gamma_{2,2}^*=\gamma e^{i\omega},\;\omega\in(0,\pi), \quad\mbox{and}\quad \lambda_{1,1}\quad \mbox{is real},$$
where $0<\gamma<|\gamma_{2,3}|$ and $|\lambda_{1,2}|<|\lambda_{1,1}|$. Note that the second option is reduced to
the first one by the inversion of time. Therefore, we assume below that in the saddle-focus case $\lambda_{1,1}$ is complex and $\gamma_{2,1}$ is real.
\item Double-focus case: here
$$\lambda_{1,1}=\lambda_{1,2}^*=\lambda e^{i\omega_1},\;\omega_1\in(0,\pi),\quad\mbox{and}\quad 
\gamma_{2,1}=\gamma_{2,2}^*=\gamma e^{ i\omega_2},\;\omega_2\in(0,\pi),$$
where $\lambda>|\lambda_{1,3}|$, $0<\gamma<|\gamma_{2,3}|$.
\end{itemize}

Below, we denote $\lambda_{1,1}$ and $\gamma_{2,1}$ by $\lambda$ and $\gamma$ if they are real.
Let $d_{cs}$ denote the number of the center-stable multipliers of $O_1$ and $d_{cu}$ be the number of the center-unstable multipliers of $O_2$. We have $d_{cs}=d_{cu}=1$ in the saddle case, $d_{cs}=2, d_{cu}=1$ in the saddle-focus 
case, and $d_{cs}=d_{cu}=2$ in the double-focus case.

\begin{figure}[!h]
\begin{center}
\includegraphics[width=0.8\textwidth]{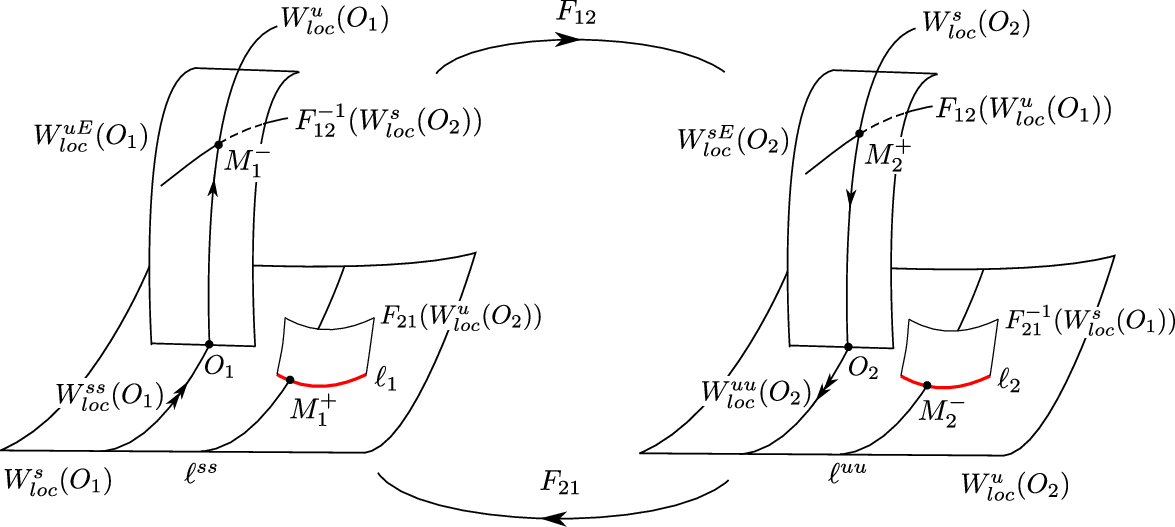}
\end{center}
\caption{A heterodimensional cycle of coindex $1$ satisfying conditions C1 - C3, which consists of two periodic orbits containing $O_1$ and $O_2$, a fragile heteroclinic orbit containing $M^-_1$ and $M^+_2$, and a robust heteroclinic orbit containing $M^-_2$ and $M^+_1$. Here $\ell^{ss}$ is a strong-stable leaf through $M^+_1$, $\ell^{uu}$ is a strong-unstable leaf through $M^-_2$, $\ell_1=W^s_{loc}(O_1)\cap F_{21}(W^u_{loc}(O_2))$ and $\ell_2=F^{-1}_{21}\ell_1$.}
\label{fig:hdc1}
\end{figure}

Recall (see e.g. \cite{sstc1,GST08}) 
that the first-return map $F_{1}$ near the point $O_1$ has a $d_1$-dimensional local unstable manifold 
$W^u_{loc}(O_1)$ which is tangent at $O_1$ to the eigenspace corresponding to the multipliers 
$\gamma_{1,1},\dots,\gamma_{1,d_1}$ and $(d-d_1)$-dimensional local stable manifold 
$W^s_{loc}(O_1)$ which is tangent at $O_1$ to the eigenspace corresponding to the multipliers 
$\lambda_{1,1},\dots,\lambda_{1,d-d_1}$. In $W^s_{loc}(O_1)$ there is a $(d-d_1-d_{cs})$-dimensional strong-stable $C^r$-smooth invariant manifold $W^{ss}_{loc}(O_1)$ which is tangent at $O_1$ to the eigenspace corresponding to the strong-stable multipliers $\{\lambda_{1,d-d_1}, \dots, \lambda_{1,d-d_1-d_{cs}}\}$. This manifold is a leaf of the strong-stable $C^r$-smooth foliation $\mathcal{F}^{ss}$ of $W^s_{loc}(O_1)$. There also exists a $(d_1+d_{cs})$-dimensional extended-unstable invariant manifold $W^{uE}_{loc}(O_1)$ corresponding to the center-stable multipliers and the multipliers $\gamma_{1,1},\dots ,\gamma_{1,d_1}$. Such manifold is not unique but all of them contain 
$W^u_{loc}(O_1)$ and are tangent to each other at the points of $W^u_{loc}(O_1)$ (see Figure \ref{fig:hdc1}). 

Similarly, the first-return map $F_{2}$ near the point $O_2$ has a $d_2$-dimensional local unstable manifold 
$W^u_{loc}(O_2)$ and $(d-d_2)$-dimensional local stable manifold 
$W^s_{loc}(O_2)$. In $W^u_{loc}(O_2)$ there is a $(d_2-d_{cu})$-dimensional strong-unstable invariant manifold $W^{uu}_{loc}(O_2)$ which is tangent at $O_1$ to the eigenspace corresponding to the strong-unstable multipliers 
$\{\gamma_{2,d_{cu}+1}, \dots, \gamma_{2,d_2}\}$. This manifold is a leaf of the strong-unstable $C^r$-smooth foliation $\mathcal{F}^{uu}$ of $W^u_{loc}(O_2)$. There also exists a $(d-d_2+d_{cu})$-dimensional extended-stable invariant manifold $W^{sE}_{loc}(O_2)$ corresponding to the center-unstable multipliers and the multipliers 
$\lambda_{2,1},\dots ,\lambda_{2,d-d_2}$. Any two such manifolds contain $W^s_{loc}(O_2)$ and are tangent to each other at the points of $W^s_{loc}(O_2)$. 

\subsection{Transition maps and geometric non-degeneracy conditions}\label{sec:gcon}
For each of the heteroclinic orbits $\Gamma^0$ and $\Gamma^1$, a transition map between neighborhoods of $O_1$
and $O_2$ is defined, as follows.

Consider, first, the case of discrete time, i.e., let $f$, be a diffeomorphism. 
Take four points $M^+_1\in\Gamma^1\cap W^s_{loc}(O_1)$, 
$M^-_1\in\Gamma^0\cap W^u_{loc}(O_1)$, $M^+_2\in\Gamma^0\cap W^s_{loc}(O_2)$, and 
$M^-_2\in\Gamma^1\cap W^u_{loc}(O_2)$. Note that $M^-_2$ and $M^+_1$ belong to the same robust heteroclinic orbit $\Gamma^1$ and $M^-_1$ and $M^+_2$ belong to the same fragile heteroclinic orbit $\Gamma^0$. Thus,  there exist positive integers $n_1$ and $n_2$ such that 
$f^{n_{1}}(M^-_1)=M^+_2$ and $f^{n_{2}}(M^-_2)=M^+_1$. We define the {\em transition maps} from a small neighborhood of $M^-_1$ to a small neighborhood of $M^+_2$ and from a small neighborhood of 
$M^-_2$ to a small neighborhood of $M^+_1$ as, respectively, $F_{12}= f^{n_1}$ and $F_{21}= f^{n_2}$ (see Figure \ref{fig:hdc1}).

In the continuous-time case, when $f$ is a flow, we take the points $M^+_1\in\Gamma^1\cap W^s_{loc}(O_1)$
and $M^-_1\in\Gamma^0\cap W^u_{loc}(O_1)$ on the cross-section $U_{01}$ and the points $M^+_2\in\Gamma^0\cap W^s_{loc}(O_2)$ and 
$M^-_2\in\Gamma^1\cap W^u_{loc}(O_2)$ on the cross-section $U_{02}$. Then
the transition map $F_{12}$ is defined as the map 
by the orbits of the flow which start in the cross-section $U_{01}$ near $M^-_1$ and hit the cross-section $U_{02}$ near
the point $M^+_2$, and the transition map $F_{21}$ is defined as the map 
by the orbits of the flow which start in the cross-section $U_{02}$ near $M^-_2$ and hit the cross-section $U_{01}$ near
the point $M^+_1$. 
By the definition, $F_{12}(M^-_1)=M^+_2$ and $F_{21}(M^-_2)=M^+_1$.

We can now precisely describe the non-degeneracy conditions which we impose on the heterodimensional cycle $\Gamma$.

\noindent\textbf{C1.} Simplicity of the fragile heteroclinic:
$F^{-1}_{12}(W_{loc}^{sE}(O_2))$ intersects $ W^{u}_{loc}(O_1)$ transversely at the point $M^-_1$, and
$F_{12}(W_{loc}^{uE}(O_1))$ intersects $ W^{s}_{loc}(O_2)$ transversely at $M^+_2$;\\[5pt]
\noindent\textbf{C2.} Simplicity of the robust heteroclinic: the leaf of 
$\mathcal{F}^{uu}$ at the point $M^-_2$ is not tangent to
$F_{21}^{-1}(W^s_{loc}(O_1))$ and the leaf of $\mathcal{F}^{ss}$ at the point $M^+_1$ is not tangent to
$F_{21}(W^u_{loc}(O_2))$; and \\[5pt]
\noindent\textbf{C3.} $\Gamma^1\cap (W^{ss}(O_1)\cup W^{uu}(O_2))=\emptyset$, i.e., $M^+_1\notin W^{ss}(O_1)$ and $M^-_2\notin W^{uu}(O_2)$.\\

Figure \ref{fig:hdc1} provides an illustration of these conditions. Note that condition C1 does not depend on the choice of $W_{loc}^{uE}(O_1)$ and $W_{loc}^{sE}(O_2)$, as any two extended stable/unstable manifolds are tangent to each other at the points of the stable or, respectively, unstable manifold, see the discussion in the end of Section \ref{sec:local}.  Moreover, the corresponding requirement in C1 is automatically satisfied if $O_2$ has no non-center stable multipliers or $O_1$ has no non-center unstable multipliers, that is, in the case where $O_2$ has a pair of complex conjugate center-unstable multipliers and $\dim W^u(O_2)=d_2=2$, or $O_1$ has a pair of complex conjugate center-stable multipliers and $W^s(O_1)=d-d_1=2$. Similarly, the corresponding requirements of condition C2 hold automatically in these cases.

The manifolds involved in these conditions depend continuously (as $C^1$-manifolds) on $f$ in the $C^r$-topology.
This implies that conditions C2 and C3 are $C^r$-open, and condition C1 is $C^r$-open in the class of systems with the heterodimensional cycle. It is also quite standard that one can always achieve the fulfillment of C1 and C2 by adding an arbitrarily $C^r$-small perturbation to $f$ (in the smooth case one adds a local perturbation to $f$; in the analytic case one uses the scheme described in \cite{BT86,GTS07}). In the case where condition C3 is not fulfilled, we do not need to perturb the system: for the same system $f$ we can always find, close to $\Gamma^1$, another robust heteroclinic orbit which satisfies C3. To see this, note
that condition C2 ensures that the line $\ell_1=W^s_{loc}(O_1)\cap F_{21}(W^u_{loc}(O_2))$ (corresponding to robust heteroclinics, see Figure \ref{fig:hdc1}) is not tangent to 
$W^{ss}(O_1)\cup F_{21}(W^{uu}_{loc}(O_2))$, so we can always shift the position of the point $M^+_1$ on this line and hence the position of $M^-_2=F_{21}^{-1}(M^+_1)$.

\subsection{Local partial linearization and the fourth non-degeneracy condition}
\label{sec:gc4}
There is one last non-degeneracy condition, which is different for the saddle case and the other cases. To state it precisely, let us introduce $C^r$-coordinates $(x,y,z)\in \mathbb{R}^{d_{cs}}\times\mathbb{R}^{d_1}\times\mathbb{R}^{d-d_1-d_{cs}}$ in $U_{01}$ such that
the local stable and unstable manifolds get straightened near $O_1$:
$$W^s_{loc}(O_1)=\{y=0\},\quad W^u_{loc}(O_1)=\{x=0, z=0\},$$
and the extended-unstable manifold $W^{uE}_{loc}(O_1)$ is tangent to $\{z=0\}$ when $x=0, z=0$ (see Section \ref{sec:firstreturn}). Moreover, the leaves of the foliation $\mathcal{F}^{ss}$ are also straightened and are given by 
$\{x=const, y=0\}$. In particular, we have
$$W^{ss}_{loc}(O_1)=\{x=0, y=0\}.$$

We also introduce $C^r$-coordinates 
$(u,v,w)\in\mathbb{R}^{d_{cu}}\times\mathbb{R}^{d-d_2}\times\mathbb{R}^{d_2-d_{cu}}$ in $U_{02}$ such that the local stable and unstable manifolds are straightened near $O_2$:
$$W^s_{loc}(O_2)=\{u=0, w=0\},\quad W^u_{loc}(O_2)=\{v=0\},$$
the extended-sable manifold $W^{sE}_{loc}(O_2)$ is tangent to $\{w=0\}$ when $u=0,v=0$, and the leaves of the foliation $\mathcal{F}^{uu}$ are also straightened and given by $\{u=const, v=0\}$. We have
$$W^{uu}_{loc}(O_2)=\{u=0,v=0\}.$$

We restrict the choice of the coordinates by a further requirement (which can always be fulfilled, see e.g. \cite{GST08})
that the first-return maps $F_{1}$ and $F_{2}$ act linearly on center-stable and, respectively, center-unstable coordinates. Namely, if we restrict these maps on $W^s_{loc} (O_1)=\{y=0\}$ and, respectively, $W^u_{loc} (O_2)=\{v=0\}$ and use the notation
$F_{1}|_{W^s_{loc}(O_1)}:(x,z)\mapsto (\bar x, \bar{z})$ and 
$F_{2}|_{W^u_{loc}(O_2)}:(u,w)\mapsto(\bar{u},\bar{w})$, then we have
\begin{equation}\label{eq:intro:2.1}\!\!\!\!
\begin{array}{ll}
\,\;\mbox{Saddle case:}&\!\bar{x} = \lambda x \quad\quad\mathrm{and}\quad\quad \bar{u}=\gamma u;\\[14pt]
\begin{array}{c}\mbox{Saddle-focus}\\\mbox{case:}\end{array}&\!\bar{x} = \lambda\begin{pmatrix}
\cos \omega &\sin\omega \\
-\sin \omega & \cos\omega
\end{pmatrix}\!x \quad\mathrm{and}\quad \bar{u}=\gamma u;\\[15pt]
\begin{array}{c}\mbox{Double-focus}\\\mbox{case:}\end{array}&\!\!\bar{x} = \lambda\!\begin{pmatrix}
\!\cos \omega_1 &\!\!\sin\omega_1 \\
\!-\sin \omega_1 &\!\!\cos\omega_1
\end{pmatrix}\!x \;\;\mathrm{and}\;\; \bar{u}=\gamma\!\begin{pmatrix}
\!\cos \omega_2 &\!\!\sin\omega_2 \\
\!-\sin \omega_2 &\!\!\cos\omega_2
\end{pmatrix}\! u.\!\!\!\!\!\!\!\!\!\!\!
\end{array}
\end{equation}
We denote, in these coordinates, $M^+_1=(x^+,0,z^+)$ and $M^-_2=(u^-,0,w^-)$ (so condition C3 reads $x^+\neq 0$ and $u^-\neq 0$).

Recall that $F_{21}$ takes a point with coordinates $(u,v,w)$ to a point with coordinates $(x,y,z)$, where $x$ and $u$ are the center-stable and center-unstable coordinates near the points $O_1$ and $O_2$, respectively.
By condition C2, the line $\ell_{1}$ is not tangent to the foliation $\mathcal{F}^{ss}$
and the line $\ell_{2}= F^{-1}_{21}\ell_1$ is not tangent to the foliation 
$\mathcal{F}^{uu}$.
In the saddle case this means that these curves are transverse to these foliations (see Figure \ref{fig:hdc1}), so they are parametrized by the coordinates
$x$ (the line $\ell_1$) and $u$ (the line $\ell_2$). Therefore, as $F_{21}|_{\ell_2}$ acts as a diffeomorphism 
$\ell_2\to\ell_1$, we have
\begin{equation}\label{eq:intro:b}
b= \left.\frac{\partial x}{\partial u}\right|_{M^-_2}\neq 0.
\end{equation}

\noindent\textbf{C4.1} (saddle case). The quantity $\alpha:=b u^-/x^+$ satisfies
\begin{equation}\label{eq:intro:4}
|\alpha|\neq 1.
\end{equation}

Note that in the saddle case conditions C1 and C2 are equivalent (see \cite{Tu96}) to the requirement that the heteroclinic cycle $\Gamma$ is a partially hyperbolic set with the 1-dimensional central direction field which includes
the center-stable eigenvector at $O_1$ and the center-unstable eigenvector at $O_2$. As we show, $\alpha$ determines the behavior in the central direction: the first-return maps near $\Gamma$
are contracting in the central direction when $|\alpha|<1$ and expanding when $|\alpha|>1$ (see Lemma \ref{lem:conefields}). Note that $\alpha$ is an invariant of smooth coordinate transformations which keep
the foliations $\mathcal{F}^{ss}$ and $\mathcal{F}^{uu}$ locally straightened and keep the action of the local maps $F_{1}$ and $F_{2}$ in the central direction linear, as in \eqref{eq:intro:2.1}. Indeed, any such transformation
is linear in the central directions in a small neighborhoods $U_{01}$ and $U_{02}$ of the points $O_1$ and $O_2$, i.e., 
the coordinates $x$ and $u$ are only multiplied to some constants $c_x$ and $c_u$. As a result, the coefficient
$b$ is replaced by $b c_x/c_u$, and $x^+$ and $u^-$ are replaced by $c_x x^+$ and $c_u u^-$, so $\alpha$ remains unchanged. Similarly, the invariant $\alpha$ does not depend on the choice of the points $M^+_1$ and $M^-_2$ on the given heteroclinic orbit $\Gamma^1$. 

\begin{figure}[!h]
\begin{center}
\includegraphics[width=0.4\textwidth]{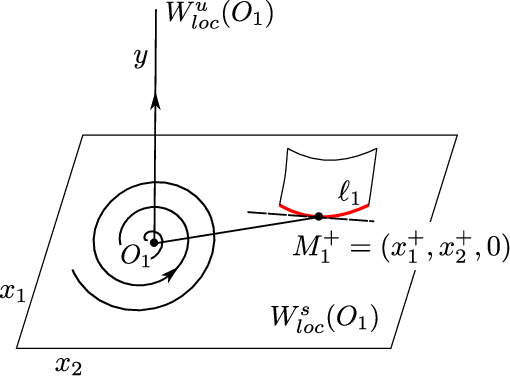}
\end{center}
\caption{An illustration of condition C4.2. The vector $(x^+_1,x^+_2)$ is not parallel to the $x$-vector component (the dashed line) of the tangent to $\ell_1$.}
\label{fig:C42}
\end{figure}

In the saddle-focus and double-focus cases, 
the partial hyperbolicity is not assumed, and no condition similar to C4.1 is needed. However, we need another condition:

\noindent\textbf{C4.2} (saddle-focus and double-focus cases). When the center-stable multipliers $\lambda_{1,1}$ 
and $\lambda_{1,2}$ are complex and $x\in \mathbb R^2$, the $x$-vector component of the tangent to $\ell_1$ at the point $M^1_+$ is not parallel to the vector $x^+$ (see Figure \ref{fig:C42}). When the center-unstable multipliers 
$\gamma_{2,1}$ 
and $\gamma_{2,2}$ are complex and $u\in \mathbb R^2$, the $u$-vector component of the tangent to $\ell_2$ at the point $M^-_2$ is not parallel to the vector $u^+$.

Like in condition C4.1, coordinate transformations that keep the action of the local maps $F_{1}$ and $F_{2}$ in $x$ and $u$ linear
are also linear in $x$ and $u$, respectively. This immediately implies that condition C4.2 is invariant with respect to the choice
of the linearizing coordinates. It also does not depend on the choice of the points $M^+_1$ and $M^-_2$. 

The heterodimensional cycles satisfying conditions C1-C4.1,2 will be further called {\em non-degenerate}.

\subsection{Finite-parameter unfoldings}\label{sec:gfamily}
The perturbations we use to prove Theorem \ref{thm:persis_hdd} are done within 
finite-parameter families $f_\varepsilon$ which we assume to be of class $C^r$ ($r=2,\dots, \infty,\omega$) jointly with respect to coordinates and parameters $\varepsilon$.

Let $f_0=f$; for any sufficiently small $\varepsilon$ the hyperbolic points $O_1$ and $O_2$ exist and depend smoothly on $\varepsilon$.
The corresponding multipliers also depend smoothly ($C^{r-1}$) on $\varepsilon$. We define
\begin{equation}\label{tht}
\theta(\varepsilon)=-\dfrac{\ln|\lambda|}{\ln|\gamma|}.
\end{equation}
In the saddle-focus and double-focus cases, an important role is also played by the frequencies $\omega(\varepsilon)$ and, respectively,
$\omega_{1,2}(\varepsilon)$. The values of $\theta$ as well as $\omega_{1,2}$ are moduli of topological conjugacy of diffeomorphisms with non-degenerate heterodimensional cycles (see \citep{NPT83,St82}).

The local stable and unstable manifolds of $O_{1,2}$, as well as 
their images by the transition maps $F_{12}$ and $F_{21}$, also depend smoothly on $\varepsilon$. The fragile heteroclinic $\Gamma^0$ is not, in general, preserved
when $\varepsilon$ changes. To determine whether the fragile heteroclinic disappears or not, one introduces a {\em splitting parameter} 
$\,\mu$, a continuous functional such that for any system $g$ from a small $C^r$-neighborhood of $f$ the absolute value of $\mu(g)$ equals to the distance between $W^s_{loc}(O_2)$ and $F_{12} (W^u_{loc}(O_1))$; the fragile heteroclinic persists for those $g$ for which $\mu(g)=0$. The codimension-1 manifold
$\mu=0$ separates the neighborhood of the system $f$ into two connected components; we define $\mu$ such that it changes sign
when going from one component to the other. 

A one-parameter family $f_\varepsilon$ is called a {\em generic one-parameter unfolding} of $f$ if $\mu(f_\varepsilon)$ depends on $\varepsilon$ smoothly and
$\frac{d\mu}{d\varepsilon}\neq 0$.
This means that we can make $\mu(f_\varepsilon)=\varepsilon$ by a smooth change of parameters.  

We also need to consider families depending on two or more parameters, i.e., $\varepsilon=(\varepsilon_1,\varepsilon_2,\dots)$. We call
the family $f_\varepsilon$ a {\em proper unfolding}, if $\frac{d\mu}{d\varepsilon}\neq 0$ (so the set 
$\mu(\varepsilon)=0$ forms a smooth
codimension-1 manifold ${\cal H}_0$ in the space of parameters $\varepsilon$)  and, the following conditions hold for the subfamily corresponding 
to $\varepsilon\in {\cal H}_0$:
\begin{itemize}
\item in the saddle case, $\frac{d\theta}{d\varepsilon}\neq 0$, where the
derivative is taken over $\varepsilon\in {\cal H}_0$ (this implies that we can make a smooth change of parameters in the family $f_\varepsilon$ such that 
$\mu(\varepsilon)=\varepsilon_1$ and $\theta(\varepsilon)=\varepsilon_2$);
\item in the saddle-focus case, the condition is that the functions 
$\theta(\varepsilon)$ and $\frac{1}{2\pi}\omega(\varepsilon)$ and $1$ are linearly independent in a neighborhood of 
$\varepsilon=0$ on ${\cal H}_0$;
\item in the double-focus case, the condition is the linear independence of $\theta(\varepsilon)$, $\frac{1}{2\pi}\omega_1(\varepsilon)$, $\frac{1}{2\pi}\omega_2(\varepsilon)\theta(\varepsilon)$ and $1$ in a neighborhood of 
$\varepsilon=0$ on ${\cal H}_0$\footnote{Note that if this condition holds, it also holds for the system obtained after a time reversal. Indeed, this operation interchanges $O_1$ and $O_2$, and hence the propriety condition for the resulting system becomes the linear independence of $\theta^{-1}(\varepsilon)$, $\frac{1}{2\pi}\omega_2(\varepsilon)$, $\frac{1}{2\pi}\omega_1(\varepsilon)\theta^{-1}(\varepsilon)$ and $1$, which is equivalent to the linear independence of $\theta(\varepsilon)$, $\frac{1}{2\pi}\omega_1(\varepsilon)$, $\frac{1}{2\pi}\omega_2(\varepsilon)\theta(\varepsilon)$ and $1$.}.
\end{itemize}
Note that the linear independence conditions for the saddle-focus and double-focus case are only used to ensure
that the corresponding quantities can be made {\em rationally independent} by an arbitrarily small change of
$\varepsilon$. However, we formulate the propriety conditions in this way in order to make the class of proper families open.

With the above definitions, the formulation of our main result, Theorem
\ref{thm:persis_hdd_2p} as given in Section \ref{sec:intro}, is now
complete. The proof goes differently in different cases: for the saddle case the theorem follows from the results described in Sections \ref{sec:types} and \ref{sec:msaddle},
and in the saddle-focus and double focus case it follows from the results of 
Section \ref{sec:mcom}.

\subsection{Three types of heterodimensional cycles in the saddle case}\label{sec:types}
In the saddle case, the proof of Theorem \ref{thm:persis_hdd_2p} is most involved: not because of technicalities, but because the dynamics emerging at
perturbations of the non-degenerate heterodimensional cycles depend, in the saddle case, very essentially on the type of the cycle.
According to that, we introduce three types of the heterodimensional cycles in the saddle case, as follows.

First, note (by counting dimensions) that in the saddle case condition C1 implies that the intersection 
of $F^{-1}_{12}(W_{loc}^{sE}(O_2))$
and $W^{uE}_{loc}(O_1)$ is a smooth curve, which we denote as $\ell^0$ (see Figure \ref{fig:hdc2}). At
 $\varepsilon=0$, this curve goes through the point $M^-_1$
and its image $F_{12}(\ell^0)$ goes through the point $M^+_2$. The tangent space 
$\mathcal T_{M^-_1}\ell^0$ lies in
$\mathcal T_{M^-_1}W^{uE}(O_1) = \{z=0\}$ and, by C1, 
$\mathcal T_{M^-_1}\ell^0\not\subset \mathcal T_{M^-_1}W_{loc}^{u}(O_1)=\{x=0, z=0\}$, which implies that
$\mathcal T_{M^-_1}\ell^0$ has a non-zero projection to the $x$-axis. Thus, the curve $\ell^0$ is parametrized by the coordinate $x$.
Similarly, the curve $F_{12}(\ell^0)$ is parametrized by coordinate $u$. The restriction of $F_{12}$ to $\ell^0$ is a diffeomorphism, so 
\begin{equation}\label{eq:a11}
a=\left.\frac{du}{dx}\right|_{M^-_1}\neq 0.
\end{equation}

\begin{figure}[!h]
\begin{center}
\includegraphics[width=0.8\textwidth]{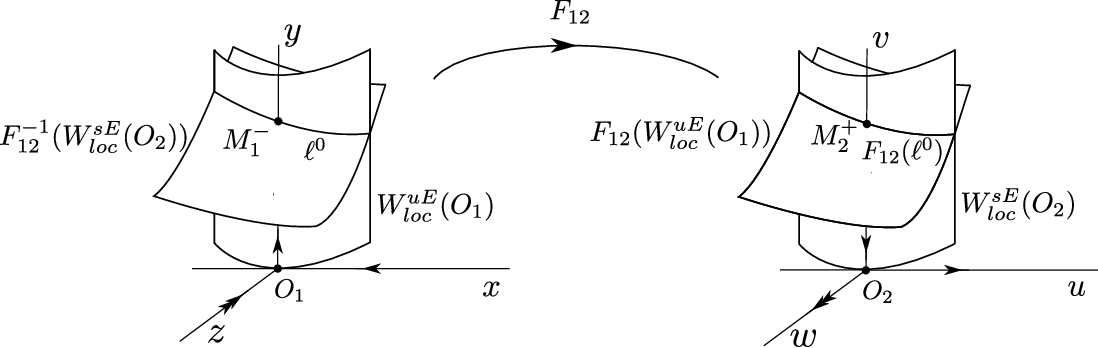}
\end{center}
\caption{Condition C1 in straightened coordinates.}
\label{fig:hdc2}
\end{figure}

We say that a heterodimensional cycle $\Gamma$ in the saddle case is of
\begin{itemize}
\item  
type I\footnote{Type-I and type-II cycles correspond to twisted, and, respectively, non-twisted cycles in \citep{BDK12}. Strictly speaking, our definition is somewhat more general, as the notion of twisted and untwisted cycles was introduced 
in \citep{BDK12} only in the case when the local maps and transition maps satisfy certain additional restrictions.}, if there exist points ~~ $M^+_1(x^+,0,z^+)\in\Gamma^1\cap U_{01}$ ~~ and ~~ $M^-_2(u^-,0,w^-)\in\Gamma^1\cap U_{02}$ ~~ such that $a x^+u^->0$;
\item 
type II, if there exist points ~~ $M^+_1(x^+,0,z^+)\in\Gamma^1\cap U_{01}$ ~~ and ~~ $M^-_2(u^-,0,w^-)\in\Gamma^1\cap U_{02}$ ~~ such that $a x^+u^-<0$;
\item 
type III, if there exist points $M^+_1\in\Gamma^1\cap U_{01}$ and $M^-_2\in\Gamma^1\cap U_{02}$ for which $a x^+u^->0$
and another pair of points $M^+_1\in\Gamma^1\cap U_{01}$ and $M^-_2\in\Gamma^1\cap U_{02}$ for which $a x^+u^-<0$.
\end{itemize}

The cycle of type III is,
by definition, a cycle which is simultaneously of type I and type II. Like in condition C4.1, one shows that the sign of $ax^+u^-$
is independent of the choice of coordinates which keep the action of the local maps $F_{1}$ in the neighborhood $U_{01}$ of $O_1$
and $F_{2}$ in the neighborhood $U_{02}$ of $O_2$ linear in the central coordinates $x$ and $u$. Thus, the above definition is
coordinate-independent. 

Notice that
$a$ is determined by a pair of points $M^-_1$ and $M^+_2$ on the fragile heteroclinic $\Gamma^0$, while $x^+$ and $u^-$ are coordinates of
points on the robust heteroclinic $\Gamma^1$. By \eqref{eq:intro:2.1} (the saddle case), if 
the central multipliers $\lambda$ and $\gamma$ are positive, the local maps $F_{1}$ and $F_{2}$ multiply $x^+$ and $u^-$
to positive factors, so the sign of $ax^+u^-$ is independent of the choice of the points $M^+_1$ and $M^-_2$ on 
$\Gamma^1$ in this case. Similarly,
it does not depend on the choice of the points $M^-_1$ and $M^+_2$ on $\Gamma^0$. On the other hand, if at least one of the central multipliers is negative, the sign of $x^+u^-$ changes when one replaces
the pair $(M^+_1,M^-_2)$ by the points $(F_{1}(M^+_1),M^-_2)$ or $(M^+_1,F_{2}^{-1}(M^-_2))$ on the same orbit. Thus, a non-degenerate heterodimensional cycle has either type I or type II, and not type III, if and only if both central multipliers are positive, 
and it has type III if and only if at least one of the central multipliers is negative.

\subsection{Main results for the saddle case}
\label{sec:msaddle}
The key observation in our proof of Theorem \ref{thm:persis_hdd_2p} in the saddle case and the fundamental reason behind the emergence
of robust heterodimensional dynamics is given by the following result proven in Section \ref{sec:blendersintype1}.
\begin{thm}\label{thm:blender_in_type1}
In the saddle case, in any neighborhood of a non-degenerate heterodimensional cycle $\Gamma$ of type I (including type III) for which the value of
$\theta=-\frac{\ln|\lambda|}{\ln|\gamma|}$ is irrational, 
there exists a  standard blender, center-stable with index $d_1$ if $|\alpha|<1$ or center-unstable with index $d_2$ if $|\alpha|>1$.
\end{thm}

The result holds true for systems $f$ of class at least $C^2$. The blender is not the one constructed in \cite{BD08} by means of
a $C^1$-small but not $C^2$-small perturbation of $f$. We do not perturb $f$, but give explicit conditions for the existence of the blender.
Moreover, the (at least) $C^2$ regularity is important for the proof, and it is not clear whether Theorem \ref{thm:blender_in_type1} holds
when $f$ is only $C^1$. Namely, it is a priori 
possible that there could exist $C^1$ systems for which a neighborhood of a heterodimensional cycle of
type I does not contain a blender even when $\theta$ is irrational.\\

Next theorem tells us when the blender of Theorem \ref{thm:blender_in_type1} is {\em activated}, implying that it gets involved in robust heterodimensional dynamics. Recall that by definition the blender exists for any system $C^1$-close to $f$.

\begin{thm}\label{thm:type1_family}
Let $\Gamma$ be a non-degenerate type-I cycle and let $\theta$ be irrational. Consider a sufficiently small $C^r$-neighborhood $\mathcal{V}$ of $f$
such that the blender given by Theorem~\ref{thm:blender_in_type1} persists for any system $g\in\mathcal{V}$.
Let $\mu$ be the splitting functional. 
\begin{itemize}
\item In the case $|\alpha|<1$, there exist constants $C_1 < C_2$ such that for all sufficiently large $m\in\mathbb{N}$ any system $g\in\mathcal{V}$ 
which satisfies $\mu \gamma^m \in  [C_1,C_2]$ has $C^1$-robust heterodimensional dynamics involving the index-$d_1$ cs-blender $\Lambda^{cs}$ of Theorem~\ref{thm:blender_in_type1} and a non-trivial, index-$d_2$ hyperbolic basic set containing $O_2$.
\item In the case $|\alpha|>1$, there exist constants $C_1 < C_2$ such that for all sufficiently large $k\in\mathbb{N}$ any system $g\in\mathcal{V}$ 
which satisfies $\mu \lambda^{-k} \in  [C_1,C_2]$ has $C^1$-robust heterodimensional dynamics involving the index-$d_2$ cu-blender $\Lambda^{cu}$ of Theorem~\ref{thm:blender_in_type1} and a non-trivial, index-$d_1$ hyperbolic basic set containing $O_1$.
\end{itemize}
\end{thm}

The theorem is proven in Section \ref{sec:sta_of_cycles} (see Proposition \ref{prop:intervalsfortype1}); Theorems \ref{thm:type1_family2} - \ref{thm:type2_2pfamily} below are proven there as well.
Note that the cases $|\alpha|<1$ and $|\alpha|>1$ are reduced to each other by the reversion of time and the interchange of the points $O_1$ and $O_2$.
Theorem \ref{thm:type1_family} immediately implies Theorem \ref{thm:persis_hdd_2p} in the case of type-I cycles. Indeed, in a proper unfolding of $f$ we can, by an arbitrarily small increment, make $\theta$ irrational while keeping 
$\mu=0$, and then put $\mu$ to an interval corresponding to the $C^1$-robust heterodimensional dynamics.

We also show (see Proposition \ref{prop:intervalsfortype1}) that {\em there exist intervals of $\mu$ for which the blender $\Lambda^{cs}$ is 
homoclinically related to $O_1$ if $|\alpha|<1$, and the blender $\Lambda^{cu}$ is homoclinically related to $O_2$
if $|\alpha|>1$}. Recall that a hyperbolic point is homoclinically related to a hyperbolic basic set of the same index if their stable and unstable manifolds intersect transversely. If the blender is homoclinically related to a saddle $O_1$ or $O_2$ and is, simultaneously, involved in robust heterodimensional dynamics with the other saddle, this would give robust heterodimensional dynamics involving both these saddles. However, the following result shows that if the central multipliers $\lambda$ and $\gamma$ are both positive, this does not happen within a small neighborhood of the cycle $\Gamma$ under consideration.

Let $f$ have a non-degenerate heterodimensional cycle $\Gamma$ of type I  (we do not
insist now that $\theta$ is irrational). Let  $\lambda>0$ and $\gamma>0$, i.e., $\Gamma$ is not type-III. Let $U$ be a small neighborhood of $\Gamma$. 
\begin{thm}\label{thm:type1_family2} 
One can choose the sign of the splitting functional $\mu$ such that for every system $g$ from a small $C^r$-neighborhood of $f$
\begin{itemize}
\item in the case $|\alpha|<1$, for $\mu(g)>0$, the set of all points whose orbits lie entirely in $U$ consists of a hyperbolic set $\Lambda$ of index $d_1$ (this set includes the cs-blender of Theorem \ref{thm:blender_in_type1} and the orbit $L_1$ of $O_1$), the orbit $L_2$ of the periodic point $O_2$, and heteroclinic orbits corresponding to the transverse intersection of $W^u(L_2)$ with $W^s(\Lambda)$, so there are no heterodimensional dynamics in $U$;\\
for $\mu(g)\leq 0$, no orbit in $W^u(L_1)$ stays entirely in $U$, except for $L_1$ itself and, at $\mu(g)=0$, the fragile heteroclinic $\Gamma^0$, so $L_1$ cannot be a part of any heterodimensional cycle in $U$ when $\mu(g)<0$;
\item in the case $|\alpha|>1$, for $\mu(g)<0$, the set of all points whose orbits lie entirely in $U$ consists of a hyperbolic set $\Lambda$ of index $d_2$ (this set includes the cu-blender of Theorem \ref{thm:blender_in_type1} and the orbit $L_2$ of $O_2$), the orbit $L_1$ of the periodic point $O_1$, and heteroclinic orbits corresponding to the transverse intersection of $W^s(L_1)$ with $W^u(\Lambda)$, so there are no heterodimensional dynamics in $U$;\\
for $\mu(g)\geq 0$, no orbit in $W^s(L_2)$ stays entirely in $U$, except for $L_2$ itself and, at $\mu(g)=0$, the fragile heteroclinic $\Gamma^0$, so $L_2$ cannot be a part of any heterodimensional cycle in $U$ when $\mu(g)>0$.
\end{itemize}
\end{thm}

This situation changes if the type-I cycle is accompanied by a type-II cycle in the following sense.
\begin{defi}[Tied cycles]\label{defi:entangled}
We say that two non-degenerate heterodimensional cycles associated with $O_1$ and $O_2$ are {\em tied} if they share the same fragile heteroclinic, and the
robust heteroclinic orbits $\Gamma^1$ and $\tilde \Gamma^1$ belonging to the corresponding cycles $\Gamma$ and
$\tilde\Gamma$ intersect the same leaf of the foliation
$\mathcal{F}^{ss}$ or the same leaf of the foliation $\mathcal{F}^{uu}$. Specifically, there exists a pair of points $M^+_1=(x^+,0,z^+) \in \Gamma^1\cap W^s_{loc}(O_1)$ and
$\tilde M^+_1=(\tilde x^+,0,\tilde z^+)\in\tilde \Gamma^1\cap W^s_{loc}(O_1)$ such that $x^+=\tilde x^+$ or a pair of points 
$M^-_2=(u^-,0,w^-) \in \Gamma^1\cap W^u_{loc}(O_2)$ and $\tilde M^-_2=(\tilde u^-,0,\tilde w^-)\in\tilde \Gamma^1\cap W^u_{loc}(O_2)$ such that $u^-=\tilde u^-$.
\end{defi}

\begin{figure}[!h]
\begin{center}
\includegraphics[width=0.8\textwidth]{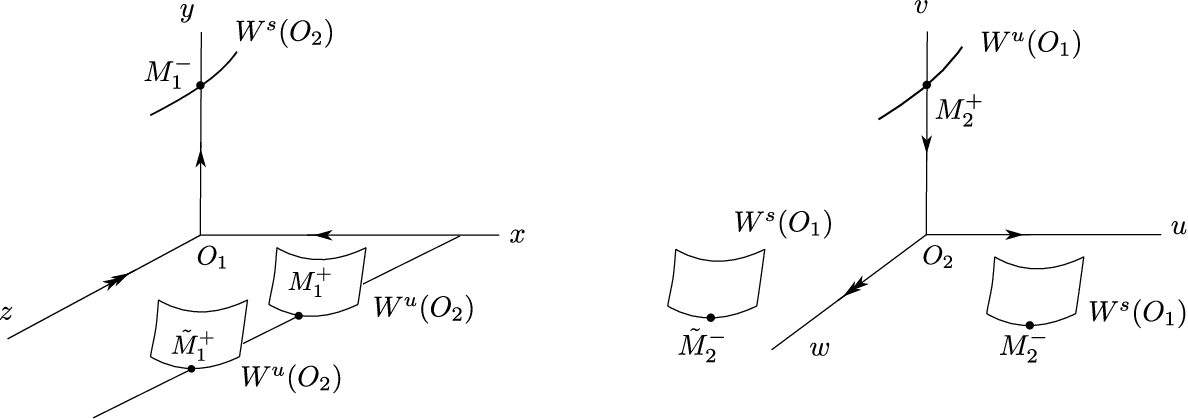}
\end{center}
\caption{A pair of tied cycles intersecting the same strong-stable leaf.}
\label{fig:hdc2}
\end{figure}

The existence of tied cycles is a $C^r$-open property in the set of systems for which the fragile heteroclinic is preserved.
Indeed, if for a system $f$ we have two points $M^+_1\in \Gamma^1$ and $\tilde M^+_1\in \tilde\Gamma^1$ lying in a common leaf $l^{ss}$ of the foliation
$\mathcal{F}^{ss}$, then there are curves $\ell_1$ and $\tilde \ell_1$ containing these points, which correspond to the transverse intersection of $W^u(O_2)$ and $W^s_{loc}(O_1)$ and which, by condition C2, are transverse to the leaf $l^{ss}$ of $\mathcal{F}^{ss}$. The transversality implies that
a $C^r$-small perturbation of $f$ does not destroy this double intersection in $l^{ss}$. The same is true if we have a double intersection with a leaf of
$\mathcal{F}^{uu}$.

\begin{thm}\label{thm:tied} Let a $C^r$ ($r\geq 2$) system $f$ have
a non-degenerate type-I cycle $\Gamma$ tied with a non-degenerate type-II cycle $\tilde \Gamma$.  Assume that 
$\theta$ is irrational. Then,
for any generic one-parameter unfolding $f_\mu$ there exists a converging to $\mu=0$ sequence of intervals $I_j$ such that $f_\mu$ at $\mu\in I_j$
has  
$C^1$-robust heterodimensional dynamics involving the blender given by Theorem \ref{thm:blender_in_type1} 
and a non-trivial hyperbolic basic set; 
of these two hyperbolic sets, the one with index $d_1$ is homoclinically related to $O_1(\mu)$ and the one with index $d_2$ is homoclinically related to $O_2(\mu)$.
\end{thm}

Observe that a type-III cycle is, by definition, a cycle of type I and II, and, obviously, it is tied with itself. Hence, applying
the above theorem, we obtain
\begin{cor}\label{cor:type3_1pfamily}
Let $\Gamma$ be a non-degenerate cycle with real central multipliers $\lambda$ and $\gamma$, at least one of which is negative. 
If $\theta=-\frac{\ln|\lambda|}{\ln|\gamma|}$ is irrational, then 
for any generic one-parameter unfolding $f_\mu$ there exist converging to zero intervals of $\mu$ corresponding to $C^1$-robust heterodimensional dynamics
involving non-trivial hyperbolic basic sets, one of which contains $O_1(\mu)$ and the other contains $O_2(\mu)$. 
\end{cor}

\begin{figure}[!h]
\begin{center}
\includegraphics[width=0.8\textwidth]{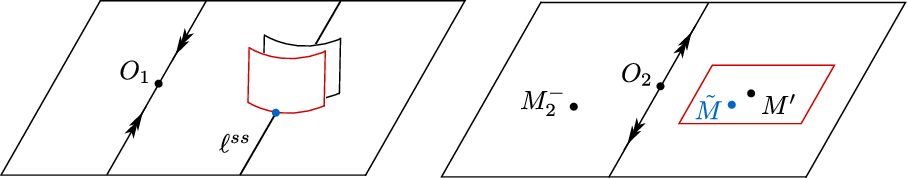}
\end{center}
\caption{Creation of a new cycle that is tied with the original one by the strong-stable leaf $\ell^{ss}$ passing through $M^+_1$.}
\label{fig:tiedhomo}
\end{figure}

\begin{rem}\label{rem:tiedcycles}
Tied cycles also occur when $O_1$ or $O_2$ have a transverse homoclinic. For example, let us have a non-degenerate heterodimensional cycle
$\Gamma$ with a fragile heteroclinic $\Gamma^0$ and a robust heteroclinic $\Gamma^1$. Assume the central multipliers are real, and let $M^\prime\in W^u_{loc}(O_2)$ be a 
point of transverse intersection of the stable and unstable manifolds of $O_2$. If we take a small piece of $W^u_{loc}(O_2)$ around $M^\prime$,
its forward images converge to the entire unstable manifold of $O_2$. Therefore, some of them must intersect transversely
the strong-stable leaf of the point $M^1_+\in \Gamma^1\cap W^s_{loc}(O_1)$ (as this leaf intersects $W^u(O_2)$ transversely at the point
$M^1_+$ by condition C2), see Figure \ref{fig:tiedhomo}. The orbit of the intersection point is a robust heteroclinic $\tilde\Gamma^1$, and the corresponding cycle
$\tilde \Gamma$ is tied with $\Gamma$. By construction, the orbit $\tilde\Gamma^1$ has 
a point $\tilde M$ in $W^u_{loc}(O_2)$ close to the homoclinic point $M^\prime$. Therefore, its $u$-coordinate is close to the $u$-coordinate $u^\prime$ of
$M^\prime$. Therefore, if $u^-u^\prime<0$, i.e., the homoclinic point $M^\prime$ and the point $M^-_2$ of $\Gamma^1$ lie in 
$W^u_{loc}(O_2)$ on opposite sides from $W^{uu}_{loc}(O_2)$, then the tied cycles $\Gamma$ and $\tilde\Gamma$ have different types, and Theorem \ref{thm:tied} is applicable.
\end{rem}

Theorem \ref{thm:persis_hdd_2p} for type-II cycles is inferred from Theorem \ref{thm:tied} by means of the following result.
\begin{thm}\label{thm:type2_2pfamily}
Let $f$ have a non-degenerate type-II cycle $\Gamma$ with irrational $\theta$.
For any generic one-parameter unfolding $f_\mu$ there exists a sequence $\mu_j \to 0$ such that $f_\mu$ at 
$\mu=\mu_j$ has a pair of tied heterodimensional cycles $\Gamma_{j,I}$ and $\Gamma_{j,II}$ of type I and type II, which are associated 
\begin{itemize}
\item with $O_1(\mu)$ and an index-$d_2$ saddle $O^\prime_2(\mu)$ which is homoclinically related to $O_2(\mu)$ if $|\alpha|<1$; or
\item with $O_2(\mu)$ and an index-$d_1$ saddle $O^\prime_1(\mu)$ which is homoclinically related to $O_1(\mu)$ if $|\alpha|>1$.
\end{itemize}
\end{thm}

In order to apply Theorem \ref{thm:tied} to the tied cycles $\Gamma_{j,I}$ and $\Gamma_{j,II}$ obtained in Theorem \ref{thm:type2_2pfamily}, we extend $f_\mu$ to a proper, at least two-parameter unfolding 
$f_\varepsilon$. In Section \ref{sec:thm:type2_2pfamily}, we make the following observation:  {\em when $r\geq 3$, 
the same family $f_\eps$ gives a proper unfolding for the cycles $\Gamma_{j,I}$ and $\Gamma_{j,II}$}, see Lemma \ref{lem:type2_sechdc} and equation \eqref{eq:type2_2pfamily:1}.
Since $f_\varepsilon$ is proper for the cycles $\Gamma_{j,I}$ and $\Gamma_{j,II}$, one can always find the values of $\varepsilon$
for which the value of $\theta$ for these cycles is irrational. Hence, applying Theorem \ref{thm:tied}, we obtain 
the result of Theorem \ref{thm:persis_hdd_2p} when $f_\varepsilon$ is at least $C^3$.  

In the $C^2$-case, the reduction of Theorem \ref{thm:persis_hdd_2p} (for type-II cycles) to Theorem \ref{thm:type2_2pfamily} requires a revision of Theorem \ref{thm:tied}, as described in Remark \ref{rem:smoothness6}. The difficulty 
is that we use, for every parameter value, the coordinates which linearize the action of the local maps $F_{1}$ and $F_{2}$
in the central direction. It is known that in the $C^2$-case the linearizing coordinate transformation is, in general, not smooth with respect
to parameters, so our technique does not allow to compute derivatives with respect to $\varepsilon$ which enter the definition of a proper unfolding. Instead, we use continuity arguments to show in the $C^2$-case that, still, for the tied cycles $\Gamma_{j,I}$ and $\Gamma_{j,II}$ the value of $\theta$ can be made irrational and the splitting parameter for these cycles
can be pushed, by a small change of $\varepsilon$, inside the open regions described by 
Remark \ref{rem:smoothness6} -- analogues of intervals $I_j$ described by Theorem \ref{thm:tied}.

Altogether, 
we prove in Section \ref{sec:thm:type2_2pfamily} Theorem \ref{thm:persis_hdd_2p} for type-II cycles in the following form.

\begin{cor}\label{cor:type2_2pfamily2}
Let $\Gamma$ be a non-degenerate type-II cycle, and let $f_\varepsilon$ be a proper, at least two-parameter unfolding. Arbitrarily close to $\varepsilon=0$
there exist open regions in the parameter space for which the corresponding system $f_\varepsilon$ has $C^1$-robust heterodimensional dynamics involving a standard blender and a non-trivial hyperbolic basic set; one of these sets is homoclinically related to $O_1$ and the other is homoclinically related to $O_2$.
\end{cor}

As we see, {\em 
the emergence of heterodimensional dynamics depends strongly on the arithmetic properties of $\theta$}, the 
modulus of topological equivalence. The following result shows that, in the saddle case, we have a clear dichotomy: 
for irrational $\theta$ we have highly non-trivial dynamics and bifurcations in any neighborhood of the heterodimensional cycle $\Gamma$, and for rational $\theta$ the dynamics in a small neighborhood of $\Gamma$ are quite simple, in general.
\begin{thm}\label{thm:rationaltheta}
Let a $C^r$ $(r\geq 2)$ system $f$ have a non-degenerate heterodimensional cycle $\Gamma$,
and let the central multipliers be real, $|\lambda|<1$ and $|\gamma|>1$.
Let $\theta=-\frac{\ln|\lambda|}{\ln|\gamma|}$ be rational, i.e., $|\gamma|=|\lambda|^{-\frac{p}{q}}$ for some coprime integers $p>0$, $q>0$.
Suppose the following conditions are fulfilled:
\begin{equation}\label{rare1}
|ab|\neq |\gamma|^{\frac{s}{q}} \quad \mbox{ for } \;\; s\in\mathbb{Z},
\end{equation}
\begin{equation}\label{rare2}
\left|\frac{u^-}{ax^+}\right|\not\in 
cl\left\{|\gamma|^{\frac{s}{q}} \dfrac{1-\lambda^{l}}{1-\gamma^{-n}} \right\}_{s\in\mathbb{Z},l\in\mathbb{N},n\in\mathbb{N}}\;.
\end{equation}
Let $U$ be a sufficiently small neighborhood of $\Gamma$ and let $\mathcal{N}$ be the set of all orbits that lie entirely in
$U$. Then, at $\mu=0$,  the set $\mathcal{N}$ is the union of $L_1$, $L_2$, $\Gamma^0$, and the orbits of transverse intersection of $W^u(L_2)$ with $W^s(L_1)$ near $\Gamma^1$.\\\\
For any generic one-parameter unfolding $f_\mu$, for any small $\mu\neq 0$, either
\begin{itemize}
\item $\mathcal{N}$ is comprised by $L_2$, an index-$d_1$ uniformly-hyperbolic compact set $\Lambda_1$ containing $L_1$, and transverse heteroclinic
connections between $W^u(L_2)$ and $W^s(\Lambda_1)$, while no heteroclinic connection between $W^u(\Lambda_1)$ and $W^s(L_2)$
exists, or 
\item $\mathcal{N}$ is comprised by $L_1$, an index-$d_2$ uniformly-hyperbolic compact set $\Lambda_2$ containing $L_2$, and transverse heteroclinic
connections between $W^u(\Lambda_2)$ and $W^s(L_1)$, while no heteroclinic connection between $W^u(L_1)$ and $W^s(\Lambda_2)$ exists.
\end{itemize}
\end{thm}
The proof of this theorem is given in Section \ref{sec:rational}. Notice that, for fixed values of $\lambda$ and $\gamma$, conditions (\ref{rare1}) and (\ref{rare2}) are fulfilled for all $ab$ and ${ax^+}/{u^-}$ except for a countable, nowhere dense set of values. Thus, the simplicity of dynamics at rational $\theta$ is indeed quite generic. It also follows from this theorem and Theorem \ref{thm:blender_in_type1}
that whenever we have a heterodimensional cycle $\Gamma$ of type I, if we change $\theta$ without destroying $\Gamma$, the blender that forms at irrational $\theta$'s immediately departs from $\Gamma$, so that for each rational $\theta$ a sufficiently small ``blender-free'' neighborhood of $\Gamma$ emerges.

\subsection{The case of nonreal central multipliers}
\label{sec:mcom}
In the remaining saddle-focus and double-focus cases, we obtain Theorem \ref{thm:persis_hdd_2p} from
\begin{thm}\label{thm:sad-foc_main}
Let a $C^r$ ($r\geq 2$) system $f$ have a heterodimensional cycle $\Gamma$ and let central multipliers be $\lambda_{1,1}=\lambda_{1,2}^*=\lambda e^{i\omega}$ and real $\gamma_{2,1}=\gamma$ (the saddle-focus case) or $\lambda_{1,1}=\lambda_{1,2}^*=\lambda e^{i\omega_1}$ and
$\gamma_{2,1}=\gamma_{2,2}^*=\gamma e^{i\omega_2}$ (the double-focus case). Assume, in the saddle-focus case, that
the numbers $\theta=-\frac{\ln\lambda}{\ln|\gamma|}$, $\frac{1}{2\pi}\omega$ and $1$ are rationally independent. In the 
double-focus case, assume that $\theta$, $\frac{1}{2\pi}\omega_1$, $\frac{1}{2\pi}\theta\omega_2$, and $1$ are rationally independent.
Then, in any neighborhood of $\Gamma$, the system $f$ has $C^1$-robust heterodimensional dynamics associated with a standard cs-blender
and a standard cu-blender.
 
Moreover, in the double-focus case, the point $O_1$ is homoclinically related to the cs-blender and
the point $O_2$ is homoclinically related to the cu-blender, and the same holds true for any $C^r$-close system $g$. In the saddle-focus case, the point $O_2$ is, for the system $f$ and for any $C^r$-close system $g$,
homoclinically related to the cu-blender,
and there exist intervals $I_j$ converging to $\mu=0$ such that if $\mu(g)\in I_j$, then
the point $O_1$ is homoclinically related to the cs-blender.
\end{thm}
\begin{rem}\label{rem:sad-foc}
By the reversion of time, we deduce a similar result for the case where $\gamma_{2,1}=\gamma_{2,2}^*=\gamma e^{i\omega}$ and 
$\lambda_{1,1}=\lambda$ is real. Namely, we obtain that the point $O_1$ is homoclinically related to a cs-blender for the system $f$ and 
for any $C^r$-close system $g$, and that there exist intervals $I_j$ converging to $\mu=0$ such that if $\mu(g)\in I_j$, then
the point $O_2$ is homoclinically related to the cu-blender.
\end{rem}

This theorem is proved in Section \ref{sec:complex_ev}. In particular, the coexistence of a cs-blender and a cu-blender for the non-perturbed system $f$ is obtained in Propositions \ref{prop:sad-foc_blender} and \ref{prop:df_blender}. Since in a proper unfolding of $f$ the rational independence conditions of 
Theorem \ref{thm:sad-foc_main} are achieved by arbitrarily small changes of parameters, the claim of Theorem \ref{thm:persis_hdd_2p} follows immediately.

\subsection{Local stabilization of heterodimensional cycles}

Corollary \ref{cor:type3_1pfamily}, Corollary \ref{cor:type2_2pfamily2}, and Theorem \ref{thm:sad-foc_main} immediately imply that if the heterodimensional cycle is not type-I, then it is locally $C^r$-stabilized as claimed in Theorem \ref{thm:stanew} (in the $C^1$ case, this result is obtained in \citep{BDK12}, but certain perturbations essential for their construction are large in the $C^2$-topology). On the other hand, Theorem \ref{thm:type1_family2} shows that any type-I heterodimensional cycle cannot be locally $C^r$-stabilized, concluding Theorem \ref{thm:stanew}. A similar result 
on the impossibility of the local $C^1$-stabilization for type-I cycles can be inferred from \citep{BD12}.

If we replace at least one of the periodic orbits in the heterodimensional cycle by a non-trivial, transitive, compact hyperbolic set $\Lambda$, then even if all periodic orbits in $\Lambda$ have central multipliers real and positive, one can, in the generic situation, find a periodic orbit in $\Lambda$ such that its homoclinic points accumulate to it from both sides in the central direction. Then, by Remark \ref{rem:tiedcycles}
and Theorem \ref{thm:tied} we obtain
\begin{cor}\label{cor:stabilization_lambda}
If at least one of the periodic orbits in a heterodimensional cycle belongs to a non-trivial, transitive, compact hyperbolic set, then the cycle is $C^r$-stabilizable for any $r=2,\dots,\infty,\omega$ 
\end{cor}

\section{First-return maps in the saddle case}\label{sec:firstreturn}
Recall that we consider the behavior in a small neighborhood of the heterodimensional cycle $\Gamma$. We take periodic points $O_1$ and $O_2$ in $\Gamma$
and consider, if $f$ is a discrete dynamical system, small neighborhoods $U_{01}$ and $U_{02}$ of $O_1$ and $O_2$. If $f$ is a flow, then $U_{01}$ and $U_{02}$
are small codimension-1 cross-section to the flow through the points $O_1$ and $O_2$. In both cases, the local maps $F_1$ and $F_2$ act on $U_{01}$
and $U_{02}$, which are defined by the orbits of the system $f$ near the orbits $L_1$ and $L_2$ of the points $O_1$ and $O_2$. We take a pair of points
$M_1^-\in U_{01}$ and $M_2^+\in U_{02}$ on the fragile heteroclinic $\Gamma^0$ and a pair of points 
$M_2^-\in U_{02}$ and $M_1^+\in U_{01}$ on the robust heteroclinic $\Gamma^1$. The orbits near the fragile heteroclinic define the transition map
$F_{12}$ from a small neighborhood of $M_1^-$ in $U_{01}$ to a small neighborhood of $M_2^+$ in $U_{02}$, and the orbits near the robust heteroclinic define the transition map
$F_{21}$ from a small neighborhood of $M_2^-$ in $U_{02}$ to a small neighborhood of $M_1^+$ in $U_{01}$. 

Our immediate goal is to study the
first-return maps near the heterodimensional cycle.
Such map is a composition $T_{k,m}:=F_{21}\circ F_2^m\circ F_{12}\circ F_1^k$, which takes points from a small neighborhood $M_1^+$ in $U_{01}$ back to a vicinity of $M_1^+$. Here, $k$ and $m$ are sufficiently large positive 
integers such that $F_1^k$ takes points from a small neighborhood of $M_1^+$ in $U_{01}$ to a small neighborhood of $M_1^-$ 
in $U_{01}$ and
$F_2^m$ takes points from a small neighborhood of $M_2^+$ in $U_{02}$ to a small neighborhood of $M_2^-$ in $U_{02}$. 

First, we use results in \citep{GST08} to obtain formulas for the iterations $F_1^k$ and $F_2^m$.
The formulas are written in the so-called ``cross-form'', see \eqref{eq:maps:F_1^k} and \eqref{eq:maps:F_2^m}. Next, we use the transversality conditions C1 and C2 and write the transition maps $F_{12}$ and $F_{21}$ 
also in the cross-form, see \eqref{eq:maps:F_12_cross} and \eqref{eq:maps:F_21_cross}. Then, combining the 
cross-form formulas for these four maps, we obtain a formula for the map $T_{k,m}$ in some rescaled variables and show that it is indeed a return map to $U_{01}$ for suitable values of $k$ and $m$, see \eqref{eq:maps:T_k,m_cross_0mu}. Finally, we conclude this section by proving the partial hyperbolicity of $T_{k,m}$ in Lemma \ref{lem:conefields}.

\subsection{Local maps}
First, we discuss necessary estimates for $F_1^k$ and $F_2^m$. By \citep[Lemmas 5 and 6]{GST08}, one can choose local coordinates 
$(x,y,z)\in \mathbb{R}\times\mathbb{R}^{d_1}\times\mathbb{R}^{d-d_1-1}$ in $U_{01}$ such that $O_1$ is at the origin and the map $F_1$ takes the form
\begin{equation}\label{eq:maps:1}
\begin{aligned}
\bar{x}&=\lambda x+ g_1(x,y,z),\\
\bar{y}&=P_1 y + g_2(x,y,z),\\
\bar{z}&=P_2 z +g_3(x,y,z),
\end{aligned}
\end{equation}
where $\lambda=\lambda_{1,1}$, and the eigenvalues of the matrices $P_1$ and $P_2$ are $\gamma_{1,1},\gamma_{1,2},\dots,\gamma_{1,d_1}$ and $\lambda_{1,2}\dots \lambda_{1,d-d_1}$, respectively (see \eqref{eq:coindex1}). The functions $g_{1,2,3}$ vanish along with their first derivatives at the origin
and satisfy the identities
\begin{equation}\label{eq:maps:2}
g_{1,3}(0,y,0)=0,\quad\quad g_2(x,0,z)=0, \quad\quad g_1(x,0,z)=0,  \quad\quad
\dfrac{\partial g_{1,3}}{\partial x}(0,y,0)=0,
\end{equation}
for all sufficiently small $x$, $y$, and $z$. As discussed in Section \ref{sec:intro2}, the first two identities in \eqref{eq:maps:2} imply that the local manifolds $W^s_{loc}(O_1)$ and $W^u_{loc}(O_1)$ are straightened and given by the equations $\{y=0\}$ and $\{x=0,z=0\}$, respectively. The third identity shows that the strong-stable foliation $\mathcal{F}^{ss}$ (which enters condition C2) is straightened, its leaves are given by equations $\{x=const,y=0\}$, and the map $F_1$ restricted to $W^s_{loc}(O_1)$ is linear in $x$. The forth identity implies that the extended unstable manifold
$W^{uE}_{loc}(O_1)$ (which is in condition C1) is tangent to $z=0$ at the points of $W^u_{loc}(O_1)$.

Similarly, we introduce coordinates $(u,v,w)\in\mathbb{R}\times\mathbb{R}^{d-d_1-1}\times\mathbb{R}^{d_1}$ in $U_{02}$ with $O_2$ at the origin such that $F_2$ takes the form
\begin{equation}\label{eq:maps:3}
\begin{aligned}
\bar{u}&=\gamma u + \hat g_1(u,v,w),\\
\bar{v}&=Q_1 v+ \hat g_2(u,v,w),\\
\bar{w}&=Q_2 w +\hat g_3(u,v,w),
\end{aligned}
\end{equation}
where $\gamma=\gamma_{2,1}$, and the eigenvalues of the matrices $Q_1$ and $Q_2$ are $\lambda_{2,2},\lambda_{2,2},\dots,\lambda_{2,d-d_1-1}$ and $\gamma_{2,2},\dots,\gamma_{2,d_1+1}$, respectively (see \eqref{eq:coindex1}). The functions $\hat g_{1,2,3}$ vanish along with their first derivatives at the origin and satisfy
\begin{equation}\label{eq:maps:4}
\hat g_{1,3}(0,v,0)=0,\quad\quad \hat g_2(u,0,w)=0, \quad\quad \hat g_1(u,0,w)=0, \quad\quad
\dfrac{\partial \hat g_{1,3}}{\partial u}(0,v,0)=0,
\end{equation}
for all sufficiently small $u$, $v$ and $w$.

It is shown in \citep{GST08} (see remarks after Lemma 6 there), that  the transformations that bring $F_1$ to form (\ref{eq:maps:1}) and $F_2$ to form \eqref{eq:maps:3} are of class $C^r$, so we 
do not loose regularity when using these coordinates. Moreover, these coordinate transformations depend uniformly-continuously on the system $f$.
When we consider parametric families $f_\eps$, and $f$ is analytic or $C^\infty$ with respect to coordinates and parameters, the coordinate transformations are also analytic or $C^\infty$ with respect to $\varepsilon$. In the case of finite smoothness $r$, we loose, in general, two derivatives with respect to $\varepsilon$. Namely, the second derivative of the transformation
is $C^{r-2}$-smooth with respect to both the coordinates and parameters. Therefore, the matrices $P_{1,2}$ and $Q_{1,2}$ 
in \eqref{eq:maps:1} and \eqref{eq:maps:3} are $C^{r-2}$-functions of $\varepsilon$, and the functions $g_{1,2,3}$ and $\hat g_{1,2,3}$, as well as their derivatives with respect to $(x,y,z)$ or $(u,v,w)$ up to order 2, are $C^{r-2}$-functions of the coordinates and $\varepsilon$. If $r\geq 3$, this gives us at least 1 continuous derivative
with respect to $\varepsilon$. In the $C^2$ case, we can only assume continuity of $P_{1,2}$, $Q_{1,2}$ with respect to $\varepsilon$, the same goes for $g_{1,2,3}$, 
$\hat g_{1,2,3}$, and their derivatives with respect to the coordinates.
The eigenvalues $\lambda$ and $\gamma$ do not depend on coordinate transformations; as they are the eigenvalues of the first derivative of $F_1$ or $F_2$ (at $O_1$ and, respectively, $O_2$), they are at least $C^1$ with respect to $\varepsilon$ in any case. 

Take any point $(x,y,z)$ in $U_{01}$ and let $(\tilde x,\tilde y, \tilde z)=F_1^k(x,y,z)$. It is known (see e.g. \citep{ASh73,Sh67}) that the value of $(\tilde x,y,\tilde z)$ is uniquely defined by $(x,\tilde y, z)$ for all $k\geqslant 0$. By \citep[Lemma 7]{GST08}, when identities \eqref{eq:maps:2} are fulfilled, the relation between the coordinates can be written as
\begin{equation}\label{eq:maps:F_1^k}
\begin{aligned}
\tilde x&=\lambda^k x + p_1(x,\tilde y,z),\\
y &= p_2(x,\tilde y,z),\qquad \tilde z =p_3(x,\tilde y,z),
\end{aligned}
\end{equation}
where 
\begin{equation}\label{eq:maps:F_1^k:derivatives}
\|p_{1,3}\|_{C^1}=o(\lambda^k), \qquad\qquad \|p_2\|_{C^1}=o(\hat\gamma^{-k}),
\end{equation}
for some constant $\hat\gamma\in (1, |\gamma_{1,1}|)$. These estimates are uniform for all systems $C^2$-close to $f$; when we consider
parametric families $f_\varepsilon$, the functions $p_{1,2,3}$ depend on $\varepsilon$ uniformly-continuously, along with their first derivatives with respect to
$(x,\tilde y,z)$. In the case $r\geq 3$, we have the same $o(\lambda^k)$ and $o(\hat\gamma^{-k})$ estimates for the derivatives of $p_{1,2}$ and, respectively, $p_3$
with respect to parameters $\varepsilon$, see \cite[Lemma 7]{GST08} for a detailed discussion.

Likewise, for any $(u,v,w)\in U_{02}$ we have $(\tilde u,\tilde v,\tilde w)=F^m_2(u,v,w)$ if and only if 
\begin{equation}\label{eq:maps:F_2^m}
\begin{aligned}
u&=\gamma^{-m} \tilde u +q_1(\tilde u,v,\tilde w),\\
\tilde v &= q_2(\tilde u,v,\tilde w),\qquad  w = q_3(\tilde u,v,\tilde w),
\end{aligned}
\end{equation}
where 
\begin{equation}\label{eq:maps:F_2^m:derivatives}
\|q_{1,3}\|_{C^1}=o(\gamma^{-m}),\qquad\qquad \|q_2\|_{C^1}=o(\hat\lambda^m),
\end{equation}
for some constant $\hat\lambda\in (|\lambda_{2,1}|,1)$; the functions $q_{1,2,3}$ depend uniformly-continuously on the system $f$, and the estimates
(\ref{eq:maps:F_2^m:derivatives}) hold uniformly for all systems $C^2$-close to $f$. When $r\geq 3$, estimates \eqref{eq:maps:F_2^m:derivatives}
also hold for the derivatives with respect to parameters $\varepsilon$.

\subsection{Transition maps}
Next, we consider the transition maps $F_{12}$ and $F_{21}$. We use the following notation for the coordinates of the points $M_{1,2}^\pm$: 
$$M^+_1=(x^+,0,z^+),\quad M^-_1=(0,y^-,0),\quad M^+_2=(0,v^+,0), \quad M^-_2=(u^-,0,w^-).$$
We can write the Taylor expansion of the transition map $F_{12}:(\tilde x,\tilde y,\tilde z)\mapsto (u,v,w)$ near $M^-_1$ as
\begin{equation}\label{eq:maps:F_12}
\begin{aligned}
u &= a'_1 +a'_{11} \tilde x + a'_{12}(\tilde y-y^-) + a'_{13}\tilde z + O(\tilde x^2 + (\tilde y-y^-)^2+\tilde z^2), \\
v-v^+ &= a'_2+a’_{21} \tilde x + a'_{22}(\tilde y-y^-) + a'_{23}\tilde z + O(\tilde x^2 + (\tilde y-y^-)^2+\tilde z^2), \\
w &= a'_3+a'_{31} \tilde x + a'_{32}(\tilde y-y^-) + a'_{33}\tilde z + O(\tilde x^2 + (\tilde y-y^-)^2+\tilde z^2),
\end{aligned}
\end{equation}
where $a'_i$ and $a'_{ij}$ $(i,j=1,2,3)$ are some constants. 
Since $F_{12}(M^-_1)=M^+_2$ for the system $f$, it follows that $a'_{1,2,3}$ vanish, but when we perturb $f$, these coefficients can become non-zero (though small). 

Recall that $W^{sE}_{loc}(O_2)$ is tangent at $M_2^+$ to $w=0$ and $W^u_{loc}(O_1)$ is given by $(\tilde x=0,\tilde z	=0)$.
Thus, the transversality condition C1 writes as the uniqueness of the (trivial) solution to the system
$$ u = a'_{12}(\tilde y -y^-),\qquad v-v^+= a'_{22}(\tilde y -y^-),\qquad 0 = a'_{32}(\tilde y-y^-),$$
i.e., $a'_{32}$ is invertible. It follows that $\tilde y - y^-$ can be expressed, from the last equation of \eqref{eq:maps:F_12}, as a smooth function of
$(w,\tilde x,\tilde z)$, so the map $F_{12}$ can be written in the cross-form as
\begin{equation}\label{eq:maps:F_12_cross}
\begin{aligned}
u &=\hat \mu+ a \tilde x + a_{12}w + a_{13}\tilde z + O(\tilde x^2 + w^2+\tilde z^2), \\
v-\hat v^+ &=a_{21} \tilde x + a_{22}w + a_{23}\tilde z + O(\tilde x^2 + w^2+\tilde z^2), \\
\tilde y-\hat y^- &=a_{31} \tilde x + a_{32}w + a_{33}\tilde z + O(\tilde x^2 + w^2+\tilde z^2),
\end{aligned}
\end{equation}
where the coefficients $\hat v^+$, $\hat y^-$, $\hat\mu$, $a$, and $a_{ij}$ change uniformly-continuously when the system $f$ is perturbed, and for the original system $f$, we have $\hat v^+=v^+$, $\hat y^-=y^-$, and $\hat\mu=0$. Since it does not cause ambiguity, in further references to \eqref{eq:maps:F_12_cross} we use $v^+$ and $y^-$ instead of $\hat v^+$ and $\hat y^-$. Note that the coefficient $a\neq 0$ is exactly the derivative defined in \eqref{eq:a11}.

Since $W^u_{loc}(O_1)$ is given by $\{(\tilde x=0,\tilde z=0)\}$ and 
$W^s_{loc}(O_2)$ is given by $\{(u=0,w=0)\}$, a fragile heteroclinic corresponding to the intersection of $F_{12}(W^u_{loc}(O_1))$ and $W^s_{loc}(O_2)$
persists for a perturbation of $f$ if and only if $\hat\mu=0$. More precisely, $\hat\mu$ is the $u$-coordinate of the point of the intersection of 
$F_{12}(W^u_{loc}(O_1))$ with $\{w=0\}$. This intersection
is transverse by condition C1, from which one infers that the ratio of $|\hat\mu|$ to the distance between $F_{12}(W^u_{loc}(O_1))$ and $W^s_{loc}(O_2)$ tends (uniformly in some $C^2$-neighborhood of $f$) to a finite non-zero limit value when $\hat\mu\to 0$. Recall that we defined the splitting parameter $\mu(f)$, whose absolute value equal the distance between $F_{12}(W^u_{loc}(O_1))$ and $W^s_{loc}(O_2)$ and which enters
Theorems \ref{thm:type1_family} - \ref{thm:rationaltheta}. By scaling the variable $u$, we can always obtain
\begin{equation}\label{mmoh}
\lim_{\mu\to 0} \dfrac{\hat\mu}{\mu(f)}  =1, \;\;\mbox{ or }\;\; \hat\mu=\mu+o(\mu).
\end{equation}
If $r\geq 3$, then when we consider parametric families $f_\eps$,
the coefficients of \eqref{eq:maps:F_12_cross} are at least $C^1$ with respect to $\eps$. When the family is generic or proper, $\mu=\mu(f_\eps)$ is one of the parameters,
and we have (for $r\geq 3$)
$$\left.\frac{\partial \hat \mu}{\partial \mu}\right|_{\mu=0} =1.$$
It follows that in this case we can make a smooth change of parameters such that $\hat\mu=\mu$.

The Taylor expansion of the other transition map $F_{21}:(\tilde u,\tilde v, \tilde w)\mapsto ( x, y, z)$, which is defined for $(\tilde u, \tilde v, \tilde w)$ 
near $M^-_2=(u^-,0,w^-)$ and takes values $( x, y, z)$ near $M^+_1=(x^+,0,z^+)$, is given by
\begin{equation}\label{eq:maps:F_21}
\begin{aligned}
 x-x^+ =& b'_{11} (\tilde u-u^-) + b'_{12}\tilde v + b'_{13}(\tilde w-w^-) + O((\tilde u-u^-)^2 + \tilde v^2+(\tilde w-w^-)^2), \\
 y =& b'_{21} (\tilde u-u^-) + b'_{22}\tilde v + b'_{23}(\tilde w-w^-) + O((\tilde u-u^-)^2 + \tilde v^2+(\tilde w-w^-)^2), \\
 z-z^+ =& b'_{31} (\tilde u-u^-) + b'_{32}\tilde v + b'_{33}(\tilde w-w^-) + O((\tilde u-u^-)^2 + \tilde v^2+(\bar w-w^-)^2),
\end{aligned}
\end{equation}
where $b'_{ij}$ are some constants. Arguing as for the map $F_{12}$ above, one can use the assumption that $F^{-1}_{21}(W^s_{loc}(O_1))\pitchfork \mathcal{F}^{uu}\neq\emptyset$ from condition C2 to deduce that $\det b'_{23}\neq 0$ (note that the leaf of the foliation $\mathcal{F}^{uu}$ through $M^-_2$ is given by $(\tilde u=u^-,\tilde v=0)$
and $W^s_{loc}(O_1)$ is given by $ \{y=0\}$ here). Consequently, $F_{21}$ can be written in the following cross-form:
\begin{equation}\label{eq:maps:F_21_cross}
\begin{aligned}
 x -x^+ &= b (\tilde u-u^-) + b_{12}\tilde v + b_{13} y + O((\tilde u-u^-)^2 + \tilde v^2+ y^2), \\
\tilde w-w^- &= b_{21} (\tilde u-u^-) + b_{22}\tilde v + b_{23} y + O((\tilde u-u^-)^2 + \tilde v^2+ y^2), \\
 z -z^+ &= b_{31} (\tilde u-u^-) + b_{32}\tilde v + b_{33} y  + O((\tilde u-u^-)^2 + \tilde v^2+ y^2),
\end{aligned}
\end{equation}
where $b$ is the derivative defined by \eqref{eq:intro:b}. All the coefficients in \eqref{eq:maps:F_21_cross} change uniformly-continuously when the system $f$ is perturbed; if $r\geq 3$, then they are at least $C^1$ with respect to the perturbation parameters $\eps$.
Note that since $b\neq 0$, we can express $\tilde u-u^-$ as a function of $( x -x^+,\tilde v, y)$ and rewrite \eqref{eq:maps:F_21_cross} as
\begin{equation}\label{f21crc}
\begin{aligned}
\tilde u - u^- &= b^{-1} ( x -x^+ - b_{13} y) + O(\|\tilde v\|+ ( x-x^+)^2+ y^2), \\
\tilde w - w^- &= O(| x -x^+|+\|\tilde v\|+\| y\|), \\
 z - z^+ &= b_{31}b^{-1}(x-x^+-b_{13}y)+b_{33}y + O(\|\tilde v\|+ ( x-x^+)^2+ y^2).
\end{aligned}
\end{equation}

\subsection{First-return maps and cone field lemma}
We can now find a formula for the first-return map $T_{k,m}=F_{21}\circ F_2^m\circ F_{12}\circ F_1^k$. Combining \eqref{eq:maps:F_1^k} and \eqref{eq:maps:F_12_cross}, we obtain a formula for the map $F_{12}\circ F_1^k$. Namely, by substituting the 
last equation of \eqref{eq:maps:F_12_cross} ($\tilde y$ as a function of $\tilde x$, $w$, and $\tilde z$) into the first and the last equations of \eqref{eq:maps:F_1^k},
we can express $\tilde x$ and $\tilde z$ as functions of $(x,z,w)$:
$$\tilde x=\lambda^k x + o(\lambda^k), \qquad \tilde z=o(\lambda^k).$$
After that, we substitute these formulas into the rest equations in \eqref{eq:maps:F_1^k} and \eqref{eq:maps:F_12_cross}, and find that there exist smooth functions
$$\tilde h_1(x,z,w)=O(\|w\|)+o(\lambda^k), \qquad \tilde h_2(x,z,w)=o(\hat\gamma^{-k}), \qquad  \tilde h_3(x,z,w)=O(\|w\|+|\lambda|^k),$$
such that, for sufficiently large $k$, a point $(x,y,z)$ from a small neighborhood of $M_1^+$ is taken by the map $F_{12}\circ F_1^k$ to a point $(u,v,w)$ in
a small neighborhood of $M_2^+$ if and only if
\begin{equation}\label{eq:maps:112}
\begin{aligned}
u &=\hat\mu+ a \lambda^k x + \tilde h_1(x,z,w)  \\
y &=\tilde h_2(x,z,w),\qquad v- v^+ =  \tilde h_3(x,z,w).
\end{aligned}
\end{equation}

Similarly, combining \eqref{eq:maps:F_2^m} and \eqref{f21crc} yields that there exist smooth functions 
\begin{equation}\label{htha}
\hat h_{0i}( x-x^+,\bar y)=O( (x-x^+)^2+ y^2) \;\;\ (i=1,2), \qquad
\hat h_{1,3}( x, v, y)=o(\gamma^{-m}), \qquad \hat h_2( x, v, y)=o(\hat\lambda^m),
\end{equation}
such that, for sufficiently large $m$, a point $(u,v,w)$ from a small neighborhood of $M_2^+$ is taken
by the map $F_{21}\circ F_2^m$ to a point $( x, y, z)$ in a small neighborhood of $M_1^+$ if and only if
\begin{equation}\label{eq:maps:221}
\begin{aligned}
u &= \gamma^{-m}(u^- + b^{-1} ( x -x^+ - b_{13}  y + \hat h_{01}( x- x^+, y)))  + \hat h_1( x,v, y), \\[5pt]
 z - z^+ &= b_{31}b^{-1}(x-x^+-b_{13}y)+b_{33}y + \hat h_{02}( x-x^+, y)   + \hat h_2( x,v, y),\\[5pt]
w &= \hat h_3( x, v, y).
\end{aligned}
\end{equation}

Let us now consider the first-return map $T_{k,m}=F_{21}\circ F_2^m\circ F_{12}\circ F_1^k$ for sufficiently large positive integers $k$ and $m$ that takes a point $(x,y,z)$ from a small neighborhood of $M_1^+$ to its image $(\bar x,\bar y,\bar z)$. To obtain a formula for this map, we replace $(x,y,z)$ in \eqref{eq:maps:221} by $(\bar x,\bar y,\bar z)$ and combine it with \eqref{eq:maps:112}. More specifically, from the system comprised by the last equations of \eqref{eq:maps:112} (for $v$) and \eqref{eq:maps:221} (for $w$) we can
express $v$ and $w$ as smooth functions of $(x,\bar x, z, \bar y)$. By substituting these functions into the rest of the equations 
\eqref{eq:maps:112} and \eqref{eq:maps:221}, we find that there exist functions 
$$h_1(x,\bar x, z, \bar y)=o(\lambda^k)+o(\gamma^{-m}), \qquad h_2(x,\bar x, z, \bar y)=o(\hat\lambda^m), \qquad h_3(x,\bar x, z, \bar y)=o(\hat\gamma^{-k}),$$
such that $T_{k,m}(x,y,z)=(\bar x,\bar y,\bar z)$ if and only if
\begin{equation}\label{eq:maps:F_k,m_cross}
\begin{aligned}
\bar x -x^+ &=b\gamma^m\hat\mu - bu^- + ab\lambda^k\gamma^m x  + b_{13}\bar y -  \hat h_{01} (\bar x - x^+,\bar y) 
+  \gamma^m h_1(x,\bar x, z, \bar y),\\
\bar z - z^+ &=b_{31}b^{-1}(x-x^+-b_{13}\bar y)+b_{33}\bar y+\hat h_{02}(\bar x-x^+,\bar y)   + h_2(x,\bar x,z,\bar y),\\
y&=h_3(x,\bar x, z, \bar y).
\end{aligned}
\end{equation}
After the coordinate transformation
\begin{equation}\label{eq:maps:coortransform1}
\begin{array}{l}
X=x-x^+- b_{13}y, \qquad Y=y,\\
Z=z-z^+ -b_{31}b^{-1}(x-x^+-b_{13}y)-b_{33}y  -\hat h_{02}(x-x^+,y)
\end{array}
\end{equation}
the map $T_{k,m}$ assumes the form
\begin{equation}\label{eq:maps:T_k,m_cross}
\begin{array}{l}
\bar X  = b \gamma^m\hat\mu+ a b \lambda^k\gamma^m x^+ - b u^- + a b \lambda^k\gamma^m X  + 
\hat\phi_0(\bar X, \bar Y) + 
\gamma^m\hat\phi_1(X,\bar X,  \bar{Y},Z),\\[5pt]
Y = \hat\phi_2(X,\bar X, \bar Y, Z), \qquad
\bar Z = \hat\phi_3(X, \bar X, \bar Y, Z),
\end{array}
\end{equation}
where
\begin{equation}\label{oriph0}
\hat\phi_0=O(\bar X^2+\bar Y^2),
\end{equation}
i.e., $\hat\phi_0$ vanishes at $(\bar X,\bar Y)=0$ along with the first derivatives, and
\begin{equation}\label{eq:maps:nonlinearterms_original}
\|\hat\phi_1\|_{C^1}=o(\lambda^{k})+o(\gamma^{-m}), \qquad \|\hat\phi_2\|_{C^1} =o(\hat\gamma^{-k}), \qquad 
\|\hat\phi_3\|_{C^1}=o(\hat\lambda^{m}).
\end{equation}
By construction, these estimates are uniform for all systems $C^2$-close to $f$, and when the system is at least $C^3$-smooth, the same estimates are also true for the first derivatives with respect to parameters $\eps$. We note also that in the coordinates (\ref{eq:maps:coortransform1}) we have the local
stable manifold of $O_1$ and the image of the local unstable manifold of $O_2$ straightened:
\begin{equation}\label{straightyz}
W^s_{loc}(O_1): \{ Y=0\} \;\;\mbox{ and }\;\; F_{21}(W^u_{loc}(O_2)):\{Z=0\}
\end{equation}
(one can see this from that the equation for $F_{21}(W^u_{loc}(O_2))$ in the $(x,y,z)$ coordinates
can be obtained by taking the limit $m\to\infty$ in (\ref{eq:maps:221})).

We stress that the first-return map must take points from a small neighborhood of $M_1^+$ in $U_{01}$ to a small neighborhood of $M_1^+$, which corresponds to small values of 
$(X,Y,Z)$ and $(\bar X,\bar Y,\bar Z)$. We further denote this neighborhood where we want the first-return maps to be defined by
\begin{equation}\label{eq:maps:T_k,m_cross_domain}
\Pi=[-\delta,\delta]\times[-\delta,\delta]^{d_1}\times [-\delta,\delta]^{d-d_1-1}.
\end{equation}
For $(X,\bar Y,Z)\in \Pi$, we can assure that $Y$ and $\bar Z$ are small in (\ref{eq:maps:T_k,m_cross}) by taking $k$ and $m$ sufficiently large. 
However, to have $\bar X$ small (i.e. $\bar X=O(\delta)$)
one needs some additional restriction on possible values of $k$ and $m$. 

In the next section we will consider the first-return maps near the unperturbed cycles $\Gamma$, which corresponds to $\mu=\hat\mu=0$ in (\ref{eq:maps:T_k,m_cross})
(see \eqref{mmoh}). In this case, the $O(\delta)$ smallness of both $\bar X$ and $X$ in the first equation of (\ref{eq:maps:T_k,m_cross}) implies that the map $T_{k,m}$ at $\mu=0$ acts from $\Pi$ to an $O(\delta)$-neighborhood of $M_1^+$ in $U_{01}$ only when there is a certain balance between $k$ and $m$, namely
\begin{equation}\label{eq:km}
ab\lambda^k\gamma^m=\alpha+O(\delta),
\end{equation}
where $\alpha=bu^-/x^+\neq 0$ is the quantity introduced in condition \eqref{eq:intro:4}. 

Thus, $\lambda^k\gamma^m$ must be uniformly bounded in this case, and hence $\hat\phi_1=o(\gamma^{-m})$ in (\ref{eq:maps:nonlinearterms_original}). Consequently, the derivative of the right-hand side of the first equation of \eqref{eq:maps:T_k,m_cross} with respect to $\bar X$ is of order $O(\delta)+o(1)_{k,m\to\infty}$, so $\bar X$ can be expressed as a function of 
$( X, \bar Y, Z)$. Therefore, for sufficiently large $k,m$ such that \eqref{eq:km} is satisfied, formula \eqref{eq:maps:T_k,m_cross} for $T_{k,m}$ at $\mu=0$ implies that for a point $(X,Y,Z)\in\Pi$ we have $(\bar X,\bar Y,\bar Z)=T_{k,m}(X,Y,Z)$ if and only if the points are related by the cross-map $T^\times_{k,m}:(X,\bar Y, Z)\mapsto (\bar X,Y,\bar Z)$ given by
\begin{equation}\label{eq:maps:T_k,m_cross_0mu}
\begin{aligned}
\bar X &=a b \lambda^k\gamma^m x^+ - b u^- + a b \lambda^k\gamma^m X+ \phi_1(X,  \bar{Y},Z),\\
&=a b \lambda^k\gamma^m x^+ - b u^- + (\alpha+O(\delta)) X  + \phi_1(X,  \bar{Y},Z),\\
Y&=\phi_2(X,  \bar{Y},Z),\qquad
\bar Z = \phi_3(X,  \bar{Y},Z),
\end{aligned}
\end{equation}
where 
\begin{equation}\label{eq:maps:nonlinearterms_0mu}
\begin{array}{ll}
\phi_1=O(\delta^2)+o(1)_{k,m\to\infty},\qquad & \dfrac{\partial\phi_1}{\partial(X,\bar Y,Z)}= O(\delta)+o(1)_{k,m\to\infty},\\
 \|\phi_2\|_{C^1} =o(\hat\gamma^{-k}), \qquad & \|\phi_3\|_{C^1}=o(\hat\lambda^{m}).
\end{array}
\end{equation}

The following result characterizes the action of the derivative $\D T_{k,m}$ of the maps $T_{k,m}$ at $\mu=0$. We will further use the notation $(\Delta X,\Delta Y,\Delta Z)$ for vectors in the tangent space to $\Pi$.

\begin{lem}\label{lem:conefields} Let $\mu=0$.
Given any positive $K<1$, for all sufficiently small $\delta$ and large $(k,m)$ satisfying \eqref{eq:km}, the cone fields on $\Pi$
\begin{eqnarray}
&&\mathcal{C}^{cu} = \{(\Delta X,\Delta Y,\Delta Z): \|\Delta Z\| \leq K (|\Delta X|+\|\Delta Y\|)\},\label{eq:conefields:cu}\\
&&\mathcal{C}^{uu} = \{(\Delta X,\Delta Y,\Delta Z): \max\{|\Delta X|,\|\Delta Z\|\}\leq K  \|\Delta Y\|\},\label{eq:conefields:u}
\end{eqnarray}
are forward-invariant in the sense that if a point $M\in \Pi$ has its image $\bar M=T_{k,m}(M)$ in $\Pi$, then the cone at $M$ 
is mapped into the cone at $\bar M$ by $\D T_{k,m}$; and the cone fields
\begin{eqnarray}
&&\mathcal{C}^{cs} = \{(\Delta X,\Delta Y,\Delta Z): \|\Delta Y\|\leq K (|\Delta X|+\|\Delta Z\|)\},\label{eq:conefields:s}\\
&&\mathcal{C}^{ss}=\{(\Delta X,\Delta Y,\Delta Z):\max\{|\Delta X|,\|\Delta Y\|\}\leq K  \|\Delta Z\|\},\label{eq:conefields:ss}
\end{eqnarray}
are backward-invariant in the sense that if a point $\bar{M}\in\Pi$ has its preimage $M=T^{-1}_{k,m}(\bar M)$ in $\Pi$, 
then the cone at $\bar{M}$ is mapped into the cone at $M$ by $\D T^{-1}_{k,m}$. Moreover, vectors in $\mathcal{C}^{uu}$ and, if $|\alpha|>1$, also in $\mathcal{C}^{cu}$ are expanded by $\D T_{k,m}$; vectors in $\mathcal{C}^{ss}$ and, if $|\alpha|<1$, also in $\mathcal{C}^{cs}$ are contracted by $\D T_{k,m}$.
\end{lem}

\noindent{\it Proof.} Let us establish the backward invariance of $\mathcal{C}^{ss}$ and $\mathcal{C}^{cs}$. Take $\bar{M}=(\bar X, \bar Y, \bar Z)\in\Pi$ such that $T_{k,m}^{-1}(\bar{M})=M\in \Pi$ for some $(k,m)$. Take a vector $(\Delta \bar{X},\Delta \bar{Y},\Delta \bar{Z})$
in the tangent space at the point $\bar{M}$. 
Let $(\Delta X,\Delta Y,\Delta Z)=\D T^{-1}_{k,m} (\Delta \bar X,\Delta \bar Y,\Delta \bar Z)$.

It follows from \eqref{eq:maps:T_k,m_cross_0mu} and \eqref{eq:maps:nonlinearterms_0mu} that
\begin{equation}\label{eq:conefields:00}
\begin{array}{l}
\Delta \bar{X}= (\alpha+O(\delta)+o(1)_{k,m\to\infty}) \Delta X +  (O(\delta)+o(1)_{k,m\to\infty})\Delta\bar{Y} + (O(\delta)+o(1)_{k,m\to\infty})\Delta Z,\\
\Delta Y= o(\hat{\gamma}^{-k}) \Delta X +o(\hat{\gamma}^{-k})\Delta\bar{Y} + o(\hat{\gamma}^{-k})\Delta Z,\\
\Delta \bar{Z}= o(\hat{\lambda}^{m}) \Delta X +o(\hat{\lambda}^{m})\Delta\bar{Y} + o(\hat{\lambda}^{m})\Delta Z.
\end{array}
\end{equation}
Thus, there exists a constant $C$ such that if $k$ and $m$ are large enough, then
\begin{align}
(|\alpha| - C\delta) |\Delta X| &\leq |\Delta \bar X| +  C \delta \|\Delta\bar Y\| + C\delta \|\Delta Z\|,\label{eq:conefields:3}\\
|\Delta \bar X|  &\leq (|\alpha| + C\delta) |\Delta X| +  C \delta \|\Delta\bar Y\| + C\delta \|\Delta Z\|,\label{eq:conefields:3b}\\
\|\Delta Y\|&= o(\hat{\gamma}^{-k}) (|\Delta X| + \|\Delta\bar Y\| + \|\Delta Z\|),\label{eq:conefields:4}\\
\|\Delta \bar Z\|&= o(\hat{\lambda}^{m}) (|\Delta X| + \|\Delta\bar Y\| + \|\Delta Z\|).\label{eq:conefields:5}
\end{align}

Let $(\Delta \bar{X},\Delta \bar{Y},\Delta \bar{Z})\in \mathcal{C}^{ss}$. Since $K<1$,
we have
$$\max\{|\Delta \bar X|,\|\Delta \bar Y\|\}\leq \|\Delta \bar Z\|.$$
Now, it follows from \eqref{eq:conefields:5} that
\begin{equation}\label{contrcon}
\|\Delta \bar Z\| =o(\hat{\lambda}^{m})(|\Delta X|+\|\Delta Z\|)
\end{equation}
and, hence,
$$|\Delta \bar X|+\|\Delta \bar Y\| =o(\hat{\lambda}^{m})(|\Delta X|+\|\Delta Z\|).$$
We substitute these estimates into  \eqref{eq:conefields:3} and  \eqref{eq:conefields:4} and obtain
$$|\Delta X| = O(\delta) \|\Delta Z\| \qquad {\rm and } \qquad \|\Delta Y\|= o(\hat{\gamma}^{-k}) \|\Delta Z\|,$$
i.e., for any fixed choice of the constant $K$, 
if $k$ and $m$ are large enough and $\delta$ is small enough, the vector $(\Delta X,\Delta Y,\Delta Z)$ lies in $\mathcal{C}^{ss}$ at the point $M$, as required.
Equation (\ref{contrcon}) implies the contraction in $\mathcal{C}^{ss}$ if $m$ is large enough.

Similar arguments are applied when $(\Delta \bar{X},\Delta \bar{Y},\Delta \bar{Z})\in \mathcal{C}^{cs}$. Here, we have 
$$\|\Delta \bar Y\|\leq |\Delta \bar X|+\|\Delta \bar Z\|.$$
Substituting this into \eqref{eq:conefields:3b} and \eqref{eq:conefields:5} gives
\begin{equation}\label{contrcon2}
\begin{array}{l}
 |\Delta \bar X|  \leq (|\alpha| + O(\delta)) |\Delta X| +  O(\delta) \|\Delta Z\|,\\
\|\Delta \bar Z\| =o(\hat{\lambda}^{m})(|\Delta X|+\|\Delta Z\|),
\end{array}
\end{equation}
for sufficiently small $\delta$ and sufficiently large $m$, and, hence,
$$\|\Delta \bar Y\|\leq (|\alpha| + 1) |\Delta X| +  \|\Delta Z\|.$$
We substitute the last estimate into \eqref{eq:conefields:4} and obtain
$$\|\Delta Y\|= o(\hat{\gamma}^{-k}) (|\Delta X|+\|\Delta Z\|),$$
i.e., for any fixed choice of the constant $K$, 
if $k$ is large enough, the vector $(\Delta X,\Delta Y,\Delta Z)$ lies in $\mathcal{C}^{cs}$ at the point $M$, as required.
Equation (\ref{contrcon2}) implies the contraction in $\mathcal{C}^{cs}$ if $|\alpha|<1$ and $m$ is large enough and $\delta$ is
small enough.

The proof of the forward invariance of $\mathcal{C}^{uu}$ and $\mathcal{C}^{cu}$ is done in the same way, as everything is symmetric here
with respect to the change of $T_{k,m}$ to $T_{k,m}^{-1}$. \qed

\section{Blenders near type-I heterodimensional cycles. Proof of Theorem \ref{thm:blender_in_type1}}\label{sec:blendersintype1}
In this section we obtain a detailed description of the blenders that may appear near a  type-I heterodimensional cycle, which is summarized in Propositions  \ref{prop:cublendertype1} and \ref{prop:csblendertype1} below. These two propositions together imply Theorem \ref{thm:blender_in_type1}. Throughout this section, we do not perturb the system $f$ in this theorem, i.e., the fragile heteroclinic is not split, so $\hat\mu=0$ in \eqref{eq:maps:112}. We consider non-degenerate type-I cycles, which means that $|\alpha|=|bu^-/x^+|\neq 1$ and $ax^+ u^->0$ in \eqref{eq:maps:T_k,m_cross}. We also assume that
$\theta=-\ln|\lambda|/\ln|\gamma|$ is irrational.

\subsection{Blenders with activating pairs}

\begin{defi}[Proper crossing]\label{defi:surface}
Consider a cube $Q=\{(x_1,\dots,x_d)\mid x_i\in I_i\}\subset \mathbb{R}^d,$ 
where $I_i$ are closed intervals in $\mathbb{R}$. A $k$-dimensional disc $S$ is said to {\em cross} $Q$ if the intersection $S\cap Q$ is given by $(x_{i_{k+1}},\dots,x_{i_d})=s(x_{i_1},\dots,x_{i_k})$, where  $(i_1,\dots,i_d)$ is some permutation of $(1,\dots,d)$, and $s$ is a smooth function defined on $I_{i_1}\times \dots\times I_{i_k}$. It crosses $Q$ {\em properly} with respect to a cone field $\mathcal{C}$ defined on $Q$ if the tangent spaces of $S\cap Q$ lie in $\mathcal{C}$.
\end{defi}
For example, a disc $S$ crossing $\Pi$ properly with respect to $\mathcal{C}^{ss}$ in Lemma \ref{lem:conefields} is the graph of some smooth function $(X,Y)=s(Z)$, which is defined on $[-\delta,\delta]^{d-d_1-1}$ and whose tangent spaces lie in $\mathcal{C}^{ss}$.

For an invariant cone field $\mathcal{C}$, denote by $\dim (\mathcal{C})$ the largest possible dimension of a linear subspace of the tangent space that can be contained in $\mathcal{C}$ (at each point where $\mathcal{C}$ is defined).

\begin{defi}[Activating pair]\label{defi:actcondition}
A pair $(Q ,\mathcal{C})$ consisting of a cube and a cone field is called an {\em activating pair} for a cs-blender $\Lambda$ if $\dim (\mathcal{C})=\dim W^s(\Lambda)-1$ and if any disc $S$ crossing $Q$ properly with respect to $\mathcal{C}$ intersects $W^u(\Lambda)$. Similarly, a pair $(Q ,\mathcal{C})$ is an activating pair for a cu-blender $\Lambda$ if $\dim (\mathcal{C})=\dim W^u(\Lambda)-1$ and any disc $S$ crossing $Q$ properly with respect to $\mathcal{C}$ intersects $W^s(\Lambda)$. The cube $Q$ is called an {\em activating domain}.
\end{defi}

\begin{rem}\label{rem:actpair}
Every disc crossing $Q$ properly with respect to $\mathcal{C}$ belongs to the set $\mathcal{D}^{ss}$ of Definition \ref{defi:blender_old} in the case of a cs-blender, or to the set $\mathcal{D}^{uu}$ in the case of a cu-blender. Thus, the activating pair specifies the position of robust non-transverse intersections. The cone field $\mathcal{C}$ is the field $\mathcal{C}^{ss}$ (for the cs-blender) or $\mathcal{C}^{uu}$ (for the cu-blender) from the definition of the standard blender in the Appendix.
\end{rem}

Note that any disc $C^1$-close to $S$ in the above definition is also a properly crossing disc, and hence intersects the invariant manifold of the blender. In particular, if $\Lambda$ is a cs-blender and $S$ is a piece of the stable manifold $W^s(\Lambda')$ of some hyperbolic basic set $\Lambda'$, then we obtain robust non-transverse intersection of $W^u(\Lambda)$ with $W^s(\Lambda')$. If at the same time we have $W^s(\Lambda) \pitchfork W^u(\Lambda')\neq\emptyset$, then robust heterodimensional dynamics emerge. Similar arguments hold in the case of a cu-blender. 

Recall that $\Pi$ defined in \eqref{eq:maps:T_k,m_cross_domain} is a cube around the robust heteroclinic  point $M^+_1$.
\begin{prop}\label{prop:cublendertype1} Let $\Gamma$ be a non-degenerate cycle of type I, with irrational $\theta$.
If $|\alpha|<1$, then there exists arbitrarily close to $\Gamma$ an index-$d_1$ cs-blender $\Lambda^{cs}$ with an activating pair $(\Pi' ,\mathcal{C}^{ss})$, where $\Pi'$ is a subcube of $\Pi$ with the center at $M^+_1$, and 
\begin{equation}\label{eq:sscone}
\mathcal{C}^{ss}=\{(\Delta X,\Delta Y,\Delta Z): \max(|\Delta X|,\|\Delta Y\|)< \dfrac{q|\alpha|}{4} \|\Delta Z\|\}.
\end{equation}
\end{prop}

When $|\alpha|>1$, one immediately gets a cu-blender of index $d_1+1$ by applying the above proposition due to the system obtained by the time 
reversal. However, the corresponding activating domain would lie near $O_2$ but not $O_1$ (the time-reversal also interchanges $O_1$ and $O_2$).
This would not match with the generalization, Proposition \ref{prop:blender}, which is the main tool used for showing the existence of blenders in the saddle-focus and double-focus cases. So, we prove the following propoasition about the cu-blender with an activating domain $\Pi'$ near $O_1$ as in the case
$|\alpha|<1$.

\begin{prop}\label{prop:csblendertype1} Let $\Gamma$ be a non-degenerate cycle of type I, with irrational $\theta$.
If $|\alpha|>1$, then there exists, arbitrarily close to $\Gamma$, an index-$(d_1+1)$ cu-blender $\Lambda^{cu}$ with an activating pair $(\Pi' ,\mathcal{C}^{uu})$, where $\Pi'$ is a subcube of $\Pi$ with the center at $M^+_1$, and
\begin{equation}\label{eq:uucone} 
\mathcal{C}^{uu} = \{(\Delta X,\Delta Y,\Delta Z): \max\{|\Delta X|,\|\Delta Z\|\}\leq \dfrac{q|\alpha|^{-1}}{4}  \|\Delta Y\|\}.
\end{equation}
\end{prop}

Below, in the proof of these propositions, we establish the existence of blenders in the sense of the ``operational'' Definition \ref{defi:blender_old}. However, the construction is explicit, and
the obtained blenders are standard in the sense of Definiton \ref{defi:blender}, as we show in Proposition {prop:stanble}.

\subsection{Finding cs-blenders when $|\alpha|<1$. Proof of Proposition \ref{prop:cublendertype1}}
Recall that the cube $\Pi$ defined in \eqref{eq:maps:T_k,m_cross_domain} is the region for which we consider the first-return maps $T_{k,m}$.
We start with finding hyperbolic sets for the collection of the return maps in $\Pi$, which can be candidates for blenders. Consider the set of pairs $(k,m)$ of sufficiently large integers:
\begin{equation}\label{eq:L_N}
\mathcal P_N=\{(k,m): k>N, m>N\,\,\mbox{and}\,\, 
|a b\lambda^{k}\gamma^{m}x^+ - b u^-|\leq \dfrac{2}{3}(1-|\alpha|)\delta\}.
\end{equation}
Note that this set is non-empty for any $N$ and $\delta>0$: because $\theta=-\ln|\lambda|/\ln|\gamma|$ is irrational, the set
$\lambda^{k}\gamma^{m}$ is dense among positive reals, hence, since $ab x^+$ and $bu^-$ have the same sign by the assumption of the theorem (we consider the cycles of type I), $a b\lambda^{k}\gamma^{m}x^+$ can be made as close as we want to $b u^-$ for arbitrarily large $k$ and $m$.

By construction, the estimate (\ref{eq:km}) is satisfied for every $(k,m)\in \mathcal P_N$. Therefore, for every such $(k,m)$
the relation between the coordinates $(X,Y,Z)\in \Pi$ of a point in the domain of definition of the first-return map $T_{k,m}$
and the coordinates $(\bar X,\bar Y,\bar Z)$ of its image by $T_{k,m}$ is given by \eqref{eq:maps:T_k,m_cross_0mu}.

Let us denote by $T^\times_{k,m}:(X,\bar Y, Z)\mapsto (\bar X,Y,\bar Z)$ the cross-map \eqref{eq:maps:T_k,m_cross_0mu}. Since $|\alpha|<1$, comparing \eqref{eq:L_N} with the $X$-equation in \eqref{eq:maps:T_k,m_cross_0mu}, one sees that for every $(k,m)\in \mathcal P_N$, with $N$ sufficiently large and $\delta$ sufficiently small, we have $\bar X\in [-\delta,\delta]$. On the other hand, it is obvious that $(Y,\bar Z)$-coordinates lie in $[-\delta,\delta]^{d-1}$ for all sufficiently large $k$ and $m$, due to the strong contraction given by  \eqref{eq:maps:T_k,m_cross_0mu}. Thus, for any $(k,m)\in \mathcal P_N$ with sufficiently large $N$ and $m$, the cross-map $T^\times_{k,m}$ satisfies
$$T^\times_{k,m}(\Pi)\subset \Pi.$$
It also follows from $|\alpha|<1$ and \eqref{eq:maps:T_k,m_cross_0mu} that the cross-map is contracting on $\Pi$, i.e.,
$$\left\|\frac{\partial (\bar X, Y, \bar Z)}{\partial (X, \bar Y, Z)}\right\| < 1.$$
By Shilnikov lemma on the fixed  point in a direct product of metric spaces (Theorem 6.2 in \citep{Sh67}), these two facts 
immediately imply
\begin{lem}\label{lem:hypset00}
For any sequence $\{(k_s,m_s)\}_{s\in \mathbb{Z}}$ of pairs $(k,m)$ from $\mathcal{P}_N$, 
there exists a unique sequence of
points $\{M_s=(X_s,Y_s,Z_s)\}_{s\in \mathbb{Z}}$ in $\Pi$ such that $M_{s+1}=T_{k_s,m_s} M_s$ for every ${s\in \mathbb{Z}}$.
\end{lem}
\begin{proof}
By definition of the cross-map, a sequence of points $\{M_s=(X_s,Y_s,Z_s)\}$ in $\Pi$ satisfies
\begin{equation}\label{eq:coding}
M_{s+1}=T_{k_s,m_s} M_s
\end{equation}
if and only if
$$
(X_{s+1},Y_{s},Z_{s+1})=T^\times_{k_s,m_s} (X_s, Y_{s+1}, Z_s).
$$
Thus, the sought sequence $\{(X_s,Y_s,Z_s)\}_{s\in \mathbb{Z}}$ is a fixed point of the map
$\{(X_s,Y_s,Z_s)\}_{s\in \mathbb{Z}}\mapsto \{(\tilde X_s,\tilde Y_s,\tilde Z_s)\}_{s\in \mathbb{Z}}$
(acting on the space of sequences of points in $\Pi$) which is defined by the rule
$$(\tilde X_{s+1},\tilde Y_{s},\tilde Z_{s+1})=T^\times_{k_s,m_s} (X_s, Y_{s+1}, Z_s).$$
Obviously, this map is a contraction (because each of the maps $T^\times_{k_s,m_s}$ is a contraction on $\Pi$),
so the fixed point indeed exists and is unique.
\end{proof}

We call the sequence  $\{(k_s,m_s)\}_{s\in \mathbb{Z}}$ the {\em coding} of the point $M_0$. Lemma \ref{lem:hypset00}
establishes the existence of an invariant set $\Lambda$ whose intersection with $\Pi$ is in {\em one-to-one correspondence} with the set of all codings formed
from the pairs $(k,m)\in\mathcal{P}_N$. Namely, a point $M_0$ lies in $\Lambda\cap \Pi$ if and only if there is a coding $\{(k_s,m_s)\in \mathcal{P}_N\}_{s\in \mathbb{Z}}$ such that the points defined by
$$
M_s=
\begin{cases}
T_{k_{s-1},m_{s-1}}\circ\dots\circ T_{k_0,m_0}M_0 \mbox{ for } s>0\\
T^{-1}_{k_{-s},m_{-s}}\circ\dots\circ T^{-1}_{k_{-1},m_{-1}}M_0 \mbox{ for } s<0
\end{cases}
$$
all lie in $\Pi$. The set $\Lambda$ is  the union of the orbits of all such points $M_0\in\Pi$ by the  system $f$. 

We are interested in a certain closed subset ${\Lambda}_{\mathcal{J}}\subset\Lambda$ which corresponds to the set of codings with pairs $(k,m)$ taken from some finite subset $\mathcal{J}\subset \mathcal{P}_N$ (so the dynamics on ${\Lambda}_{\mathcal{J}}$ correspond to the full shift with this finite set of symbols). It  follows from Lemma \ref{lem:conefields} that this set is a locally maximal, transitive hyperbolic set and, also, partially-hyperbolic with a 1-dimensional central (weakly-hyperbolic) direction. In summary, we have
\begin{lem}\label{lem:hypset2}
Let $|\alpha|<1$ and $\mathcal{J}$ be any finite subset of $\mathcal{P}_N$. Then, there exists a hyperbolic basic set ${\Lambda}_{\mathcal{J}}$ of index $d_1$ near the heterodimensional cycle $\Gamma$ such that it is in one-to-one correspondence with the set of codings $\{(k_s,m_s)\in\mathcal{J}\}_{s\in\mathbb Z}$.
\end{lem}

\begin{rem}\label{rem:strips}
The intersection $\Lambda_{\mathcal{J}}\cap\Pi$ is located in a finite union of ``horizontal strips'' in $\Pi$. Indeed, one sees from \eqref{eq:maps:T_k,m_cross_0mu} that the domains of $T_{k,m}$ with $(k,m)\in \mathcal P_N$ are given by the strips:
\begin{equation}\label{eq:strips}
\sigma_{k,m}=\{(X,Y,Z)\mid X\in [-\delta,\delta], Y\in \phi_2(X, [-\delta,\delta]^{d_1},Z),Z\in [-\delta,\delta]^{d-d_1-1}  \}.
\end{equation}
Since $\phi_2=o(\hat\gamma^{-k})$, these strips accumulate on $\{Y=0\}$ and strips corresponding to different $(k,m)$ are disjoint. By construction, we have 
\begin{equation}\label{eq:sigmaJ}
\Lambda_{\mathcal{J}}\cap\Pi\subset \bigcup_{(k,m)\in\mathcal J}\sigma_{k,m}=:\Sigma_\mathcal{J}.
\end{equation}
Every point whose backward orbit stays in $\Sigma_\mathcal{J}$ lies in $W^u(\Lambda_\mathcal{J})$, and every point whose forward orbit stays in $\Sigma_\mathcal{J}$ lies in $W^s(\Lambda_\mathcal{J})$. This property holds for any $C^1$-close system since the hyperbolic set persists under $C^1$-small perturbations.
\end{rem}

\begin{rem}\label{rem:transintersection} For every point $M\in \Lambda_{\mathcal{J}}$, we can define
its local stable manifold as a connected piece of  $W^s(\Lambda_{\mathcal{J}})\cap \Pi$ through $M$,
and the local unstable manifold as a connected piece of  $W^u(\Lambda_{\mathcal{J}})\cap \Pi$ through $M$. If $M\in \sigma_{k,m}$, then $W^s_{loc}(M)\subset \sigma_{k,m}$. Moreover, since the tangent
to $W^s_{loc}(M)$ at any point lies in the stable cone $\mathcal{C}^{cs}$, it follows that 
$W^s_{loc}(M)$ is a ``horizontal disc'' of the form $Y=\xi^s_M(X,Z)$, where the smooth 
function $\xi^s_M$ is defined for all $X\in [-\delta,\delta], Z\in [-\delta,\delta]^{d-d_1-1}$. Similarly, $W^u_{loc}(M)$
is a ``vertical disc'' of the form $(X,Z)=\xi^u_M(Y)$, where $\xi^u_M$ is defined for all 
$Y\in [-\delta,\delta]^{d_1}$. In particular, it follows from (\ref{straightyz}) that for each $M$ there exist transverse intersections
$$W^u_{loc}(M)\cap W^s_{loc} (O_1)\neq\emptyset\quad\mbox{and}\quad
W^s_{loc}(M) \cap F_{21}(W^u_{loc}(O_2))\neq\emptyset.$$
\end{rem}

Now let us find a subset $\mathcal{J}$ such that the corresponding hyperbolic set $\Lambda_\mathcal{J}$ is a cs-blender. We will define a cube $\Pi'\subset\Pi$ and take discs which cross $\Pi'$ properly with respect to the cone field $\mathcal{C}^{ss}$ (see Definition \ref{defi:surface}) given by Lemma \ref{lem:conefields} (with an appropriate $K$ in its definition to be determined). We will show that any such disc $S$ intersects the unstable manifold $W^u(\Lambda_\mathcal{J})$ (and this property persists at $C^1$-small perturbations of the system).

Recall that by Remark \ref{rem:strips}, any point whose entire backward orbit by $T_{k,m}$ stays in $\Sigma_{\mathcal{J}}$ belongs to $W^u(\Lambda_\mathcal{J})$. Therefore, we will have the desired property that 
$S\cap W^u(\Lambda_\mathcal{J})\neq \emptyset$ by showing below that any disc $S$ crossing
$\Pi'$ properly with respect to $\mathcal{C}^{ss}$ contains a points whose entire backward orbit stays in $\Sigma_{\mathcal{J}}$.

Take any positive $q<(1-|\alpha|)/2$, and define
\begin{equation}\label{eq:delta'}
\delta'=q\delta.
\end{equation}

\begin{lem}\label{lem:crossingsurfaces}
There exists a finite subset $\mathcal{J_{\delta'}}\subset \mathcal{P}_{N}$ such that, if a $(d-d_1-1)$-dimensional disc $S$ crosses 
\begin{equation}\label{eq:actdomian}
\Pi': =\Pi \cap\{-\delta'\leq X \leq\delta'\}
\end{equation}
properly with respect to the cone field $\mathcal{C}^{ss}$ defined in \eqref{eq:sscone},
then one can find a pair $(k,m)\in\mathcal{J_{\delta'}}$ such that the preimage $T^{-1}_{k,m}(S)$ contains a disc crossing $\Pi' $ properly with respect to $\mathcal{C}^{ss}$.
\end{lem}

\begin{figure}[!h]
\begin{center}
\includegraphics[width=0.4\textwidth]{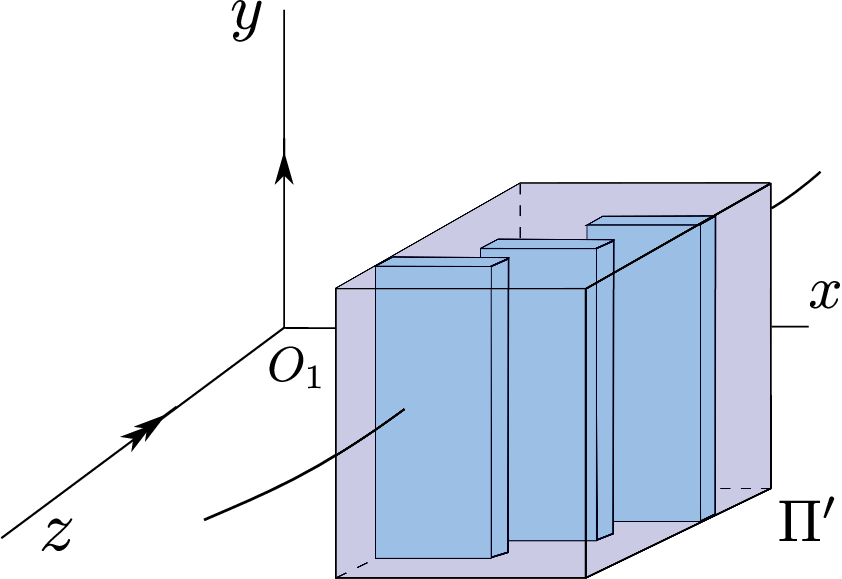}
\end{center}
\caption{The blue regions depict intersections $\Pi'\cap T_{k,m}(\Pi') $ with $(k,m)\in \mathcal{J}_{\delta'}$. For any disc (or curve in this three dimensional figure) crossing $\Pi'$ properly, it must intersect at least one of the three blue regions. Hence its preimage contains a piece crossing $\Pi' $.}
\label{fig:cubecovering}
\end{figure}

The proof of the lemma is postponed until the end of this section. See Figure \ref{fig:cubecovering} for an illustration. The lemma shows that given any properly crossing disc $S$, there is a sequence of discs defined by $S_{i+1}=T^{-1}_{k_i,m_i}(S_i)\cap \Pi'$ with $S_0:=S$ and some sequence $(k_i,m_i)\in \mathcal{J}_{\delta'}$. By construction, this gives a sequence of nested closed sets
$$\hat{S}_i:=T_{k_0,m_0}\circ\dots\circ T_{k_{i-1},m_{i-1}}(S_i)\subset S.$$
The intersection of all $\hat{S}_i$ contains a point $M$ whose backward orbit stays in $\Sigma_{J_{\delta'}}$ defined by \eqref{eq:sigmaJ} (in fact, this point is unique due to the contraction in $Z$-directions.) Hence, we have $M\in W^u(\Lambda_{\mathcal{J}_{\delta'}})$, where $\Lambda_{\mathcal{J}_{\delta'}}$ is the hyperbolic set given by Lemma \ref{lem:hypset2} associated with the set $\mathcal{J}_{\delta'}$. Obviously, Lemma \ref{lem:crossingsurfaces} holds for any $C^1$-close system for the same set $\mathcal{J}_{\delta'}$. Therefore, we  obtain a cs-blender of Definition \ref{defi:blender_old}. 

Since the cone field $\mathcal{C}^{ss}$ is backward-invariant by Lemma \ref{lem:conefields}, the preimage of any disc proper with respect to $\mathcal{C}^{ss}$ is still proper. Hence, in order to prove Lemma \ref{lem:crossingsurfaces} we only need to show that any properly crossing disc has a preimage which crosses $\Pi'$. The key procedure in controlling the preimage of the disc $S$ is to show that the $X$-coordinates of its points always lie in the domain $\sigma_{k,m}$ of some return map $T_{k,m}$ (the so-called ``covering property'', as described in the definition of a standard blender in the Appendix, see more discussion in \citep{BDV,NP12,BKR14}). This amounts to prove the Lemma \ref{lem:I_k,m}, which we formulate below.

Denote the length of an interval $E$ by $|E|$, and the affine part of the $X$-equation in \eqref{eq:maps:T_k,m_cross_0mu} by
\begin{equation}\label{eq:I_k,m:linearX}
R_{k,m}(X):=a b \lambda^k\gamma^m x^+ - b u^- + a b \lambda^k\gamma^m X.
\end{equation}
Recall that $\delta$ is the size of the domain of the return maps $T_{k,m}$ as in \eqref{eq:strips}, and $\delta'$ is the size of the smaller cube $\Pi'\subset \Pi$, which depends on $\delta$ as given by \eqref{eq:delta'}.

\begin{lem}\label{lem:I_k,m}
Let $|\alpha|<1$. There exists, for all sufficiently large $N$ and sufficiently small $\delta$, a finite subset $\mathcal{J_{\delta'}}=\{(k_j,m_j)\}^n_{j=1}\subset \mathcal{P}_{N}$ such that the intervals $E_j:=R_{k_j,m_j}([-\delta',\delta'])$ satisfy
\begin{equation}\label{eq:I_k,m:covering}
\bigcup^n_{j=1}int(E_j)\supset  [-\delta',\delta']
\end{equation}
and
\begin{equation}\label{eq:I_k,m:overlap}
|E_j\cap E_{j+1}|> \dfrac{\delta'|\alpha|}{2}, \quad j=1,\dots n-1.
\end{equation}
\end{lem}

\begin{proof}
One notes that since $\theta$ is irrational, the set of numbers $k\theta -m$ where $k$ and $m$ are even positive integers is dense in $\mathbb{R}$. Thus, given any $\rho \in \mathbb{R}$, we can find
a sequence $\{(k_i(\rho),m_i(\rho))\}_{i\in \mathbb{N}}$ of pairs of positive even integers such that $k_i(\rho),m_i(\rho)\to \infty$ and
\begin{equation*}\label{eq:log}
k_i(\rho)\theta-m_i(\rho)\to - \ln\left(\dfrac{b u^-+\delta\rho}{a b x^+}\right)\ln^{-1}|\gamma|
\end{equation*}
as $i\to\infty$ (the logarithm is defined for small $\delta$ since $ax^+u^->0$ and $b\neq 0$ by our assumptions). Since $k_i$ and $m_i$ are even,
we have $|\lambda|^{k_i}=\lambda^{k_i}$ and $|\gamma|^{m_i}=\gamma^{m_i}$ even when $\lambda$ or $\gamma$ are negative, so we obtain that
\begin{equation}\label{eq:hypset:1}
\dfrac{1}{\delta}(a b \lambda^{k_i(\rho)}\gamma^{m_i(\rho)}x^+ - b u^-) \to \rho
\end{equation}
as $i\to \infty$. It follows from \eqref{eq:L_N} that for any 
\begin{equation}\label{eq:I_hyp}
\rho\in\left(-\dfrac{2}{3}(1-|\alpha|), \dfrac{2}{3}(1-|\alpha|)\right)=:I_{hyp},
\end{equation}
we have
\begin{equation}\label{eq:L_Nrho}
\mathcal P_N(\rho):=\{(k_i(\rho),m_i(\rho)):k_i(\rho),m_i(\rho)>N\}\subset \mathcal P_N.
\end{equation} 
In what follows we construct $\mathcal J_{\delta'}$ by taking pairs $(k,m)$ from $\mathcal P_N(\rho)$ for a finite set of values of $\rho$.

Recall $\alpha=bu^-/x^+$. Note from \eqref{eq:I_k,m:linearX} and \eqref{eq:hypset:1} that, for any $(k_i(\rho),m_i(\rho))\in \mathcal{P}_N(\rho)$, we have 
$$R_{k_i(\rho),m_i(\rho)}(X)=\rho\delta + \left(\alpha+\dfrac{\rho\delta}{x^+}\right)X+o(1)_{i\to \infty}.$$
Consequently, the end points of $E_{k_i(\rho),m_i(\rho)}=R_{k_i(\rho),m_i(\rho)}([-\delta',\delta])$ are
$$a_{i}(\rho)=\rho\delta-\delta'\left(\alpha+\dfrac{\rho\delta}{x^+}\right)+o(1)_{i\to \infty}
\quad\mathrm{and}\quad
b_{i}(\rho)=\rho\delta+\delta'\left(\alpha+\dfrac{\rho\delta}{x^+}\right)+o(1)_{i\to \infty}.
$$
The length of each $E_{k_i(\rho),m_i(\rho)}$ satisfies
\begin{equation}\label{eq:I_k,m:length}
|E_{k_i(\rho),m_i(\rho)}|=|b_{i}-a_{i}|=2\delta'\left|\alpha+\dfrac{\rho\delta}{x^+}\right|+o(1)_{i\to \infty}>\delta'|\alpha|
\end{equation}
for all sufficiently large $i$ and sufficiently small $\delta$. 

We can now construct the desired set $\mathcal{J}_{\delta'}$ of pairs. Assume $\alpha>0$ for certainty. The case of negative $\alpha$ is completely parallel.

Let $n$ be the smallest integer not less than ${4}/{|\alpha|}+1$. Take
$$\rho_j=-\dfrac{\delta'}{\delta}+\dfrac{(j-1)2\delta'}{(n-1)\delta}\quad j=1,\dots, n.$$ The length of the intersections between two consecutive intervals $E_{k_i(\rho_j),m_i(\rho_j)}$ and $E_{k_i(\rho_{j+1}),m_i(\rho_{j+1})}$ are computed as 
\begin{equation*}
\begin{aligned}
b_i(\rho_j)-a_i(\rho_{j+1})
&=(\rho_j - \rho_{j+1})\delta + \delta'\left(2\alpha + \dfrac{(\rho_{j}+\rho_{j+1})\delta}{x^+}\right) + o(1)_{i\to\infty}\\
&=-\dfrac{2\delta'}{n-1} + \delta'(2\alpha+O(\delta')) + c_i,
\end{aligned}
\end{equation*}
where $c_i=o(1)_{i\to\infty}$. Let us fix a large $i=i^*$ such that for all $j=1,\dots,n-1$ and for all sufficiently small $\delta$ (hence sufficiently small $\delta'$ by \eqref{eq:delta'}), the last two terms in the above equation satisfy $\delta'(2\alpha+O(\delta')) + c_{i^*}>\delta'\alpha$. By the choice of $n$, this implies 
$$b_{i^*}(\rho_j)-a_{i^*}(\rho_{j+1})>\dfrac{\delta'\alpha}{2}.$$
Since by construction one has $-\delta'\in int(E_{k_{_{i^*}}(\rho_1),m_{_{i^*}}(\rho_1)})$ and $-\delta'\in int(E_{k_{i^*}(\rho_n),m_{i^*}(\rho_n)})$, it follows that $\mathcal{J}_{\delta'}=\{(k_j,m_j):=(k_{_{i^*}}(\rho_j),m_{_{i^*}}(\rho_j))\}_{j=1}^n$.
\end{proof}

\noindent{\it Proof of Lemma \ref{lem:crossingsurfaces}}.
As we explained before Lemma \ref{lem:I_k,m}, it suffices to show that for any proper disc $S$ there exists a pair $(k,m)\in \mathcal{J}_{\delta'}$
such that $T^{-1}_{k,m}(S)$ crosses $\Pi'$. As $S$ crosses
$\Pi'$, it is a graph of a smooth function $s=(s_X,s_Y):  [-\delta,\delta]^{d-d_1-1} \to [-\delta',\delta']\times
[-\delta,\delta]^{d_1}$. We need to find some $(k,m)\in\mathcal{J_{\delta'}}$ such that, for any $Z\in[-\delta,\delta]^{d-d_1-1}$, there exist $X\in[-\delta',\delta'],Y\in [-\delta,\delta]^{d_1}$ and $\bar{Z}\in [-\delta,\delta]^{d-d_1-1}$ which satisfy $(X,Y,Z)=T^{-1}_{k,m}(s_X(\bar{Z}),s_Y(\bar{Z}),\bar{Z})$. By formulas \eqref{eq:maps:T_k,m_cross_0mu} and \eqref{eq:I_k,m:linearX}, this is equivalent to solving the system of equations
\begin{eqnarray}
s_X(\bar{Z})&=& R_{k,m}(X) + {\phi}_1(X,s_Y(\bar{Z}),Z),\label{eq:xingsurfs:1}\\
Y&=&\phi_2(X,s_Y(\bar{Z}),Z),\label{eq:xingsurfs:2}\\
\bar{Z}&=& {\phi}_3(X,s_Y(\bar{Z}),Z).\label{eq:xingsurfs:3}
\end{eqnarray}

According to the estimates \eqref{eq:maps:nonlinearterms_0mu}, $Y$ lies in $[-\delta,\delta]^{d_1}$ for all sufficiently large $k$ as required; also, $\bar{Z}$ can be expressed from the last equation as a function 
of $(X,Z)$:
$$\bar{Z}=\tilde{\phi}_3(X,Z)$$
where $\|\tilde{\phi}_3\|_{C^1}=o(\hat\lambda^m)$. Hence, to solve the above system it suffices to find $X\in[-\delta',\delta']$ satisfying
\begin{equation}\label{eq:xingsurfs:4}
R_{k,m}(X) - s_X(\tilde{\phi}_3(X,Z)) + \tilde{\phi}_1(X,Z)=0,
\end{equation}
where
\begin{equation}\label{eq:xingsurfs:5}
\tilde{\phi}_1(X,Z)=\phi_1(X,s_Y(\tilde{\phi}_3(X,Z)),Z)=O(\delta^2)+o(1)_{k,m\to\infty}.
\end{equation}
We claim that one can choose $(k,m)\in \mathcal J_{\delta'}$ such that for any fixed $Z\in [-\delta,\delta]^{d-d_1-1}$, the left-hand side of \eqref{eq:xingsurfs:4} takes both positive and negative values when $X$ runs over $[-\delta',\delta']$. The lemma then follows by the intermediate value theorem.

Let us prove this claim. On one hand, the proper crossing with respect to $\mathcal{C}^{ss}$ means that the range of $s_X$ lies in $[-\delta',\delta']$, and, according to \eqref{eq:sscone}, the total change of $s_X$ is bounded by 
$$\max_{Z_1,Z_2\in[-\delta,\delta]^{d-d_1-1}} \{|s_X(Z_1)-s_X(Z_2)|\}<\dfrac{\delta'|\alpha|}{4}.$$
Since $\tilde{\phi}_3([-\delta',\delta'],Z)\subset [-\delta,\delta]^{d-d_1-1}$ for large $m$, we obtain that $s_X(\tilde{\phi}_3([-\delta',\delta'],Z))$ lies in some interval in $[-\delta',\delta']$ of length less than $\delta'|\alpha|/4$. 

On the other hand, Lemma \ref{lem:I_k,m} shows that the interiors of the intervals $R_{k,m}([-\delta',\delta'])$ with $(k,m)\in\mathcal{J}_{\delta'}$ cover the interval $[-\delta',\delta']$ with overlaps larger than $\delta'|\alpha|/2$. It follows that there exist pairs $(k,m)\in\mathcal{J}_{\delta'}$ with arbitrarily large $k,m$ such that the $\delta'|\alpha|/8$-neighbourhood of $s_X(\tilde{\phi}_3([-\delta',\delta'],Z))$ lies in $int(R_{k,m}([-\delta',\delta']))$ for all $Z\in [-\delta,\delta]^{d-d_1-1}$. Thus, the difference $R_{k,m}(X) - s_X(\tilde{\phi}_3(X,Z))$ changes its sign when $X$ runs the interval $[-\delta',\delta']$. In fact, there exists $C>0$ such that $R_{k,m} - s_X$ runs from $C\delta'$ to $-C\delta'$. Therefore, by \eqref{eq:xingsurfs:5}, the claim is proved when $\delta$ is sufficiently small and $k,m$ are sufficiently large.
\qed

The preceding lemma gives us a cs-blender of index $d_1$ with an activating pair $(\Pi' ,\mathcal{C}^{ss})$, and this blender can be taken arbitrarily close to the heterodimensional cycle by taking $N$ of $\mathcal{P}_N$ sufficiently large. Proposition \ref{prop:cublendertype1} follows immediately.

\subsection{Finding  cu-blenders when $|\alpha|>1$. Proof of Proposition \ref{prop:csblendertype1}}
As mentioned before, instead of using the symmetry of the problem to obtain cu-blenders in the $|\alpha|>1$ case,  we make computations directly for the maps $T^{-1}_{k,m}$ in order to  get an activating pair with the same activating domain $\Pi'$ defined in \eqref{eq:actdomian}.

\begin{proof}[Proof of Proposition \ref{prop:csblendertype1}] Inverting the first equation in \eqref{eq:maps:T_k,m_cross_0mu}, we obtain the following formula for the map $T_{k,m}^{-1}$:
\begin{equation}\label{eq:csble-t1:1}
\begin{array}{l}
X=(\alpha^{-1}+O(\delta))\bar X + (a \lambda^k\gamma^m)^{-1}u^- - x^+ + \psi_1(\bar X ,\bar{Y},Z),\\[5pt]
Y=\psi_2(\bar X,\bar{Y},Z),\\[5pt]
\bar{Z}= \psi_3(\bar X ,\bar{Y},Z),
\end{array}
\end{equation}
where 
\begin{equation}\label{eq:csble-t1:2}
\begin{array}{ll}
\psi_1=O(\delta^2)+o(1)_{k,m\to\infty},\qquad & \dfrac{\partial\psi_1}{\partial(X,\bar Y,Z)}= O(\delta)+o(1)_{k,m\to\infty},\\
 \|\psi_2\|_{C^1} =o(\hat\gamma^{-k}), \qquad & \|\psi_3\|_{C^1}=o(\hat\lambda^{m}).
\end{array}
\end{equation}
Observe that this formula has the same form as \eqref{eq:maps:T_k,m_cross_0mu}, with the replacement of
$\alpha$ by $\alpha^{-1}$ and the term $(ab\lambda^k\gamma^m x^+-bu^-)$ by the term
$((a \lambda^k\gamma^m)^{-1}u^- - x^+)$. So, since $|\alpha^{-1}|<1$, we obtain the result by repeating the same arguments we used for the case $|\alpha|<1$ - we only need to replace the cone field
$\mathcal{C}^{ss}$ by $\mathcal{C}^{uu}$, as we work with the inverse map $T_{k,m}^{-1}$.
\end{proof}

\subsection{A criterion for the existence of blenders}

Note that in proving Propositions  \ref{prop:cublendertype1} and \ref{prop:csblendertype1}, we did not use the full strength of the estimates in \eqref{eq:maps:nonlinearterms_0mu}. We essentially used the fact that the functions $\phi$ along with their first derivatives go to 0 as $\delta\to 0$ and $k,m\to\infty$. Additionally, we used the fact that $\phi_1=o(\delta)$ (hence $\tilde\phi_1=o(\delta)$ in \eqref{eq:xingsurfs:5}) in the last line in the proof of Lemma \ref{lem:crossingsurfaces}. With the above observation, we finish this section by the following summary which is used later for the saddle-focus and double-focus cases.

Take small $\delta>0$ and let $\Pi=[-\delta,\delta]\times [-\delta,\delta]^{d_Y}\times [-\delta,\delta]^{d_Z}$ be a cube in $\mathbb R^{d}$ for some $d_Y,d_Z\in \mathbb{N}$ satisfying $1+d_Y+d_Z=d$. Let $\{ T_n\}_{n=1}^{\infty}$ be a sequence of maps whose domains of definition $\sigma_n$ lie in $\Pi$ and  which satisfy 
$$\sigma_i\cap \sigma_j=\emptyset \quad\mbox{and}\quad T_i(\sigma_i)\cap T_j(\sigma_j)\cap \Pi =\emptyset,\qquad i\neq j.$$
Consider the map $T:\bigcup_n \sigma_n\to \mathbb{R}^d $ defined by $T(X,Y,Z)= T_n(X,Y,Z)$ if $(X,Y,Z)\in \sigma_n$.

\begin{prop}\label{prop:blender}
Suppose for every $n$ and every $(X,Y,Z)\in\Pi$ we have $(\bar X,\bar Y,\bar Z)=T_n (X,Y,Z)\in \Pi$ if and only if 
\begin{equation}\label{eq:blender:1}
\begin{aligned}
\bar X &= A_{n} X + B_{n} + \phi_1(X,\bar{Y},Z;n), \\
Y &=\phi_2(X,\bar{Y},Z;n), \\
\bar Z &= \phi_3(X,\bar{Y},Z;n),
\end{aligned}
\end{equation}
where the coefficients $A_n$ and $B_n$ and the functions $\phi$ defined on $\Pi$ may depend on $\delta$. Let  $\|\phi_{2,3}\|_{C^1}=o(1)_{\delta\to 0}+o(1)_{n\to\infty}$, whereas $\phi_{1}=o(\delta)+o(1)_{n\to\infty}$ and the derivative of $\phi_1$ is $o(1)_{\delta\to 0}+o(1)_{n\to\infty}$.

Assume that there exists a neighborhood $\Delta$ of zero independent of $\delta$ such that the set $\{B_n(\delta)\}$ is dense in $\Delta$ for each fixed small $\delta$. Assume also that the set $\{A_n(\delta)\}$ stays bounded away from $0$, $\pm 1$, and $\pm \infty$ for all sufficiently small $\delta$. Then, there exists $q\in(0,1)$ such that, for all sufficiently small $\delta$, the map $T$ in the cube
$\Pi'=[-q\delta,q\delta]\times [-\delta,\delta]^{d_Y}\times [-\delta,\delta]^{d_Z}\subset\Pi$
 has 
%
\begin{itemize}
\item a cs-blender of index $d_Y$ if $|A_n|<1$ for all $n$, or 
\item a cu-blender of index $(d_Y+1)$ if $|A_n|>1$ for all $n$.
\end{itemize}
The blender has an activating pair $(\Pi',\mathcal{C})$, where
$\mathcal{C}$ is a  field of cones  around $Z$-coordinates if the blender is center-stable, or around $Y$-coordinates if it is center-unstable. The cone field $\mathcal{C}$ can be chosen independent of $\delta$.
\end{prop} 

Note that in this paper, the maps $T_n$ are the first-return maps $T_{k,m}$. So, they are either iterations of a diffeomorphism or flow maps, depending on whether the system is discrete-time or continuous-time. The original system has a blender equal to the orbit of the blender of $T$. Like it is shown in Appendix for the blender of Proposition \ref{prop:cublendertype1}, the blenders of Proposition \ref{prop:blender}
are standard in the sense of Definition \ref{defi:blender}. 

\section{Local stabilization of heterodimensional cycles in the saddle case}\label{sec:sta_of_cycles}
In this section we prove Theorems \ref{thm:type1_family} - \ref{thm:type2_2pfamily} and Corollary \ref{cor:type2_2pfamily2}. We conduct the proofs only for the $|\alpha|<1$ case -- the $|\alpha|>1$ case is dealt with by using the symmetry argument (i.e., by considering the system obtained by the reversion of time) except for Theorem \ref{thm:type1_family}, where Proposition \ref{prop:csblendertype1} is additionally used. Depending on the situations, we embed $f$ into one- or two-parameter families that generically and properly unfold the heterodimensional cycle of $f$. In the remaining part of this paper, we denote the continuations of hyperbolic objects (e.g. $O_{1,2}$ and $\Lambda^{cs}$) after a small perturbation by the same letters and omit the term `continuation'.

\subsection{Unfolding type-I cycles. Proof of Theorem \ref{thm:type1_family}}\label{sec:thm:type1_family}
Recall that $\delta'$ is the size of the subcube $\Pi'\subset \Pi$ and it is related to the size $\delta$ of $\Pi$ via \eqref{eq:delta'}. 
Let us define two sequences of intervals
\begin{equation}\label{eq:type1_pert:includeO_1:0}
I^u_m=\left(\gamma^{-m} u^- -  \dfrac{1}{2}|b ^{-1}\gamma^{-m}|\delta',\gamma^{-m}u^-+ \dfrac{1}{2}|b ^{-1}\gamma^{-m}|\delta'\right).
\end{equation}
and
\begin{equation}\label{eq:type1_pert:includeO_2:1} 
I^s_k=\left(-a \lambda^k x^+ -\dfrac{1}{2}|a \lambda^k|\delta',-a \lambda^k x^+ +\dfrac{1}{2}|a \lambda^k|\delta'\right).
\end{equation}
The goal of this section is to prove the following general result for all saddle heterodimensional cycles, which not only implies Theorem \ref{thm:type1_family} but also will be used in the proof of Theorem \ref{thm:tied}
in Section \ref{sec:thm:tied}.

\begin{prop}\label{prop:inclusionlem}
Let the system $f$ have a non-degenerate saddle heterodimensional cycle (of type I or II).
\begin{itemize}
\item  Suppose $f$ has an index-$d_1$ cs-blender $\Lambda^{cs}$ with the activating pair $(\Pi',\mathcal{C}^{ss})$ defined in \eqref{eq:actdomian} and \eqref{eq:sscone}. Moreover, let the properties described in Remark \ref{rem:transintersection} hold for $\Lambda^{cs}$,  namely, $W^s_{loc}(\Lambda^{cs})$
intersects transversely any vertical disc that crosses $\Pi$ properly with respect to $\mathcal{C}^{uu}$, and $W^u_{loc}(\Lambda^{cs})$ 
intersects transversely $W^s(O_1)$. Then, 
\begin{itemize}
\item when $\hat\mu\in I^u_m$, the index-$d_1$ saddle $O_1$ is homoclinically related to $\Lambda^{cs}$, and
\item when $\hat\mu\in I^s_k$, there exist robust heterodimensional dynamics involving $\Lambda^{cs}$ and a non-trivial hyperbolic basic set containing $O_2$.
\end{itemize}
\item Suppose $f$ has an index-$d_2$ cu-blender $\Lambda^{cu}$ with the activating pair 
$(\Pi',\mathcal{C}^{uu})$, where $\mathcal{C}^{uu}$ is defined in \eqref{eq:uucone}. Moreover, let $W^u_{loc}(\Lambda^{cu})$
intersect transversely any horizontal disc that crosses $\Pi$ properly with respect to $\mathcal{C}^{ss}$, and $W^s_{loc}(\Lambda^{cu})$ 
intersect transversely $W^u(O_2)$. Then,
\begin{itemize}
\item when $\hat\mu\in I^s_k$, the saddle $O_2$ is homoclinically related to $\Lambda^{cu}$, and
\item when $\hat\mu\in I^u_m$, there exist robust 
heterodimensional dynamics involving $\Lambda^{cu}$ and a non-trivial hyperbolic basic set containing $O_1$.
\end{itemize}
\end{itemize}
\end{prop}

For type-I cycles, Propositions \ref{prop:cublendertype1} and   \ref{prop:csblendertype1} guarantee the existence of a cs-blender $\Lambda^{cs}$ when $|\alpha|<1$, and a cu-blender $\Lambda^{cu}$ when $|\alpha|>1$. Thus, Proposition \ref{prop:inclusionlem} gives

\begin{prop}\label{prop:intervalsfortype1}
Let the cycle $\Gamma$ have type I and let $\theta$ be irrational. 
\begin{itemize}
\item Let $|\alpha|<1$. Then the blender $\Lambda^{cs}$ has the same index as $O_1$ and, when 
$\hat\mu\in I^u_m$, the saddle $O_1$ is homoclinically related to $\Lambda^{cs}$. If $\hat\mu\in I^s_k$, then there exists a persistent non-transverse
heteroclinic connection between $W^u(\Lambda^{cs})$ and $W^s(O_2)$ and a transverse heteroclinic intersection of 
$W^u(O_2)$ and $W^s(\Lambda^{cs})$, i.e., there exist robust heterodimensional dynamics involving $\Lambda^{cs}$ and a non-trivial hyperbolic basic set containing $O_2$.
\item Let $|\alpha|>1$. Then the blender $\Lambda^{cu}$ has the same index as $O_2$ and, when 
$\hat\mu\in I^s_k$, the saddle $O_2$ is homoclinically related to $\Lambda^{cu}$. If $\hat\mu\in I^u_m$, then there exists a persistent non-transverse
heteroclinic connection between $W^u(O_1)$ and $W^s(\Lambda^{cu})$ and a transverse heteroclinic intersection of 
$W^u(\Lambda^{cu})$ and $W^s(O_1)$, i.e., there exist robust heterodimensional dynamics involving $\Lambda^{cu}$ and a non-trivial hyperbolic basic set containing $O_1$.
\end{itemize}
\end{prop}
This implies Theorem \ref{thm:type1_family} as follows.
\begin{proof}[Proof of Theorem \ref{thm:type1_family}]
Recall that $\hat\mu\sim \mu$ by \eqref{mmoh}, and all coefficients in formulas (\ref{eq:type1_pert:includeO_1:0}) and (\ref{eq:type1_pert:includeO_2:1}) for the intervals $I^u_m$ and $I^s_k$ depend continuously on the system. Thus, there exists $\kappa>0$ such that for any system $g$ in a small neighborhood of $f$, we have $\hat\mu \in I^u_m$ if
$$\mu(g)\gamma(g)^m \in[u^-_0-\dfrac{1}{2}b^{-1}_0\delta' + \kappa,u^-_0+\dfrac{1}{2}b^{-1}_0\delta' -\kappa]
$$
and $\hat\mu\in I^u_m$ if
$$\mu(g)\lambda(g)^k\in[-a_0x^+_0-\dfrac{1}{2}a_0\delta' + \kappa,-a_0x^+_0+\dfrac{1}{2}a_0\delta' -\kappa],$$
where $a_0$, $b_0$, $x^+_0$, and $u^-_0$ are the values of the corresponding coefficients for $f$. So, Proposition \ref{prop:intervalsfortype1} implies Theorem \ref{thm:type1_family} immediately.
\end{proof}

Let us now prove Proposition \ref{prop:inclusionlem}. The first step is to investigate the iterations of $W^u_{loc}(O_1)$ and $W^s_{loc}(O_2)$. Whether their iterates intersect the region `reserved' for the emergence of blenders is crucial for the stabilization of heterodimensional cycles. 

\begin{lem}\label{lem:type1_pert:includeO_1}
For all sufficiently large $k$ and $m$, and all sufficiently small $\delta$ (and hence small $\delta'$ by \eqref{eq:delta'}), the follow results hold:
\begin{itemize}
\item When $\hat\mu\in I^u_m$, the image $S^u_m:=F_{21}\circ F_2 \circ F_{12}(W^u_1)$ of
\begin{equation*}\label{eq:w^u}
W^u_1:=W^u_{loc}(O_1)\cap\{ \|y-y^{-}\| \leq \delta \}
\end{equation*}
is a `vertical' disc of the form $(X,Z)=s(Y)$ for some smooth function $s$. The disc $S^u_m$ crosses the cube $\Pi'$ (defined in \eqref{eq:actdomian}) properly with respect to the cone field $\mathcal{C}^{uu}$ defined in \eqref{eq:uucone}. In particular, $S^u_m$ intersects $W^s_{loc}(O_1):\{Y=0\}$ transversely, i.e., $O_1$ has a transverse homoclinic orbit when
$\hat\mu\in I^u_m$, and hence is contained in a non-trivial hyperbolic basic set.
\item When $\hat\mu\in I^s_k$, the preimage $S^{s}_k:=F^{-k}_{1}\circ F^{-1}_{12}({W}^s_2)$ of
\begin{equation*}\label{eq:w^s}
{W}^s_2:=W^s_{loc}(O_2)\cap\{ \|v-v^{+}\| \leq \delta\}
\end{equation*}
is a `horizontal' disc of the form $(X,Y)=s(Z)$ for some smooth function $s$. The disc $S^s_k$ crosses $\Pi'$ properly with respect to the cone field $\mathcal{C}^{ss}$ defined in \eqref{eq:sscone}. In particular, $S^s_k$ intersects $F_{21}(W^u_{loc}(O_2)):\{Z=0\}$ transversely (see (\ref{straightyz})), i.e., $O_2$ has a transverse homoclinic orbit when
$\hat\mu\in I^s_k$, and hence is contained in a non-trivial hyperbolic basic set.
\end{itemize}
\end{lem}

\begin{proof}
{\textit{(The first statement.)}} By formulas \eqref{eq:maps:F_12_cross} and \eqref{eq:maps:221}, for any $(0,y,0)\in W^u_1$, we have $F_{21}\circ F^m_2\circ F_{12}(0,y,0)=(\bar{x},\bar{y},\bar{z})$ if and only if
\begin{equation*}\label{eq:type1_pert:includeO_1:1}
\begin{aligned}
\bar{x}-x^+ &= b\gamma^m\hat\mu-b u^- + b_{13}\bar{y} +O((b\gamma^m\hat\mu-b u^-)^2+\bar y^2)+o(1)_{m\to\infty}, \\
\bar{z}-z^+ &= \hat b_{31}\gamma^m\hat\mu-b_{31}u^-   +  b_{33}\bar{y} +O((b\gamma^m\hat\mu-b u^-)^2+\bar y^2) +o(1)_{m\to\infty},
\end{aligned}
\end{equation*}
which after the coordinate transformation \eqref{eq:maps:coortransform1} recasts as
\begin{equation}\label{eq:type1_pert:includeO_1:2}
\begin{aligned}
\bar X &= b\gamma^m\hat\mu-b u^-  +O((b\gamma^m\hat\mu-b u^-)^2+\bar Y^2)+o(1)_{m\to\infty}, \\
\bar Z &=  O((b\gamma^m\hat\mu-b u^-)^2+\bar Y^2)+o(1)_{m\to\infty}.
\end{aligned}
\end{equation}
By \eqref{eq:type1_pert:includeO_1:0} we have
$$|b\gamma^m\hat\mu-b u^-| <\dfrac{\delta'}{2} ,$$
which for $\|\bar Y\|\leq \delta$ implies
$$|\bar X|<\dfrac{\delta'}{2}+O(\delta^2)+o(1)_{m\to\infty}< \delta'
\quad\mbox{and}\quad
\|\bar Z\|=O(\delta^2)+o(1)_{m\to\infty}<\delta.
$$
This means that $S^u_m$ crosses $\Pi'$. One also finds from \eqref{eq:type1_pert:includeO_1:2} that 
$$\dfrac{\partial (\bar X,\bar Z)}{\partial \bar Y} =O(\delta) +o(1)_{m\to\infty},$$
which can be made sufficiently small so that the tangent spaces of $S^u_m\cap\Pi$ lie in $\mathcal{C}^{uu}$. So, the crossing is also proper with respect to $\mathcal{C}^{uu}$. This completes the proof of the first statement.\\

{\textit{(The second statement.)}}
According to \eqref{eq:maps:F_12_cross}, the preimage $F^{-1}_{12}({W}_2^s)$ is given by
\begin{equation*}
\begin{aligned}
\tilde x&= -a^{-1}\hat\mu -{a_{13}}{a^{-1}_{11}}\tilde z + O(\tilde z^2+\hat\mu^2), \\
\tilde y&=y^- -{a_{31}}{a^{-1}}\hat\mu +(a_{33} - {a_{13}a_{31}}{a^{-1}}) \tilde z + O(\tilde z^2+\hat\mu^2),
\end{aligned}
\end{equation*}
By \eqref{eq:maps:F_1^k}, we find in coordinates \eqref{eq:maps:coortransform1} the disc
$F^{-k}_1\circ F^{-1}_{12}({W}^s_2)\cap \Pi$ as  
\begin{equation}\label{eq:type1_pert:includeO_2:2}
\begin{aligned}
X&= -a^{-1}\lambda^{-k}(\hat\mu + O(\hat\mu^2)) - x^+ + o(1)_{k\to\infty}, \\
Y&=o(\hat\gamma^{-k}),
\end{aligned}
\end{equation}
where $o(\cdot)$ and $O(\cdot)$ terms are functions of $(\hat\mu,Z)$, and the first derivatives satisfy
\begin{equation}\label{eq:type1_pert:includeO_2:3}
\dfrac{\partial (X,Y)}{\partial Z}=o(1)_{k\to\infty}.
\end{equation}

This disc crosses $\Pi'$ if $|X|<\delta'$ and $\|Y\|<\delta$. One readily finds from \eqref{eq:type1_pert:includeO_2:2} that this happens when $k$ is sufficiently large and
$$|-a^{-1}\lambda^{-k}\hat\mu - x^+|<\dfrac{\delta'}{2},$$
which gives the intervals \eqref{eq:type1_pert:includeO_2:1}. By \eqref{eq:type1_pert:includeO_2:3}, the crossing is proper with respect to $\mathcal{C}^{ss}$ for large $k$.
\end{proof}  

\begin{proof}[Proof of Proposition \ref{prop:inclusionlem}]
We only prove for the cs-blender  $\Lambda^{cs}$. In the case of the cu-blender, the arguments are completely analogous. By assumption, the local stable manifold of $\Lambda^{cs}$
intersects transversely (for $\mu=0$, hence for all small $\mu$; see Figure \ref{fig:forwardWu1}) any vertical disc that crosses $\Pi$ properly with respect to $\mathcal{C}^{uu}$, so it intersects
transversely $F_{21}\circ F^m_2 \circ F_{12}(W^u_{loc}(O_1)$ when $\hat\mu\in I^u_m$ (by the first item of Lemma \ref{lem:type1_pert:includeO_1}).
Our assumption on $\Lambda^{cs}$ also gives us the existence of a transverse intersection of $W^u_{loc}(\Lambda^{cs})$ and $W^s_{loc}(O_1)$, which proves 
that $O_1$ and $\Lambda^{cs}$ are homoclinically related for $\hat\mu\in I^u_m$.

Now, let $\hat\mu\in I^s_k$. With the formula for $F_{21}(W^u_{loc}(O_2))$ in \eqref{straightyz}, the existence of a non-trivial hyperbolic basic set containing $O_2$ is given by the second item of Lemma \ref{lem:type1_pert:includeO_1} (the existence of a transverse
homoclinic to $O_2$). Also by the same result, we immediately find that 
$O_2$ activates $\Lambda^{cs}$, meaning the existence of a persistent non-transverse heteroclinic connection between $W^s(O_2)$ and $W^u(\Lambda^{cs})$ (see Figure \ref{fig:backwardWs2}). It remains to note that 
the non-empty transverse intersection of $W^s(\Lambda^{cs})$ with $W^u(O_2)$ is given by the assumption on $\Lambda^{cs}$.

\begin{figure}[!h]
\begin{center}
\includegraphics[width=0.85\textwidth]{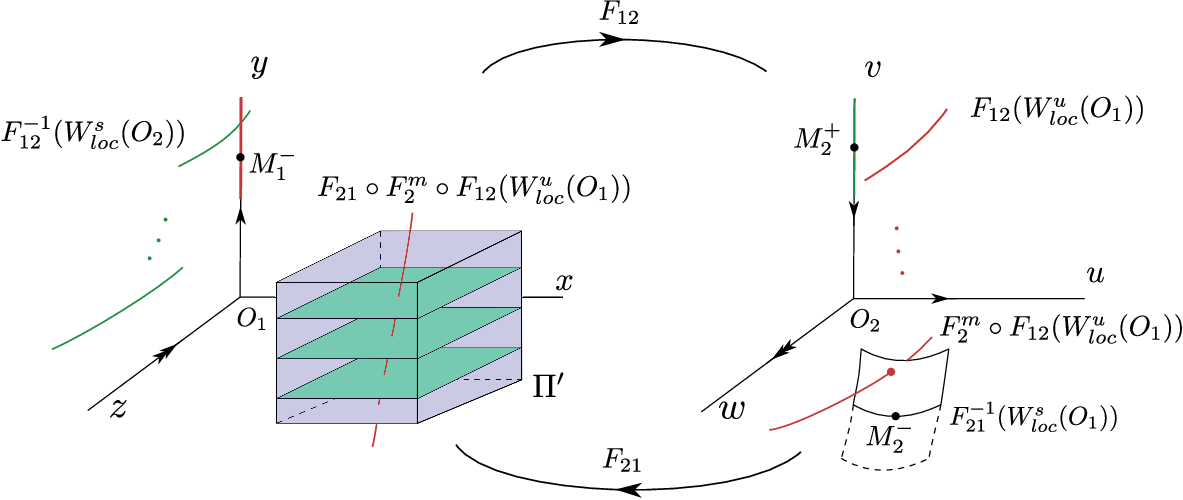}
\end{center}
\caption{For a cs-blender $\Lambda^{cs}$ satisfying the properties in Remark \ref{rem:transintersection}, there exist positive $\hat\mu$ values such that some forward iterate of $W^u_{loc}(O_1)$ intersects $W^s(\Lambda^{cs})$ while the backward iterates of $W^{s}_{loc}(O_2)$ leave the neighborhood of the cycle. Here the green planes in $\Pi' $ represent the local stable manifolds of points in the blender $\Lambda^{cs}$.}
\label{fig:forwardWu1}
\end{figure}

\begin{figure}[!h]
\begin{center}
\includegraphics[width=0.85\textwidth]{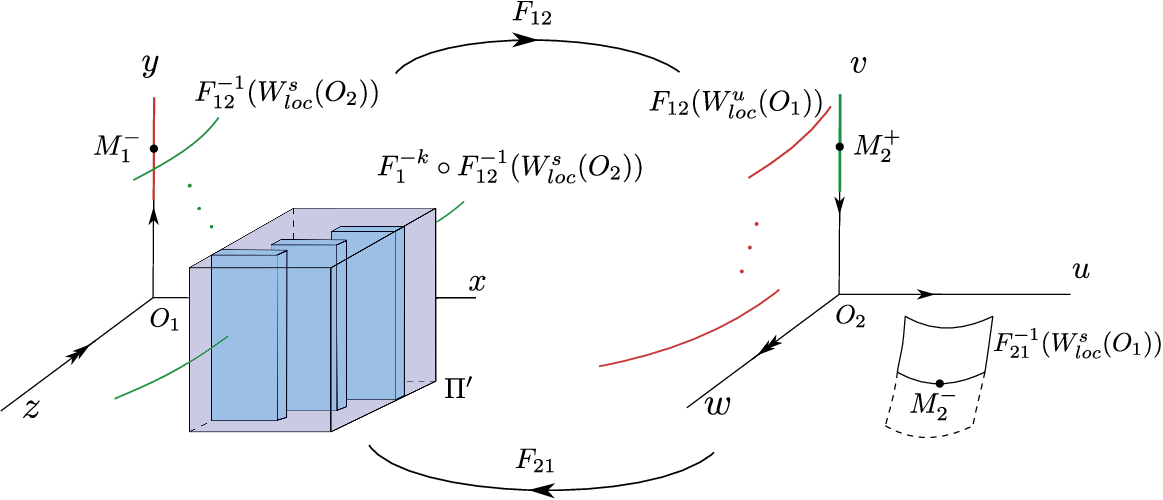}
\end{center}
\caption{
For a cs-blender $\Lambda^{cs}$ satisfying the properties in Remark \ref{rem:transintersection}, there exist negative $\hat\mu$ values such that some backward iterate of 
$W^s_{loc}(O_2)$ intersects $W^u(\Lambda^{cs})$ and the forward iterates of $W^{u}_{loc}(O_1)$ leave the neighborhood of the cycle.}
\label{fig:backwardWs2}
\end{figure}
\end{proof}

\subsection{Nonstabilizability of type-I cycles. Proof of Theorem \ref{thm:type1_family2}}\label{sec:thm:type1_family2}
By the assumption that the cycle is of type-I, we have $\lambda>0,\gamma>0$ and $ax^+u^->0$ at $\mu=0$, and hence  at all sufficiently small $\mu$. We further assume that $u^->0$ at $\mu=0$, which means that $ax^+>0$ too. It is clear from the proof that when $u^-<0$, all results remain true after reversing the sign of $\mu$. As mentioned before, we only consider the $|\alpha|<1$ case. 

We first prove the second part of the statement describing the situation at $\mu\leq 0$. It suffices to show that, at $\mu\leq 0$, the image 
$F^m\circ F_{12}(W^u_{loc}(O_1))$ does not enter the domain of $F_{21}$, i.e., it is outside a small neighborhood of $M^-_2=(u^-,0,w^-)$ for any sufficiently large $m$. Since $W^u_{loc}(O_1)$ is given by the equation
$(x=0,z=0)$, this is equivalent (see \eqref{eq:maps:F_12_cross} and \eqref{eq:maps:F_2^m}) to showing that the system of equations
$$
\begin{array}{ll}
 u=\hat\mu+a_{12}w +O(w^2), & u=\gamma^{-m} \tilde u +o(\gamma^m),\\
 v=a_{22}w +O(w^2),& \tilde v = o(\hat\lambda^m),\\
 y-y^-= a_{32}w+O(w^2),& w = o(\gamma^m),
\end{array}
$$
does not have a solution $(\tilde u,\tilde v,\tilde w)$ near $(u^-,0,w^-)$ for small $(y-y^-)$ and $w$, for any sufficiently large $m$. A straightforward computation reduces the above system to
\begin{equation*}
  \tilde u =  \gamma^{m}\hat\mu  + o(1)_{m\to\infty}.
\end{equation*}
Since we have here $\gamma>0$, the above equation together with the relation \eqref{mmoh} implies that $\tilde u\leq 0$ for $\mu\leq 0$ and sufficiently large $m$. Thus, by the assumption $u^->0$, the point $(\tilde u,\tilde v,\tilde w)$ never enters a small neighborhood of $(u^-,0,w^-)$ indeed.

Now we switch to the case $\mu>0$. First, note that like above, one can easily see (e.g. from equation \eqref{eq:type1_pert:includeO_2:2} for $F^{-k} \circ F^{-1}_{12} (W^s_{loc}(O_2))$) 
that no point of $W^s_{loc}(O_2)$ has an orbit that lies entirely in $U$ for $\mu>0$. In particular, no heteroclinic orbits of intersection of $W^u(O_1)$ and $W^s(O_2)$ can lie entirely in $U$.
Thus, given an orbit $\mathcal{O}$ that lies entirely in $U$, we have two possibilities: ether $\mathcal{O}$ is one of the orbits $L_{1,2}$ of the points $O_{1,2}$, or it intersects the $\delta$-neighborhood $\Pi$ of the heteroclinic point $M_1^+$ in $U_{01}$. Let $M_s$ be the consecutive points of intersection of $\mathcal{O}$ with $\Pi$; we have $M_{s+1}=T_{k_s,m_s} M_s$, where $T_{k,m}$ is the first-return map
given by (\ref{eq:maps:T_k,m_cross}) and $(k_s,m_s)$ are positive integers. The sequence $M_s$ can be infinite in both directions, infinite in one direction, or finite; let us consider the case where the sequence $\{M_s\}_{s\in\mathbb{Z}}$ is infinite in both directions first.

Let us show that the set $\hat\Lambda$ comprised by such orbits $\mathcal{O}$ is indeed a hyperbolic set with index-$d_1$ when $\mu>0$ (equivalently, when $\hat\mu>0$, see \eqref{mmoh}). 
It follows from \eqref{eq:maps:T_k,m_cross} and $M_{s+1}\in \Pi$ that
\begin{equation*}\label{eq:thm5:0}
\gamma^{m_s}\hat\mu+ a x^+ \lambda^{k_s}\gamma^{m_s}  = u^- + O(\delta).
\end{equation*}
Recall that we assumed $\gamma>0$, $\lambda>0$, $ax^+>0$, $u^->0$, and $\hat\mu>0$, so we obtain
\begin{equation}\label{eq:thm5:1}
0<a  x^+ \lambda^{k_s}\gamma^{m_s} < u^- + O(\delta).
\end{equation}
In particular, the numbers $(k_s,m_s)$ are always such that $\lambda^{k_s}\gamma^{m_s}$ is uniformly bounded. Therefore, by the same argument we used to deduce formula \eqref{eq:maps:T_k,m_cross_0mu} from \eqref{eq:maps:T_k,m_cross} at $\hat\mu=0$, we obtain for our case the following formula for 
$T_{k_s,m_s}:(X,Y,Z)\mapsto (\bar X, \bar Y, \bar Z)$ (where we use coordinates (\ref{eq:maps:coortransform1})):
\begin{equation}\label{eq:maps:T_k,m_cross_mu}
\begin{array}{l}
\bar X =b \gamma^{m_s}\hat\mu + a b \lambda^{k_s}\gamma^{m_s} x^+ - b u^- +  a b \lambda^{k_s}\gamma^{m_s} X  + \tilde\phi_1(X,  \bar{Y},Z),\\[5pt]
Y=\tilde\phi_2(X,  \bar{Y},Z),\qquad \bar Z = \tilde\phi_3(X,  \bar{Y},Z),
\end{array}
\end{equation}
where 
\begin{equation*}\label{eq:maps:nonlinearterms_mu}
\begin{array}{ll}
\tilde\phi_1=O(\delta^2)+o(1)_{k,m\to\infty},\qquad & \dfrac{\partial\tilde\phi_1}{\partial(X,\bar Y,Z)}= O(\delta)+o(1)_{k,m\to\infty},\\
 \|\tilde\phi_2\|_{C^1} =o(\hat\gamma^{-k}), \qquad & \|\tilde\phi_3\|_{C^1}=o(\hat\lambda^{m}).
\end{array}
\end{equation*}
Multiplying $|b/x^+|$ to each term in (\ref{eq:thm5:1}) yields
\begin{equation}\label{eq:thm5:2}
| a b \lambda^{k_s}\gamma^{m_s}|<|\alpha|+O(\delta)<1
\end{equation}
for all sufficiently small $\delta$. So, we have exactly the same formula \eqref{eq:conefields:00} for the derivative of $T_{k_s,m_s}$
(one should only replace $\alpha$ by 
$a b \lambda^{k_s}\gamma^{m_s}$). Also note that by shrinking the size of $U$, the values $k_s,m_s$ can be made arbitrarily large. Now, we have
the result of Lemma \ref{lem:conefields}, which along with \eqref{eq:thm5:2}, gives us the existence of invariant cone fields $\mathcal{C}^{cs}$ and $\mathcal{C}^{uu}$, which
immediately implies that $\hat\Lambda$ is a hyperbolic set with index-$d_1$ (all maps $T_{k_s,m_s}$ contract in $(X,Z)$ and expand in $Y$).

The results analogous to Remarks \ref{rem:strips} and \ref{rem:transintersection} also hold. In particular, every forward orbit $\{M_s\}$ that intersects $\Pi$ infinitely many times
must belong to $W^s_{loc}(\hat\Lambda\cap \Pi)$, while every backward orbit that intersects $\Pi$ infinitely many times must belong to $W^u_{loc}(\hat\Lambda\cap \Pi)$.
Also, $W^u_{loc}(\hat\Lambda\cap \Pi)$ intersects $W^s_{loc}(O_1)$ transversely, while $W^s_{loc}(\hat\Lambda\cap \Pi)$ intersects $F_{21}(W^u_{loc}(O_2))$ transversely.

We can now go through various cases where the sequence $M_s$ of consecutive intersections of $\mathcal{O}$ with $\Pi$ is not infinite in both directions.
\begin{enumerate}
\item $\{M_s\}_{s=-\infty}^{0}$ is infinite backwards. As we just explained, this may happen only when $M_0\in W^u_{loc}(\tilde{\Lambda}\cap\Pi)\cap W^s_{loc}(O_1)$ or 
$M_0\in W^u_{loc}(\tilde{\Lambda}\cap\Pi)\cap W^s(O_2)$. The latter case does not happen since no point of $W^s_{loc}(O_2)$ has an orbit entirely contained in $U$ for $\mu>0$, as we mentioned. In the former case, $M_0$ is a hyperbolic point of index $d_1$ (since $O$ and $\tilde\Lambda$ have the same index $d_1$ and the intersection of $W^u_{loc}(\tilde{\Lambda}\cap\Pi)$ and $W^s_{loc}(O_1)$ is transverse). We include all such orbits into the set $\Lambda$ of Theorem \ref{thm:type1_family2}, along with the set $\tilde\Lambda$ and the orbit $L_1$ (and orbits in $W^u(L_1)$ which we consider below).

\item $\{M_s\}_{s=0}^{\infty}$ is infinite forwards. This happens when $M_0\in W^s_{loc}(\tilde{\Lambda}\cap\Pi)\cap F_{21}\circ F_2^m\circ F_{12}(W^u_{loc}(O_1))$ for some large
$m$ or when $M_0\in W^s_{loc}(\tilde{\Lambda}\cap\Pi)\cap F_{21}(W^u_{loc}(O_2))$. The latter case, as the intersection is transverse, is in the complete agreement with the statement
of the Theorem.
In the former case, $M_0$ is a hyperbolic point of index $d_1$ because the intersection of $W^s_{loc}(\tilde\Lambda\cap\Pi)$ and $F_{21}\circ F_2^m\circ F_{12}(W^u_{loc}(O_1))$
is transverse. Indeed, the disc $F_{21}\circ F_2^m\circ F_{12}(W^u_{loc}(O_1))$ is given by equation (\ref{eq:type1_pert:includeO_1:2}) where we have
$$b\gamma^m\hat\mu-b u^- =O(\delta)$$
because $M_0\in\Pi$. From that, exactly like in Lemma \ref{lem:type1_pert:includeO_1}, one obtains that
the tangent spaces of $F_{21}\circ F_2 \circ F_{12}((W^u_{loc}(O_1))$ lie in $\mathcal{C}^{uu}$, hence the transversality with $W^s_{loc}(\tilde{\Lambda}\cap\Pi)$.

\item The sequence $\{M_s\}$ is finite. This happens when $M_0\in W^s(O_1)\cap F_{21}(W^u_{loc}(O_2))$ or $M_0\in W^s(O_1)\cap F_{21}\circ F_2 \circ F_{12}((W^u_{loc}(O_1))$.
Since the tangents to $W^s_{loc}(O_1):\{Y=0\}$ lie in in the backward-invariant cone field $\mathcal{C}^{cs}$, all its preimages by $T_{k_s,m_s}$ also have tangents in $\mathcal{C}^{cs}$,
hence the intersections of $W^s(O_1)$ with $F_{21}(W^u_{loc}(O_2))$ and $F_{21}\circ F_2 \circ F_{12}((W^u_{loc}(O_1))$ are transverse.
\end{enumerate}
In all cases we have the agreement with the statement of the theorem.
\qed

\subsection{Unfolding a tied pair of type-I and -II cycles. Proof of Theorem \ref{thm:tied}}\label{sec:thm:tied}
As can be seen from the proof of Theorem \ref{thm:type1_family2}, the main reason that prevents the blender from having homoclinic connections to both $O_1$ and $O_2$ is that, when one saddle connects to the blender, the iterates of the local stable or unstable manifold of the other leave the small neighborhood $U$ of $\Gamma$ (see Figures \ref{fig:forwardWu1} and \ref{fig:backwardWs2}). In what follows, we show that, if there exists a type-II cycle tied with $\Gamma$, then the leaving manifold can return to $U$ by following the robust heteroclinic orbit of the type-II cycle. 

Let $\tilde\Gamma=L_1\cup L_2\cup \Gamma^0 \cup\tilde\Gamma^1$ be a heterodimensional cycle sharing the fragile heteroclinic orbit $ \Gamma^0$ with $\Gamma$. By definition, any generic one-parameter unfolding of $\Gamma$ is also a generic one for $\tilde{\Gamma}$. These two cycles have the same transition map $F_{12}$, and the local maps $F_1$ and $F_2$ are the same. Denote by $\tilde F_{21}$ the transition map from a neighborhood of 
$\tilde M^-_2=(\tilde u^-,0,\tilde w^-)\in\tilde\Gamma^1$ to a neighborhood of $\tilde M^+_1=(\tilde x^+,0,\tilde z^+)\in\tilde\Gamma^1$.
The map $\tilde F_{21}$ has the same form as \eqref{eq:maps:F_21_cross}, just one needs to replace coefficients $b_{ij}$ by some $\tilde b_{ij}$, $b$ by $\tilde b\neq 0$ and also $x^+,z^+,u^-$ by $\tilde x^+,\tilde z^+,\tilde u^-$.

Recall that the coefficients in \eqref{eq:maps:F_21_cross}
depend continuously on $\mu$; they are smooth in $\mu$ when the smoothness class of $f_\mu$ is at least $C^3$. The eigenvalues
$\lambda$ and $\gamma$ remain smooth in $\mu$ in the $C^2$-case too. 

Let us mark the values of the coefficients $a$, $\tilde b$, $\lambda$, $\gamma$, $\tilde u^-$, $\tilde x^+$ at $\mu=0$ by the subscript ``$0$'' 

\begin{lem}\label{lem:type2_sechdc}
Let $\tilde \Gamma$ be a non-degenerate type-II cycle (i.e., $a_0 \tilde x^+_0 \tilde u^-_0 <0$ and $|{\tilde b_0 \tilde u^-_0}/{\tilde x^+_0}|\neq 1$). 
Consider a generic one-parameter unfolding $\{f_{\mu}\}$ of $\tilde \Gamma$. 
Assume $\theta_0$ is irrational and take a sequence $\{(k_j,m_j)\}$ of pairs of positive even integers satisfying $k_j,m_j\to \infty$ and 
\begin{equation}\label{eq:sechdc:0a}
a_0 \tilde b_0 \lambda_0^{k_j}\gamma_0^{m_j}\to -\tilde \alpha_0=-\dfrac{\tilde b_0 \tilde u^-_0}{\tilde x^+_0}
\end{equation}
as $j\to \infty$. There exists a sequence of values $\{\mu_j\}$ satisfying
\begin{equation}\label{eq:sechdc:0b}
\mu_j= - a_0\tilde x_0^+ \lambda_0^{k_j}+o(\lambda_0^{k_j}) = \tilde u^-_0 \gamma_0^{-m_j} + o(\gamma^{-m_j}),
\end{equation}
such that the system $f_{\mu_j}$ has a new orbit $\Gamma_j^{0,new}$ of heteroclinic intersection of $W^u(O_1)$ with $W^s(O_2)$, and hence a new heterodimensional cycle $\Gamma_j^{new}=L_1\cup L_2\cup \Gamma_j^{0,new} \cup\Gamma^1$. The heteroclinic connection $\Gamma_j^{0,new}$ splits when $\mu$ varies in an $o(\mu_j)$-interval around $\mu_j$. It splits with non-zero velocity if $f_\mu$ is of class $C^r$ with $r\geq 3$.
\end{lem}
\begin{rem}\label{remgen6}
The newly created heterodimensional cycle $\Gamma_j^{new}$ is non-degenerate. Indeed,
one only needs to check that $\Gamma_j^{0,new}$
satisfies condition C1 (the other non-degeneracy conditions hold because of the non-degeneracy of the cycle $\Gamma$). Condition
C1 holds for $\Gamma_j^{0,new}$ because $\Gamma_j^{new}$ is a partially-hyperbolic set with a one-dimensional center direction.
This is true because $\mbox{cl}(\Gamma_j^{0,new})$ lies in a small neighborhood of $\tilde\Gamma$ which is a compact partially-hyperbolic set 
with a one-dimensional center direction (by the non-degeneracy assumption), and the partial hyperbolicity of a compact invariant set is inherited 
by every closed invariant set in its neighborhood.
\end{rem}
\noindent{\em Proof of Lemma \ref{lem:type2_sechdc}.}
We create the new (fragile) heteroclinic orbit $\Gamma_j^{0,new}$ by finding values of $\mu$ at which $W^s_{loc}(O_2)$ intersects $\bar W^u_1:=F_{12}\circ F_1^k\circ \tilde F_{21}\circ F_2^m\circ F_{12}(W^u_{loc}(O_1))$.

By putting $\tilde x=0$, $\tilde z=0$ in \eqref{eq:maps:F_12_cross}, we find that the image $F_{12}(W^u_{loc}(O_1))$ near the point $M_2^+$
is given by
\begin{equation}\label{eq:sechdc:2}
u = \hat\mu+  O(w), \qquad  v =v^+ + O(w).
\end{equation}
Let us now find an equation for $\tilde F_{21}\circ F_2^m\circ F_{12}(W^u_{loc}(O_1))$. We use an analogue of formula \eqref{eq:maps:221} for the map $\tilde F_{21}\circ F_2^m$ from a small neighborhood of $M_2^+$ to
a small neighborhood of $\tilde M_1^+$, which now reads
\begin{equation*}
\begin{aligned}
u &= \gamma^{-m}(u^- +\tilde  b^{-1} ( x -\tilde x^+ - \tilde b_{13}  y + \hat h_{01}( x- \tilde x^+, y)))  + \hat h_1( x,v, y), \\[5pt]
 z - \tilde z^+ &= \tilde b_{31}\tilde b^{-1}(x-\tilde x^+-\tilde b_{13}y)+\tilde b_{33}y + \hat h_{02}( x-\tilde x^+, y)   + \hat h_2( x,v, y)
\end{aligned}
\end{equation*}
with the modified functions $\hat h$ satisfying the same estimates \eqref{htha}. By substituting equation \eqref{eq:sechdc:2} into this formula, we obtain that
$\tilde F_{21}\circ F_2^m\circ F_{12}(W^u_{loc}(O_1))$ is given by
\begin{equation*}\label{eq:maps:221++}
\begin{aligned}
x &= \tilde x^+ + \tilde b (\gamma^m \hat\mu - \tilde u^-)  +  O(\|y\| + (\gamma^m \hat\mu - \tilde u^-)^2) + o(1)_{m\to\infty} , \\
z  &= \tilde z^+  + O(\|y\| + |\gamma^m \hat\mu - \tilde u^-|) + o(1)_{m\to\infty},
\end{aligned}
\end{equation*}
where $m$ must be such that $\gamma^m \hat\mu$ is sufficiently close to $\tilde u^-$.

Finally, substituting  the above two equations into \eqref{eq:maps:112} yields the equation for
$\bar W^u_1:=F_{12}\circ F_1^k\circ \tilde F_{21}\circ F_2^m\circ F_{12}(W^u_{loc}(O_1))$:
\begin{equation*}\label{eq:sechdc:6}
\begin{aligned}
u&= \hat\mu + a\tilde b\lambda^k\gamma^m\hat\mu + a\lambda^k\tilde x^+ - a\tilde b\lambda^k\tilde u^- +O(w) + \lambda^k O((\gamma^m\hat\mu-\tilde u^-)^2)+o(\lambda^k),\\
v &=v^+ + O(w) + O(\lambda^k).\\
\end{aligned}
\end{equation*}
The heteroclinic orbit $\Gamma^{0,new}$ corresponds to an intersection of $\bar W^u_1$ with $W^{s}_{loc}(O_2)$, which corresponds to letting 
$(u,w)=0$ in the above two equations, i.e.,
\begin{equation}\label{eq:sechdc:7}
0 = \hat\mu + a\tilde b\lambda^k\gamma^m\hat\mu + a\lambda^k\tilde x^+ - a\tilde b\lambda^k\tilde u^-+ 
\lambda^k O((\gamma^m\hat\mu-\tilde u^-)^2)+o(\lambda^k).
\end{equation}

We will look for solutions $\mu=O(\lambda^k)$ for large enough and even $(k,m)$ from the sequence $\{k_j,m_j\}$ satisfying \eqref{eq:sechdc:0a},
i.e.,
\begin{equation}\label{mktlog}
m - k \theta_0\to \frac{1}{\ln|\gamma_0|} \ln \left|\dfrac{\tilde u^-_0}{a_0\tilde x^+_0}\right|,
\end{equation}
where the $0$ subscript stands for the value at $\mu=0$. Since $\lambda$ and $\gamma$ depend smoothly on $\mu$, we have, in particular, {\em Lipshitz dependence of} $\theta=-\ln|\lambda|/\ln|\gamma|$
{\em on} $\mu$, implying
\begin{equation}\label{liptheta}
\theta=\theta_0 + O(\mu)= \theta_0+ O(\lambda^k).
\end{equation}
Since $a_0\tilde x^+_0\tilde u^-_0<0$ by assumption, and $k,m$ are even, the above equation along with \eqref{mktlog} implies
\begin{equation}\label{kmbo}
\lambda^k\gamma^m = |\gamma|^{m-k\theta}
=|\gamma_0|^{m-k\theta_0 +O(k\lambda^k)}
\left|\dfrac{\gamma}{\gamma_0}\right|^{-k\theta_0 +O(k\lambda^k)}
=- \;\dfrac{\tilde u^-_0}{a_0\tilde x^+_0}
+o(1)_{k,m\to\infty}.
\end{equation}
Note that we used here Lipshitz dependence of $\theta$ on $\mu$, but we only used continuity of $\gamma$'s dependence on $\mu$.

Now, we rewrite \eqref{eq:sechdc:7} as
$$0 = (\gamma^m\hat\mu -\tilde u^-)(1 + a\tilde b\lambda^k\gamma^m) + \tilde u^- + a\lambda^k\gamma^m\tilde x^+ + 
\lambda^k\gamma^m O((\gamma^m\hat\mu-\tilde u^-)^2)+o(\lambda^k\gamma^m),$$
or, by (\ref{kmbo}) and continuous dependence of coefficients on $\mu$,
$$0= (\gamma^m\hat\mu -\tilde u^-)(1 - \;\dfrac{\tilde b\tilde u^-_0}{\tilde x^+_0} +o(1)_{k,m\to\infty}) +
O((\gamma^m\hat\mu-\tilde u^-)^2) + o(1)_{k,m\to\infty}.$$
Since ${\tilde b\tilde u^-_0}/{\tilde x^+_0}=\tilde\alpha_0 \neq 1$ (as given by the non-degeneracy assumption C4.1 for the cycle $\tilde\Gamma$),
we find that the new heteroclinic connection $\Gamma_j^{0,new}$ exists when
\begin{equation}\label{muhatjg}
\hat\mu_j  = \gamma^{-m_j} (\tilde u^-_0 +o(1)_{j\to\infty}) = - \lambda^{k_j} (a_0\tilde x^+_0 + o(1)_{j\to\infty})
\end{equation}
for large enough even $(k_j,m_j)$ defined by \eqref{eq:sechdc:0a}. By \eqref{mmoh}, this can be rewritten as
\begin{equation}\label{mujg}
\mu_j  = \gamma^{-m_j} (\tilde u^-_0 +o(1)_{j\to\infty}) = - \lambda^{k_j} (a_0\tilde x^+_0 + o(1)_{j\to\infty}).
\end{equation}
Note that $\gamma$ and $\lambda$ in this formula (as well as the $o(1)$ terms) depend on $\mu$, so this is an implicit relation on $\mu_j$.

If $f_\mu$ is of class $C^3$ with respect to variables and $\mu$, then
the right-hand side of (\ref{mujg}) is smooth with respect to $\mu$,
so we obtain a unique solution $\mu_j$ for each sufficiently large $j$ by the Implcit Function Theorem. It also follows that the corresponding heteroclinic connection $\Gamma_j^{0,new}$ splits with non-zero velocity as $\mu$ varies across $\mu_j$.

If $f_\mu$ is only $C^2$, we can only guarantee a continuous dependence of the right-hand side of (\ref{mujg}) on
$\mu$. Still, it is obvious that the continuity gives the existence (though not necessarily uniqueness) of the solution $\mu_j$ for all sufficiently large $j$. Thus, we have the existence of the sought fragile heteroclinic connection $\Gamma_j^{0,new}$ for some $\mu_j\in int(\Delta_j)$ where the interval
$\Delta_j$ corresponds to
\begin{equation}\label{deltaj6}
|\mu_j  - \gamma^{-m_j} \tilde u^-_0| < \kappa_j \gamma^{-m_j}, \qquad
|\mu_j  + a_0\tilde x^+_0 \lambda^{k_j}| < \kappa_j \lambda^{k_j},
\end{equation}
for some $\kappa_j\to 0$. Note that we can choose $\kappa_j$ going to 0 sufficiently slowly, so that $\Gamma_j^{0,new}$ necessarily splits when $\mu$ gets out of $\Delta_j$.

We finish the proof of the lemma by recalling that $\lambda$ and $\gamma$ depend smoothly on $\mu$, hence
$$\lambda=\lambda_0+O(\mu), \qquad \gamma=\gamma_0+O(\mu),$$
which, upon substitution into \eqref{mujg}, gives \eqref{eq:sechdc:0b}.  
\qed

Since the cycle $\tilde\Gamma$ is tied with $\Gamma$, we have $x^+=\tilde x^+$ or $u^-=\tilde u^-$. It follows immediately from formula
(\ref{muhatjg}) that if $x^+=\tilde x^+$, then $\hat\mu_j$ lies in the interval $I^s_k$ defined in \eqref{eq:type1_pert:includeO_2:1}, and
if $u^-=\tilde u^-$, then $\hat\mu_j$ lies in the interval $I^u_m$ defined in  \eqref{eq:type1_pert:includeO_1:0}. In either case, Proposition 
\ref{prop:intervalsfortype1} implies that the blender which exists by Theorem \ref{thm:blender_in_type1} near the type-I cycle $\Gamma$ 
is robustly connected\footnote{Here by robust connection between a blender and a saddle we mean that they are homoclinically related if they have the same index, or they are involved in robust heterodimensional dynamics if they have different indices.} to one of the saddles $O_{1,2}$ when the system has the newly created fragile heteroclinic orbit $\Gamma_j^{0,new}$ of Lemma \ref{lem:type2_sechdc}.

Theorem \ref{thm:blender_in_type1} gives us the blender at $\mu=0$, which persists for all small $\mu$, and satisfies conditions of Proposition \ref{prop:inclusionlem}. The robust connection to the saddle
$O_1$ or $O_2$ also persists for all $\mu$ corresponding to the interval $I^u_m$ or $I^k_s$, i.e., it persists for all $\mu$ from the intervals $\Delta_j$
defined by (\ref{deltaj6}). Let us show that there are
sub-intervals inside $\Delta_j$ where the blender has a robust connection to the other saddle too. 

\begin{figure}[!h]
\begin{center}
\includegraphics[width=0.9\textwidth]{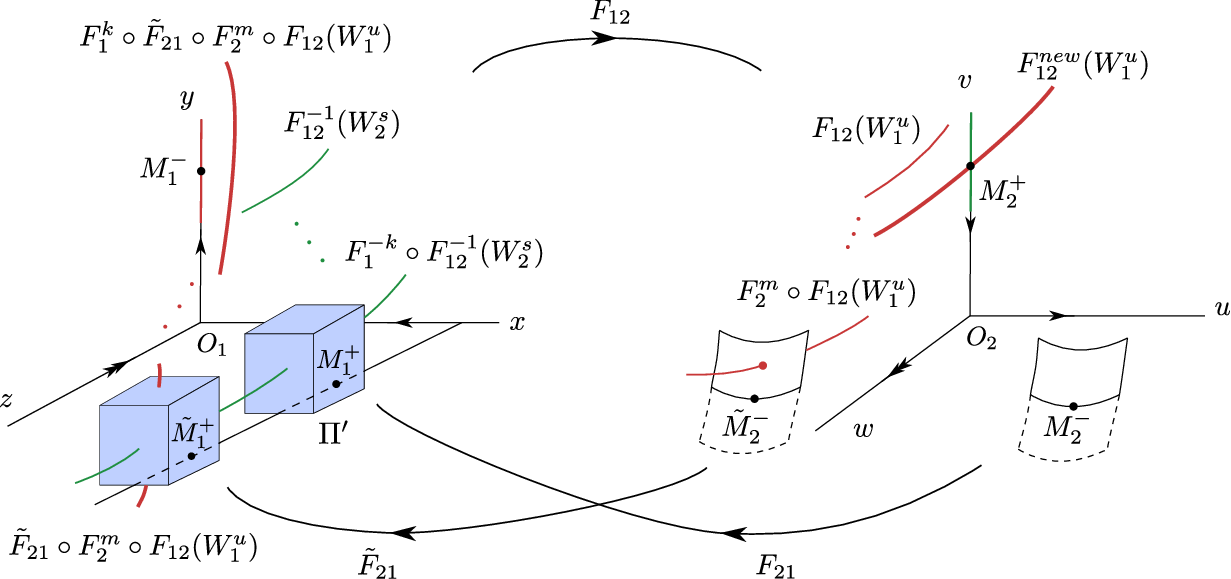}
\end{center}
\caption{A pair of tied heterodimensional cycles of type I and II with $M^+_1$ and $\tilde{M}^+_1$ lying in the same strong-stable leaf and with $x^+>0$, $u^->0$. There exist negative $\mu$ values such that some backward iterate of $W^s_{loc}(O_2)$ intersects $W^u(\Lambda^{cs})$, and the forward iterate $F_{12}\circ F^k_{1}\circ \tilde F_{21}\circ F^m_2\circ F_{12}(W^u_{loc}(O_1))=: F^{new}_{12}(W^u_{loc}(O_1))$ intersects $W^s_{loc}(O_2)$ producing a new heterodimensional cycle. After a further change in $\mu$, so that the newly obtained fragile heteroclinics splits to the right (i.e., $F^{new}_{12}(W^u_1)$ lies in $\{u>0\}$), one can find $m'$ satisfying $F^{m'}_2\circ F^{new}_{12}(W^u_1)\cap F^{-1}_{21}(W^s_{loc}(O_1))\neq \emptyset$ such that the iterate $F_{21}\circ F^{m'}_2\circ F^{new}_{12}(W^u_1)$ crosses $\Pi' $ vertically, and, therefore, it intersects $W^s(\Lambda^{cs})$. As a result, the two saddles $O_1$ and $O_2$ get connected to the blender at the same time.}
\label{fig:tied}
\end{figure}

Indeed, when $\mu$ varies from one end of $\Delta_j$ to
the other, the heteroclinic orbit $\Gamma_j^{0,new}$ splits. Therefore, we can apply Proposition \ref{prop:inclusionlem} to the heterodimensional cycle
$\Gamma_j^{new}=L_1\cup L_2\cup \Gamma^{0,new} \cup \Gamma^1$ at $\mu=\mu_j$ (that is, system $f$ in Proposition \ref{prop:inclusionlem} is our
system $f_{\mu_j}$; the coefficient $\hat\mu$ in Proposition \ref{prop:inclusionlem} measures the distance between $\bar W^u_1$ with $W^{s}_{loc}(O_2)$ when $\Gamma_j^{0,new}$ splits; note that Proposition \ref{prop:inclusionlem} is applicable because $\Gamma_j^{new}$ is non-degenerate, see Remark \ref{remgen6}). As a result, we find that
there are intervals of $\mu$ values inside $\Delta_j$ which correspond to a robust connection between
the blender and $O_1$, as well as intervals which correspond to a robust connection between
the blender and $O_2$ (the connections can be transverse, or persistent heterodimensional). In any case, we have found
intervals of $\mu$ values for which the blender is connected to both saddles $O_1$ and $O_2$; see Figure \ref{fig:tied}.
This completes the proof of Theorem \ref{thm:tied}. \qed

\begin{rem}\label{rem:smoothness6}
In the proof of this theorem we only used the smoothness with respect to $\mu$ in two places: when we claimed the Lipshitz dependence
of $\theta$ on $\mu$ (see \eqref{liptheta}), and when we inferred \eqref{eq:sechdc:0b} from \eqref{mujg}. However, we did not use formula
\eqref{eq:sechdc:0b} for $\mu_j$ when we derived Theorem \ref{thm:tied} from Lemma \ref{lem:type2_sechdc}, and used only 
relation \eqref{mujg} for $\mu_j$. The existence of such values $\mu_j$ is established for any continuous family in which $\Gamma^0$ splits,
i.e., for which the splitting functional is not constantly zero and changes sign. Therefore, Theorem \ref{thm:tied} remains intact, for example, for any
one-parameter family going through $f$ such that
\begin{itemize}[nolistsep]
\item the splitting functional $\mu$ changes sign, and
\item $\theta$ stays equal to the same irrational constant.
\end{itemize}
Even more, the result still holds if we replace the one-parameter unfolding by any connected set of systems which contains $f$ and satisfies the
above conditions -- the intervals $I_j$ in the formulation of Theorem \ref{thm:tied} should be replaced by open subsets $I_j$ converging to $f$.
\end{rem}

\subsection{Creating a tied pair of type-I and -II cycles. Proofs of Theorem \ref{thm:type2_2pfamily} and Corollary \ref{cor:type2_2pfamily2}}\label{sec:thm:type2_2pfamily}
\begin{figure}[!h]
\begin{center}
\includegraphics[width=0.9\textwidth]{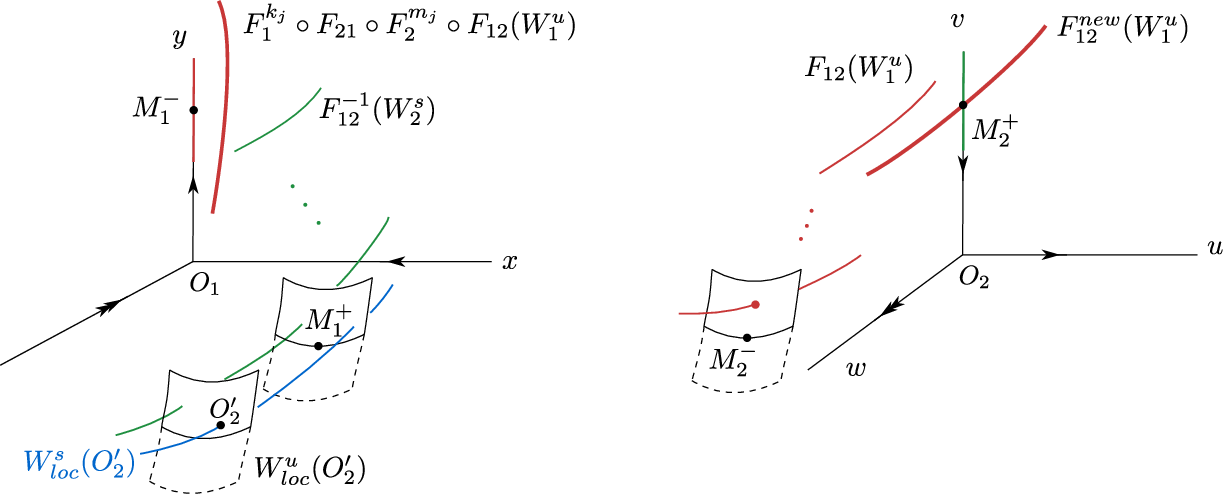}
\end{center}
\caption{At the moment of a new cycle given by Lemma \ref{lem:type2_sechdc} at $\mu=\mu_j$, we find an index-($d_1+1$) point $O_2'$ near  $M^+_1$, homoclinically related to $O_2$. After a small perturbation, we can transfer the fragile heteroclinic intersection of $W^s(O_2)$ with $W^u(O_1)$ to  $W^s(O'_2)$. We will show that there are two robust heteroclinics between $O_1$ and $O'_2$, giving rise to two heterodimensional cycles of type I and II.}
\label{fig:type2a}
\end{figure}

Theorem \ref{thm:tied} immediately shows that cycles of type III can be stabilized. In this section we prove the same for type-II
cycles in the following way. We let $|\alpha|<1$ (the case $|\alpha|>1$ is reduced to this one by the time-reversal). Then, we show that a generic one-parameter unfolding of a type-II cycle gives rise to a pair of tied cycles of type I and type II associated with $O_1$ and to a new saddle periodic point $O_2^\prime$ of index $d_1+1$, see Figure \ref{fig:type2a}. Next, we compute the value of the modulus $\theta$ for the new cycles: 
$\theta^\prime:=-\ln|\lambda|/\ln|\gamma^\prime|$, where $\gamma^\prime$ is the center-unstable multiplier of $O_2^\prime$, and show that by taking $\theta$ of the original cycle as a second parameter and varying it, we can make $\theta^\prime$ irrational, so that Theorem \ref{thm:tied} becomes applicable.

\noindent{\it Proof of Theorem \ref{thm:type2_2pfamily}.}
Assume $|\alpha|<1$. Since $\Gamma$ is a cycle of type II, Lemma \ref{lem:type2_sechdc} is applicable to it (just replace $\tilde \Gamma$
by $\Gamma$ in the formulation of the lemma, and remove the tildes from the coefficients in \eqref{eq:sechdc:0a} and
\eqref{eq:sechdc:0b}). This gives that
for a sequence $\{(k_j,m_j)\}$ of pairs of positive even integers satisfying $k_j,m_j\to \infty$ and 
\begin{equation}\label{eq:type2entangled:1}
a_0  b_0 \lambda_0^{k_j} \gamma_0^{m_j} \to -\dfrac{ b_0  u_0^-}{ x_0^+}
\end{equation}
as $j\to \infty$, there exists a sequence of values $\{\mu'_j\}$
\begin{equation}\label{eq:type2entangled:2}
\mu'_j= - a_0 x_0^+ \lambda_0^{k_j}+o(\lambda_0^{k_j})
\end{equation}
such that the system $f_{\mu'_j}$ has a new heterodimensional cycle associated with $O_1$ and $O_2$. By \eqref{mmoh}, we also
have
\begin{equation}\label{eq:tied2:0}
\hat\mu_j':=\hat\mu(\mu_j')=- a_0 x_0^+ \lambda_0^{k_j}+o(\lambda_0^{k_j})
\end{equation}
for the coefficient $\hat\mu$ in formula \eqref{eq:maps:T_k,m_cross} for the first-return map near $\Gamma$.

Denote $m'_j=m_j+m^*$ for some fixed even integer $m^*$ (to be determined below). In what follows, we find for each sufficiently large $j$ a fixed point of the first-return map $T_{k_j,m'_j}$ for the values of $\mu$ which are $o(\mu_j')$-close to $\mu_j'$. Note first that due to the smooth dependence of
$\lambda$ and $\gamma$ on the parameter $\mu$, we obtain from \eqref{eq:type2entangled:2} and \eqref{eq:type2entangled:1} that
\begin{equation}\label{eq:dependence}
\begin{array}{l}
\lambda(\mu_j)^{k_j}=(\lambda_0+O(\mu_j))^{k_j}=\lambda^{k_j}_0+O(k_j\lambda^{k_j}_0\mu_j)=
\lambda^{k_j}_0 (1+ O(k_j\lambda_0^{k_j}),\\
\gamma(\mu_j)^{m_j}=(\gamma_0+O(\mu_j))^{m_j}=\gamma^{m_j}_0+O(m_j\gamma^{m_j}_0\mu_j)=\gamma_0^{m_j}
(1+O(m_j\gamma_0^{-m_j})).
\end{array}
\end{equation}
Thus,
\begin{equation}\label{eq:tied2:1}
\lambda^{k_j}\gamma^{m'_j}\to -\dfrac{u^-_0\gamma^{m^*}}{a_0 x^+_0},
\end{equation}
so, in particular, $\lambda^{k_j}\gamma^{m'_j}$ stays uniformly bounded as $j\to\infty$. 

With this, one notices that the term $\gamma^{m} \hat\phi_1$ in formula \eqref{eq:maps:T_k,m_cross}
for the first-return map $T_{k,m}$ is $o(1)_{j\to\infty}$ when $k=k_j$ and $m=m_j'$. Hence, with the estimates in \eqref{oriph0} and \eqref{eq:maps:nonlinearterms_original}, we can write the map
$T_{k_j,m_j'}$ as
\begin{equation}\label{eqcross}
\begin{aligned}
\bar X &=b \gamma^{m_j'}\hat\mu+ a b \lambda^{k_j}\gamma^{m_j'} x^+ - b u^- + a b \lambda^{k_j}\gamma^{m_j'} X  + 
O(\bar X^2 + \bar Y^2) + o(1)_{j\to\infty},\\[5pt]
Y&=o(\hat\gamma^{-k_j}), \qquad \bar Z = o(\hat\lambda^{m_j'}).
\end{aligned}
\end{equation}
It follows that the fixed point $(X_j,Y_j,Z_j)$ of $T_{k_j,m'_j}$ satisfies
\begin{equation}\label{eq:type2entangled:4}
\begin{aligned}
X_j&=\dfrac{-(a\lambda^{k_j})^{-1}\hat\mu_j' + (a \lambda^{k_j}\gamma^{m'_j})^{-1}u^- - x^+ + O(X_j^2)+o(1)_{k_j,m_j'\to\infty}}
{1-(ab\lambda^{k_j}\gamma^{m'_j})^{-1}},\\[5pt]
Y_j&=O(\hat\gamma^{-k_j}),\qquad Z_j = O(\hat\lambda^{m_j'});
\end{aligned}
\end{equation}
note that $1-(ab\lambda^{k_j}\gamma^{m'_j})^{-1}\neq 0$, as follows from \eqref{eq:tied2:1} if $m^*$ is sufficiently large.

We need to verify that $(X_j,Y_j,Z_j)\in \Pi$. It suffices to show that $|X_j|<\delta$ for all large $k_j$ and $m'_j$.
Substituting \eqref{eq:tied2:0}, \eqref{eq:dependence}, and \eqref{eq:tied2:1} into the first equation of \eqref{eq:type2entangled:4} yields
\begin{equation}\label{eq:type2entangled:5}
X_j = x_\infty+ o(1)_{j\to\infty},
\end{equation}
where
\begin{equation}\label{xinf}
x_\infty= - x^+_0\gamma_0^{-m^*}+O(\gamma_0^{-2m^*}).
\end{equation}
Thus, by taking $m^*$ sufficiently large and $\delta$ sufficiently small, we obtain that $|X_j|<\delta$ for all sufficiently large $j$, as required.

We have shown that at $\mu=\mu_j'$ the system has a non-degenerate heterodimensional cycle associated with $O_1$ and $O_2$ and
the first-return map $T_{k_j,m_j'}$ has a fixed point in $\Pi$ given by \eqref{eq:type2entangled:5} for any
sufficiently large $m^*$. We denote this point as $O_2^\prime$. By \eqref{eqcross}, the map $T_{k_j,m_j'}$ strongly contracts in $Z$ and strongly expands in $Y$. Also, since
$$|a b \lambda^{k_j}\gamma^{m'_j}|=\left|\dfrac{b_0u^-_0\gamma_0^{m^*}}{x^+_0}\right|+o(1)_{j\to\infty}\gg 1,$$
the map is expanding in $X$ if $m^*$ is taken large enough. Thus, the point $O_2^\prime$ is a saddle of index $d_1+1$.

Arguing as in the proof of Lemma \ref{lem:conefields}, one finds, on the set of points whose images under $T_{k_j,m_j'}$ belong to $\Pi$, a forward-invariant unstable cone field around the $(X,Y)$-space and a forward-invariant strong-unstable cone-field around the $Y$-space, and also,
on the set of points whose preimages belong to $\Pi$, a backward-invariant stable cone field around the $Z$-space.
This implies that $W^u_{loc}(O_2^\prime)$ and $W^s_{loc}(O_2^\prime)$ are given by $Z=Z_j+w^u_j(X,Y)$ and 
$(X,Y)=(X_j,Y_j)+w^s_j(Z)$, respectively, for a smooth function $w^u_j$ defined for $(X,Y)\in[-\delta,\delta]\times[-\delta,\delta]^{d_1}$ and a smooth function $w^s_j$ defined for $Z\in[-\delta,\delta]^{d-d_1-1}$. Also, $W^u_{loc}(O_2^\prime)$ contains the strong-unstable manifold $W^{uu}_{loc}(O_2^\prime)$ of the form $\{X=X_j+w^{uu}_j(Y), Z=Z_j+w^u_j(X,Y)\}$ where the smooth function $w^{uu}_j$ is defined for $Y\in[-\delta,\delta]^{d_1}$. Note that it immediately follows from \eqref{eqcross} that in the limit $j\to \infty$ we have 
\begin{equation}\label{manlim}
W^u_{loc}(O_2^\prime)\to \{Z=0\}, \qquad W^s_{loc}(O_2^\prime)\to \{X=x_\infty, Y=0\}, \qquad 
W^{uu}_{loc}(O_2^\prime)\to \{X=x_\infty, Z=0\}
\end{equation}
in the $C^1$-topology.

\begin{figure}[!h]
\begin{center}
\includegraphics[width=0.6\textwidth]{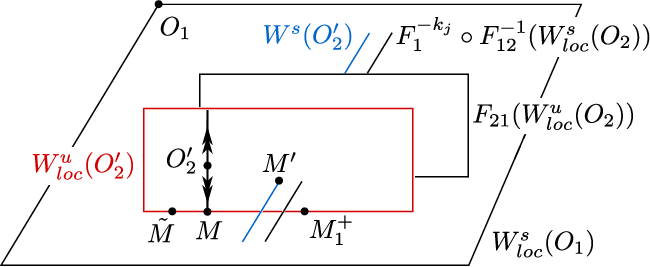}
\end{center}
\caption{Creation of the heteroclinic  point $\tilde M$ and the homoclinic point $M'$, which are on different sides of $W^{uu}_{loc}(O'_2)$.}
\label{fig:type2b}
\end{figure}

Now, since $F_{21}(W^u_{loc}(O_2))$ is given by $\{Z=0\}$ (see \eqref{straightyz}), we obtain that it intersects $W^s_{loc}(O_2^\prime)$ transversely. Since the value $\hat\mu'_j$ lies in the interval $I^s_{k_j}$ in Lemma \ref{lem:type1_pert:includeO_1}, it follows that 
$F_1^{-k_j}\circ F^{-1}_{12}(W^s_{loc}(O_2))$ is a horizontal disc that crosses $\Pi$ properly with respect to $\mathcal{C}^{ss}$ in Lemma \ref{lem:conefields}, so it intersects $W^u_{loc}(O_2^\prime)$
transversely. Thus, $O_2^\prime$ is homoclinically related to $O_2$. This, in particular, implies that $O_2^\prime$ has transverse homoclinics. See Figure \ref{fig:type2b} for an illustration.

Let us show that in $W^u_{loc}(O_2^\prime)$
there exists a homoclinic point $M'$ of transverse intersection of $W^u_{loc}(O_2^\prime)$ with $W^s(O_2^\prime)$, such that
$M'\not\in W^{uu}_{loc}(O_2^\prime)$. By (\ref{manlim}), we just need to show that the $X$-coordinate of the homoclinic point
$M'$ can be kept bounded away from $x_\infty$. Since $W^s_{loc}(O_2^\prime)$ intersects $F_{21}(W^u_{loc}(O_2))$ transversely,
the preimages $F_2^{-m}\circ F_{21}^{-1}(W^s_{loc}(O_2^\prime))$ accumulate (in $C^1$) to $W^s_{loc}(O_2)$ as $m\to\infty$. 
Hence, we can choose the homoclinic point $M'$ as close as we want to a point of transverse intersection of $W^u_{loc}(O_2^\prime)$ 
with the preimage $F_1^{-k_j}\circ F_{12}^{-1}(W^s_{loc}(O_2))$. This preimage is given by the equation 
$$x = -\frac{1}{a} \hat\mu \lambda^{-k_j} + o(1)_{k_j\to\infty}, \qquad  y=o(1)_{k_j\to\infty},$$
as follows from substitution of $(u=0,w=0)$ (the equation of $W^s_{loc}(O_2)$) into equation \eqref{eq:maps:112} for the map
$F_{12}\circ F_1^k$. In the coordinates \eqref{eq:maps:coortransform1}, we get
$$X = -\frac{1}{a} \hat\mu \lambda^{-k_j} - x^+ + o(1)_{k_j\to\infty}.$$
This implies, by \eqref{eq:tied2:0} and \eqref{eq:dependence}, that for any $\mu$ which is $o(\mu_j')$-close to $\mu=\mu_j'$
the $X$-coordinate of the intersection of $F_1^{-k_j}\circ F_{12}^{-1}(W^s_{loc}(O_2))$ with $W^u_{loc}(O_2^\prime)$ is $o(1)_{j\to\infty}$. Hence, the $X$-coordinate of the homoclinic
point $M'$ can be made as close to zero as we want, i.e., it is bounded away, as claimed, from $x_\infty$ (which is non-zero by (\ref{xinf})), see Figure \ref{fig:type2b}.

Now, let us construct a heteroclinic intersection of $W^u(O_1)$ with $W^s(O_2^\prime)$. The newly created heterodimensional cycle at $\mu=\mu'_j$ includes an orbit of a non-transverse intersection between $W^u(O_1)$ and $W^s(O_2)$. This orbit splits as $\mu$ varies 
in an interval of size $o(\mu_j')$ (by Lemma \ref{lem:type2_sechdc}). Since $O_2$ and $O_2^\prime$ are homoclinically related, the invariant manifold $W^s(O_2^\prime)$ accumulates on $W^s(O_2)$, which means that when $\mu$ varies, an orbit of a non-transverse intersection of $W^u(O_1)$ and $W^s(O_2^\prime)$ emerges. Let it happen at $\mu=\mu_j$; it is $o(\mu_j')$-close to $\mu'_j$, and by \eqref{eq:type2entangled:2} we have
\begin{equation}\label{eq:mu_j}
\mu_j= - a_0 x_0^+ \lambda_0^{k_j}+o(\lambda_0^{k_j}).
\end{equation}
Note that the orbit $\Gamma^{0}_j$  of the fragile heteroclinic intersection of $W^u(O_1)$ and $W^s(O_2^\prime)$ at $\mu=\mu_j$ satisfies condition C1 by the partial-hyperbolicity argument as in Remark \ref{remgen6}.

To finish the proof, we need to find a pair of orbits $\Gamma^1_{j,I}$ and $\Gamma^1_{j,II}$ of a transverse intersection of $W^u(O_2^\prime)$ and $W^s(O_1)$ such that the corresponding heterodimensional cycles
$\Gamma_{j,I}=L_1\cup L_2^\prime\cup \Gamma^{0}_j\cup \Gamma^{1}_{j,I}$ and
$\Gamma_{j,II}=L_1\cup L_2^\prime\cup \Gamma^{0}_j\cup \Gamma^{1}_{j,II}$
(where $L_2^\prime$ is the orbit of the periodic point $O_2^\prime$) satisfy the remaining non-degenerate conditions C2 - C4.1, are tied to each other, and have different types (I and II).

We, first, notice that $W^s_{loc}(O_1):\{Y=0\}$ transversely intersects the local strong-unstable manifold $W^{uu}_{loc}(O_2^\prime)$
(indeed, by (\ref{manlim}), $W^{uu}_{loc}(O_2^\prime)$ is $C^1$-close to $\{X=x_\infty, Z=0, Y\in[-\delta,\delta]^{d_1}\}$). Let $M$ be the point of intersection. The local unstable manifold $W^{u}_{loc}(O_2^\prime)$ is divided by $W^{uu}_{loc}(O_2^\prime)$ into two connected components. It follows from the transversality of the intersection of $W^{uu}_{loc}(O_2^\prime)$ with $W^s_{loc}(O_1)$ that the intersection
of $W^s_{loc}(O_1)$ with $W^{u}_{loc}(O_2^\prime)$ near $M$ is a curve that goes from one component to another while crossing
$W^{uu}_{loc}(O_2^\prime)$ at $M$. Choose a point $\tilde M$ on this curve such that $\tilde M$ is close to $M$ and $\tilde M$ lies in a different component of $W^{u}_{loc}(O_2^\prime)\setminus W^{uu}_{loc}(O_2^\prime)$ from the homoclinic point 
$M'\in W^u_{loc}(O_2^\prime)\cap W^s(O_2^\prime)$.

The orbit of $\tilde M$ is an orbit of transverse intersection of $W^{u}(L_2^\prime)$ and $W^s(L_1)$. The non-degeneracy condition C2
is satisfied by this orbit due to partial hyperbolicity of every orbit lying in a small neighborhood of the original non-degenerate cycle $\Gamma$
(see Remark \ref{remgen6}). Condition C3 is satisfied because we take $\tilde M\neq M$, i.e., $\tilde M$ is not in $W^{uu}_{loc}(O_2^\prime)$. However, we take $\tilde M$ sufficiently close to $W^{uu}_{loc}(O_2^\prime)$, so the $u$-coordinates of $\tilde M$ is sufficiently
close to zero, making the non-degeneracy condition C4.1 fulfilled too. The non-degenerate heterodimensional cycle $\Gamma^*$ comprised by the orbit of $\tilde M$, the periodic orbits  $L_1$ and $L_2^\prime$ and the fragile heteroclinic $\Gamma_0^j$ at $\mu=\mu_j$ is the sought cycle $\Gamma_{j,I}$ or
$\Gamma_{j,II}$ (depending on whether it is type I or type II). 

By Remark \ref{rem:tiedcycles}, the fact that $\tilde M$ and the homoclinic point $M^\prime$ belong to different components of $W^{u}_{loc}(O_2^\prime)\setminus W^{uu}_{loc}(O_2^\prime)$ implies that the cycle $\Gamma^{*}$ is tied with a cycle $\Gamma^{**}$ of a different type. So, we have the sought pair of heterodimensional cycles $\Gamma_{j,I}$ and
$\Gamma_{j,II}$, and hence the result of the theorem, once we show the non-degeneracy of $\Gamma^{**}$ below.

The orbit of transverse intersection of $W^s(O_1)$ and $W^u(O'_2)$ given by Remark \ref{rem:tiedcycles} (i.e., the robust heteroclinic of $\Gamma^{**}$) intersects 
$W^u_{loc}(O'_2)$ at some point $M''$ close to the homoclinic point $M^\prime$. By the partial hyperbolicity argument, we have 
the non-degeneracy condition C2 for this orbit; condition C3 holds because $M^\prime\not\in W^{uu}_{loc}(O_2^\prime)$. Let us establish the last non-degeneracy condition C4.1. The orbit of the point $M''$ goes close to the orbit of $M'$ and gets back to a small neighborhood 
of $O_2^\prime$ (since $M'$ is homoclinic). After that the orbit spends a long time near $L_2^\prime$, which corresponds to a large number of iterations of the first-return map $T_{k_j,m_j+m^*}$ near $O_2$, before getting to $W^s_{loc}(O_1)$. The iterations near $O_2$ create a very large expansion in the central direction, i.e., the expansion factor $b$ in \eqref{eq:intro:b} gets very large
for this orbit, making $|\alpha|\gg 1$ in \eqref{eq:intro:4}. Thus, the non-degeneracy condition C4.1 is fulfilled for this orbit, provided
$M''$ is chosen close enough to $M'$.
\qed

We can now finish the proof of Theorem \ref{thm:persis_hdd_2p} for the case of real central multipliers.

\noindent{\it Proof of Corollary \ref{cor:type2_2pfamily2}.} In a proper unfolding $f_\varepsilon$ of a heterodimensional cycle of type II,
we fix the values of all parameters except for $\mu$ and $\theta$ and show that an arbitrary small change of $\mu$ and $\theta$
can lead the system into the region of robust heterodimensional dynamics involving $O_1$ and $O_2$. So, we may from the very beginning
consider a two-parameter family $f_{\mu,\theta}$ where $\mu$ varies in a small neighborhood of zero and $\theta$ varies in a neighborhood
of some irrational $\theta=c_0$. By Theorem \ref{thm:type2_2pfamily}, on the line $\theta=c_0$ we have (see \eqref{eq:mu_j}) a sequence of values
$$\mu_j= - a_0 x_0^+  \lambda_0^{k_j}(1+o(1)_{j\to\infty})$$
for which the system has a tied pair of cycles of different types involving $O_1$ and a periodic point $O_2^\prime$, homoclinically related to $O_2$. 
The coefficients
$a_0$, $x_0^+$, $\lambda_0$ correspond to $\mu=0$ and depend continuously on the parameter $\theta$; the $o(1)$ term
also depends continuously on $\theta$ and the rate at which it tends to zero is uniform with respect to $\theta$.

Moreover, as proved in Lemma \ref{lem:type2_sechdc}, when $\mu$ varies (on the line $\theta=c_0$) within a centered at $\mu=\mu_j$ interval $\Delta_j$ of size $o(\mu_j)$, the fragile heteroclinic connection $\Gamma^0$ shared by the tied cycles splits. Namely, there exist $\mu_j^-= \lambda_0^{k_j} (- a_0 x_0^+-\kappa)$ and $ \mu_j^+= \lambda_0^{k_j} (- a_0 x_0^+ +\kappa)$ with some small $\kappa>0$ such that $\Gamma^0$ is split at $(\mu=\mu_j^\pm,\theta=c_0)$ in opposite directions. The same holds true at $\mu=\mu_j^\pm$ for any value of $\theta$ close to $c_0$. It follows that for any connected set $\cal L$ sufficiently close to the line $\theta=c_0$ in the $(\mu,\theta)$-plane such that $\cal L$ has a point at the line
$\mu= \mu_j^-$ and another point at the line 
$\mu= \mu_j^+$, there is a point in $\cal L$ corresponding to the existence of $\Gamma^0$, i.e., to the tied pair of cycles involving
$O_1$ and $O_2^\prime$. Moreover, $\Gamma^0$ splits when we move within $\cal L$ from $\mu= \mu_j^-$ to $\mu= \mu_j^+$. Therefore, we immediately get the result (the existence of robust heterodimensional
dynamics involving $O_1$ and $O_2^\prime$ - hence, $O_2$) by applying the version of Theorem \ref{thm:tied} given in Remark
\ref{rem:smoothness6}, if we can choose the connected set $\cal L$ such that the modulus 
$\theta^\prime= -\ln |\lambda|/\ln|\gamma^\prime|$ stays equal to a certain irrational constant everywhere on $\cal L$ 
(recall that we denote by $\gamma^\prime$ the central multiplier of the point $O_2^\prime$).

For that, we just need to show that arbitrarily close to $c_0$ there exist constants $c_-$ and $c_+$ such that, if $j$ is large enough,
\begin{equation}\label{thetachange}
\sup_{\mu\in[\mu_j^-, \mu_j^+]} \theta^\prime(\mu,\theta=c^-) \; < \; 
\inf_{\mu\in[\mu_j^-, \mu_j^+]} \theta^\prime(\mu,\theta=c^+).
\end{equation}
Let us compute the central multiplier $\gamma^\prime$ of $O_2^\prime$. It is, up to $o(1)_{j\to\infty}$ terms, the derivative 
${d \bar X}/{dx}$ of the map (\ref{eqcross}) at $O_2^\prime$. This gives us
$$\gamma^\prime = a b \lambda^{k_j}\gamma^{m_j'}(1+ O(\delta)) + o(1)_{j\to\infty},$$
or, according to \eqref{eq:tied2:1},
$$\gamma^\prime = - \; \frac{b_0 u^-_0}{x^+_0}\gamma^{m^*}(1 + O(\delta)) + o(1)_{j\to\infty}.$$
This gives
\begin{equation}\label{eq:type2_2pfamily:1}
\theta^\prime=\frac{1}{m^*} \theta + O((m^*)^{-2}),
\end{equation}
so, by taking $m^*$ sufficiently large, we obtain \eqref{thetachange}.
\qed 

\subsection{Heterodimensional cycles with rational $\theta$. Proof of Theorem \ref{thm:rationaltheta}}\label{sec:rational}
Recall that by a heterodimensional cycle we mean the set $\Gamma=L_1\cup L_2\cup \Gamma^0\cup\Gamma^1$, where $L_{1,2}$ are the orbits of $O_{1,2}$, and $\Gamma^0$ and $\Gamma^1$ are the fragile and robust heteroclinic orbits, respectively. In this section we prove that the complexity of dynamics in a sufficiently small neighborhood of a heterodimensional cycle $\Gamma$ in the saddle case depends essentially on whether $\theta$ is rational or not. That is, in contrast to the appearance of blenders for irrational $\theta$, we show the persistence of simple dynamics near $\Gamma$ for generic one-parameter unfoldings when the $\theta$ value (for the unperturbed system) is rational. This is done by investigating the structure of the set consisting of points whose entire orbits lie in a small neighborhood $U$ of $\Gamma$.

Note that such orbits, except for those in the stable and unstable manifolds of $L_1$ and $L_2$, must intersect the neighborhood $\Pi$ 
of the heteroclinic point $M_1^+$ in $U_{01}$ (the small neighborhood of $O_1$ contained in $U$) infinitely many times both forwards and backwards in time.
If such orbit exists, then for any two consecutive intersection points of the orbit with $\Pi$, the second one is the image of the first one under the map $T_{k,m}$ given by \eqref{eq:maps:T_k,m_cross} for some pair of integers $(k,m)$. Moreover, $k$ and $m$ have to be large enough if $U$ is small. Recall that for a saddle heterodimensional cycle, we use $\lambda$ and $\gamma$ to denote the center-stable multiplier of $O_1$ and, respectively, center-unstable multiplier of $O_2$.

\begin{lem} \label{lemkey9} Let $f_\mu$ be a generic one-parameter unfolding of the system $f_0$ with a non-degenerate heterodimensional cycle involving two saddles,
such that $\theta(f_0)=-\frac{\ln|\lambda|}{\ln|\gamma|}=\frac{p}{q}$ is rational. Assume condition \eqref{rare2} is satisfied. Let $\{(k_s,m_s)\}$ be a sequence (finite or infinite) of different pairs of integers
such that, for some small $\mu$ and some points $M_i=(X_i,Y_i,Z_i)\in \Pi$,  we have
$T_{k_i,m_i}(M_i)=(\bar X_i,\bar Y_i,\bar Z_i)=:\bar M_i\in\Pi$
for all $i$.
Then, for all sufficiently small $\delta$, either all $m_i$ are equal to each other, $m_i=m$, and 
\begin{equation}\label{meq12}
\lambda^{k_i} = O(\delta) \gamma^{-m}
\end{equation}
for all $i$, or all $k_i$ are equal to each other, $k_i=k$, and
\begin{equation}\label{keq12}
\gamma^{-m_i} = O(\delta) \lambda^{k}
\end{equation}
for all $i$.
\end{lem}
\begin{rem}\label{reminfkm}
Note that in this lemma we also allow $k_i$ or $m_i$ to be infinite
(with the convention $\gamma^{-\infty}=0$, $\lambda^\infty=0$).
Here, $k_i=\infty$ corresponds to 
$M_i\in W^s_{loc}(O_1)$ and $\bar M_i$ being a point of intersection of $W^u(L_1)$ with $\Pi$; 
more specifically, $\bar M_i \in F_{21}\circ F_2^{m_i} \circ F_{12} (W^u_{loc}(O_1))$.
The case $m_i=\infty$ corresponds to $M_i$ being a point of intersection of $W^s(L_2)$ with 
$\Pi$, i.e., $M_i \in F_1^{-k_i}\circ F_{12}^{-1} (W^s_{loc}(O_2))$ and $\bar M_i \in F_{21}(W^u_{loc}(O_2))$.
\end{rem}

\begin{proof}[Proof of Lemma \ref{lemkey9}] Consider, first, the case where we have only 2 pairs in the sequence $(k_i,m_i)$. 
By \eqref{eq:maps:T_k,m_cross}, conditions $T_{k_1,m_1}(M_1)=\bar M_1$ and $T_{k_2,m_2}(M_2)=\bar M_2$ give
\begin{equation}\label{eqfir}\!\!\!\!\!\!\!
\begin{array}{l}
\hat\mu+a  \lambda^{k_1}X_1+ a  \lambda^{k_1}x^+ - \gamma^{-m_1}  u^- 
+\gamma^{-m_1}O((\bar X_1)^2+(\bar Y_1)^2)+o(\lambda^{k_1})+o(\gamma^{-m_1})=\gamma^{-m_1}b^{-1}\bar X_1,\\
 \hat\mu+a  \lambda^{k_2}X_2+ a  \lambda^{k_2}x^+ - \gamma^{-m_2} u^- +\gamma^{-m_2}O((\bar X_2)^2+(\bar Y_2)^2)+o(\lambda^{k_2})+o(\gamma^{-m_2})=\gamma^{-m_2}b^{-1}\bar X_2.
\end{array}
\end{equation}
Assumptions that $M_{1,2}\in \Pi$ and $\bar M_{1,2}\in \Pi$ imply $|X_{1,2}| <\delta$, $|\bar X_{1,2}|<\delta$, $\|\bar Y_{1,2}\|<\delta$. Thus,
it follows from \eqref{eqfir} that the system
\begin{equation}\label{eq:theta:1}
\begin{array}{l}
\hat\mu+ a  \lambda^{k_1}(x^+ + K_1\delta ) - \gamma^{-m_1} u^- = \gamma^{-m_1}b^{-1}C_1\delta,\\
 \hat\mu+ a \lambda^{k_2}(x^+ + K_2\delta ) - \gamma^{-m_2} u^- = \gamma^{-m_2}b^{-1}C_2\delta,
\end{array}
\end{equation}
must have a solution $(\hat\mu,k_1,m_1,K_1,C_1,k_2,m_2,K_2,C_2)$ with $|K_{1,2}|<1$, $|C_{1,2}|<1$. Subtracting the second equation of \eqref{eq:theta:1} from the first one, yields
\begin{equation}\label{eq:theta:2}
\lambda^{k_1}-\lambda^{k_2}\dfrac{x^+ + K_2\delta}{x^++K_1\delta}
= \gamma^{-m_1}\dfrac{C_1\delta+bu^-}{ab(x^+ + K_1\delta)}-\gamma^{-m_2}\dfrac{C_2\delta+bu^-}{ab(x^+ + K_1\delta)}.
\end{equation}

Recall that $\lambda$, $\gamma$, $a$, $b$, $x^+$, and $u^-$ depend on $\mu$. Let us indicate their values at $\mu=0$ by the subscript ``$0$''.
Since the multipliers $\lambda$ and $\gamma$ depend smoothly on $\mu$, we have
\begin{equation}\label{eq:dependence0}
\lambda^{k}=(\lambda_0+O(\mu))^{k}=\lambda^{k}_0+O(k\lambda^{k}_0\mu),\quad
\gamma^{-m}=(\gamma_0+O(\mu))^{-m}=\gamma^{-m}_0+O(m\gamma^{-m}_0\mu).
\end{equation}
Using \eqref{mmoh} and \eqref{eq:dependence0}, one can estimate $\mu$ from \eqref{eq:theta:1} as
\begin{equation*}\label{eq:theta:mu}
\mu = O(|\lambda_0|^{k_1}+|\gamma_0|^{-m_1}) \quad\mbox{and}\quad \mu = O(|\lambda_0|^{k_2}+|\gamma_0|^{-m_2}).
\end{equation*}
Substituting this into \eqref{eq:dependence0}, and also using the continuous dependence of all the coefficients on $\mu$,
we rewrite \eqref{eq:theta:2} as
$$
\lambda^{k_1}_0-\lambda^{k_2}_0\dfrac{x^+_0 + K_2\delta}{x^+_0+K_1\delta} (1+o(1))
=
\gamma^{-m_1}_0\dfrac{C_1\delta+b_0 u^-_0}{a_0b_0(x^+_0 + K_1\delta)}(1+o(1))
-\gamma_0^{-m_2}\dfrac{C_2\delta+b_0 u^-_0}{a_0b_0(x^+_0 + K_1\delta)}(1+o(1)),
$$
where $o(1)$ denotes terms that tend to zero as $k_{1,2}\to \infty$, $m_{1,2}\to\infty$. Note that these estimates make sense since when $\delta$ is sufficiently small, the $k$ and $m$ satisfying the assumption must be sufficiently large.

It is obvious that this equation is not solvable for sufficiently small $\delta$ and sufficiently large $k_1\neq k_2$ and $m_1\neq m_2$, unless
the quantity $u^-_0/(a_0x^+_0)$ is a limit point of the set
\begin{equation}\label{setb}
\begin{array}{l}
\left\{\dfrac{\lambda_0^{k_1}-\lambda_0^{k_2}}{\gamma_0^{-m_1}-\gamma_0^{-m_2}}:\; k_{1,2}\in\mathbb{N}, m_{1,2}\in\mathbb{N},
\;k_1\neq k_2, \; m_1\neq m_2\right\}.
\end{array}
\end{equation}
Since 
\begin{equation}\label{eq:lambdagamma}
|\lambda_0|=|\gamma_0|^{-\theta_0},
\end{equation}
and $\theta_0$ is a rational number $p/q$ (with $p$ and $q$ coprime), we obtain, assuming $k_1>k_2$ and $m_1>m_2$, that
\begin{equation*}
\left|\dfrac{\lambda_0^{k_1}-\lambda_0^{k_2}}{\gamma_0^{-m_1}-\gamma_0^{-m_2}}\right|
=
\dfrac{|\lambda_0|^{k_2}}{|\gamma_0|^{-m_2}}\cdot
\left|\dfrac{\gamma_0^{-(k_1-k_2)\theta}-1}{\gamma_0^{-(m_1-m_2)}-1}\right|
=|\gamma_0|^{\frac{-k_2p+m_2q}{q}}
\dfrac{1-\lambda_0^{k_1-k_2}}{1-\gamma_0^{-(m_1-m_2)}}.
\end{equation*}
This implies that the absolute values of the limit points of the set \eqref{setb} form the set
$$\displaystyle
\mbox{cl}\,\left\{|\gamma_0|^{\frac{s}{q}} \dfrac{1-\lambda_0^{l}}{1-\gamma_0^{-n}} \right\}_{s\in\mathbb{Z},l\in\mathbb{N},n\in\mathbb{N}}
$$
and, therefore, by \eqref{rare2}, we find that $u^-_0/(a_0x^+_0)$ is not a limit point of the set \eqref{setb}.  Thus, the system \eqref{eq:theta:1} can have a solution (for sufficiently small $\delta$ and large $k_{1,2}$, $m_{1,2}$)
only if $k_1=k_2$ or $m_1=m_2$. 

If $m_1=m_2=m$, then we have
$$\lambda^{k_1}-\lambda^{k_2}\dfrac{x^+ + K_2\delta}{x^++K_1\delta} = O(\delta)\gamma^{-m}$$
from \eqref{eq:theta:2}, which implies \eqref{meq12} when $k_1\neq k_2$. 
If $k_1=k_2=k$, then we have
$$\gamma^{-m_1}\dfrac{C_1\delta+bu^-}{ab(x^+ + K_1\delta)}-\gamma^{-m_2}\dfrac{C_2\delta+bu^-}{ab(x^+ + K_1\delta)}=
O(\delta)\lambda^k$$
from \eqref{eq:theta:2}, which implies \eqref{keq12} when $m_1\neq m_2$.

This proves the lemma in the particular case of two pairs $(k_i,m_i)$.
If we have more pairs, take the first two in the sequence. As we just proved, we either have $m_1=m_2=m$, or $k_1=k_2=k$. The arguments are the same in both cases, so we consider only the case
$m_1=m_2=m$; note that we then have
\begin{equation}\label{again}
\lambda^{k_1}\ll \gamma^{-m_1}
\end{equation}
by (\ref{meq12}).
We want to prove that $m_i=m$ for all $i$ in this case, so assume, to the contrary, that $m_i\neq m_1=m$ for some $i$. Then, as we just proved, $k_i=k_1$, and
$$\gamma^{-m_1}\ll \lambda^{k_1}$$
by (\ref{keq12}). This is a contradiction with (\ref{again}), which proves
the claim.
\end{proof}

\begin{proof}[Proof of Theorem \ref{thm:rationaltheta}]
The case $k_1=m_1=\infty$ in Lemma \ref{lemkey9} corresponds (see
Remark \ref{reminfkm}) to the orbit of the fragile heteroclinic intersection of $F_{12}(W^u_{loc}(O_1))$ with $W^s_{loc}(O_2)$ (so we may think of $M_1\in W^s_{loc}(O_1)$ and $\bar M_1 \in F_{21}(W^u_{loc}(O_2))$). Such intersection exists at $\mu=0$. Note that
in this case we also have $k_2=m_2=\infty$. Indeed, if we assume $m_2\neq m_1$, then $k_2=k_1=\infty$ and, by \eqref{keq12}, 
we get $\gamma^{-m_2}=0$, i.e., $m_2=\infty$, and if we assume $k_2\neq k_1$, then $m_2=m_1=\infty$ and, by \eqref{meq12}, 
$\lambda^{k_2}=0$, i.e., $k_2=\infty$, too. This means that no other orbits in $U$ can intersect $\Pi$ in this case except for the orbits corresponding to the transverse intersection of $F_{21}(W^u_{loc}(O_2))$ with $W^s_{loc}(O_1)$. This gives us the result of Theorem \ref{thm:rationaltheta} at $\mu=0$.

At $\mu\neq 0$, let there exist an orbit $L$ in $U$, different from $L_1$, $L_2$ and from heteroclinic orbits corresponding to the intersection of $F_{21}(W^u_{loc}(O_2))$ with $W^s_{loc}(O_1)$ (we call $L$ a non-exceptional orbit). Then, $L$ intersects $\Pi$ in a sequence of points $M_s$ such that $M_{s+1}=T_{k_s,m_s} (M_s)$. If this sequence is infinite, then all $(k_s,m_s)$ consist of finite positive integers. If the sequence is finite from the left, then either the most left point 
$M_{s_l} \in F_{21}(W^u_{loc}(O_2))\cap \Pi$, or
$M_{s_l} \in F_{21}\circ F_2^m \circ F_{12}(W^u_{loc}(O_1))\cap \Pi$ for some finite $m$ -- in this case we define $k_{s_l-1}=\infty$, $m_{s_l-1}=m$.
Similarly, if this sequence is finite from the right, then either the most right point 
$M_{s_r} \in W^s_{loc}(O_1)\cap \Pi$, or
$M_{s_r} \in F_1^{-k} \circ F_{12}^{-1} (W^s_{loc}(O_2))\cap \Pi$ for some finite $k$ -- in this case we define $m_{s_r}=\infty$, $k_{s_r}=k$.

By Lemma \ref{lemkey9}, we have three possibilities.

\noindent $\bullet$ The first possibility is that $k_s$ and $m_s$ are finite and the same for all $s$  and all non-exceptional orbits, that is, $k_s=k$ and $m_s=m$ for some integers $k$ and $m$.
In this case $T_{k,m}(\Pi)$ intersects $\Pi$, which implies, by the first equation of \eqref{eq:maps:T_k,m_cross},
that $\mu =  O(|\lambda|^{k}+|\gamma|^{-m})$. By \eqref{eq:dependence0}, this implies
$\mu =  O(|\lambda_0|^{k}+|\gamma_0|^{-m})$, hence
\begin{align*}
\lambda^k\gamma^m &= \lambda_0^k (1 + O(k |\lambda_0|^k) + O(k |\gamma_0|^{-m})) \gamma_0^m((1 + O(m |\gamma_0|^{-m}) 
+ O(m |\lambda_0|^k))\\
&=\lambda_0^k \gamma_0^m (1 + o(1)_{k,m\to\infty}) + o(1)_{k,m\to\infty}.
\end{align*}
Thus, using \eqref{eq:lambdagamma}, we obtain
\begin{align*}
ab\lambda^k\gamma^m &=a_0b_0\lambda_0^k \gamma_0{m}(1+o(1)_{k,m\to\infty})+ o(1)_{k,m\to\infty}\\
&=a_0b_0\gamma_0^{\frac{mq-kp}{q}}(1+o(1)_{k,m\to\infty})+o(1)_{k,m\to\infty}.
\end{align*}
By \eqref{rare1}, we obtain that $|ab\lambda^k\gamma^m|$ stays bounded away from $1$. Since this is, up to small corrections,
the derivative $d\bar X/dX$ in \eqref{eq:maps:T_k,m_cross}, and since we have a strong contraction in $Z$ and a strong expansion in $Y$, the 
hyperbolicity of $T_{k,m}$ follows, if $k,m$ are sufficiently large and the neighborhood $U$ is sufficiently small.

By the hyperbolicity of $T_{k,m}$, it can have only one fixed point and it is the only orbit of $T_{k,m}$ that never leaves $\Pi$. Thus, in the case under
consideration, we have that $M_s=M_0$ for all $s$ and the orbit $L$ of $M_0$ is a hyperbolic periodic orbit. Any other orbits in $U$ must lie in the stable or unstable manifold of $L_{1,2}$, which includes the orbits $L_1$ and $L_2$ themselves, as well as orbits of transverse intersections of $W^u(L_2)$ with $W^s(L_1)$, of $W^u(L_2)$ with $W^s(L)$ if $L$ has index $d_1$, and of $W^s(L_1)$ with $W^u(L)$ if $L$ has index $d_2$. At the same time, no orbits from $W^u(L_1)\setminus L_1$ or $W^s(L_2)\setminus L_2$ can lie entirely in $U$
in this case (as such orbits would correspond to infinite $k$ or $m$). This is in a complete agreement with the statement of the theorem: if $L$ is of index $d_1$, then the hyperbolic set $\Lambda_1$ from the statement of the theorem is the union of $L_1$, $L$, and the heteroclinic orbits corresponding to the intersection of $W^u(L)$ with $W^s(L_1)$; otherwise, we have the set $\Lambda_2$ comprised by $L_2$, $L$, and the heteroclinic orbits corresponding to the intersection of $W^u(L_2)$ with $W^s(L)$.
 
\noindent $\bullet$ The second possibility is that $m_s$ are finite 
and the same for all $s$ and all non-exceptional orbits, that is, we have $m_s=m$ while some $k_s$ are different. In this case,
$|ab \lambda^k \gamma^m|<1$ by \eqref{meq12}, where $k$ is the minimal value of $k_s$ taken over all non-exceptional orbits.
Thus, the derivative $d\bar X/dX$ in \eqref{eq:maps:T_k,m_cross}
is small. This means that the maps $T_{k_s,m_s=m}$ are all hyperbolic, with contraction in $(X,Z)$ and expansion in $Y$,
and the set of all non-exceptional orbits which are not in $W^u(L_2)$ is uniformly hyperbolic of index $d_1$. 
The union of this set and $L_1$ is the hyperbolic set $\Lambda_1$ from the statement of the theorem. Note also that the finiteness of $m$ implies that no orbit in $W^s(L_2)\setminus L_2$ can lie entirely in $U$ in this case.

\noindent $\bullet$ The last possibility is that $k_s$ are finite and the same for all $s$ and all non-exceptional orbits, that is, we have $k_s=k$ while some $m_s$ are different.
In this case, $|ab \lambda^k \gamma^m|>1$ by \eqref{keq12}, where $m$ is the minimal value of $m_s$ taken over all non-exceptional orbits.
Thus, the derivative $d\bar X/dX$ in \eqref{eq:maps:T_k,m_cross} is very large, so all the maps $T_{k_s=k,m_s}$ are hyperbolic, with expansion in $(X,Y)$ and contraction in $Z$. The set of all non-exceptional orbits which are not in $W^s(L_1)$ is uniformly hyperbolic of index $d_2$, 
and stays at a non-zero distance from $L_1$ (the union of this set with $L_2$ is the set $\Lambda_2$ from the statement of the theorem). 
The finiteness of $k$ implies that no orbit in $W^u(L_1)\setminus L_1$ can lie entirely in $U$ in this case.

In all three cases we have a complete agreement with the statement of the theorem.
\end{proof}

\section{Saddle-focus and double-focus cases. Proof of Theorem \ref{thm:sad-foc_main}}\label{sec:complex_ev}
In this section, we consider the case where at least one of the central multipliers $\lambda_{1,1}$ and $\gamma_{2,1}$ is complex (nonreal), and prove Theorem \ref{thm:sad-foc_main}. The strategy is to detect the dynamics similar to those near type-I heterodimensional cycles; namely,
we check that the rational independence conditions of Theorem \ref{thm:sad-foc_main} imply that certain collections of the first-return maps satisfy conditions of Proposition \ref{prop:blender}.

It should be pointed out that we do not create type-I cycles via bifurcations; instead, we use the rotation brought by the complex central multipliers in order to see that conditions of Proposition \ref{prop:blender} are satisfied, both for the creation of cs-blenders and cu-blenders. The only place where we study bifurcations is the result on the stabilization of cycles in the saddle-focus case.

Let us start with the saddle-focus case. Like in Section \ref{sec:firstreturn}, we derive formulas for first-return maps $T_{k,m}$ near a saddle-focus heterodimensional cycle, see Section \ref{sec:firstreturn_sf}. An immediate observation is that the formula \eqref{eq:complex:T_k,m_cross_0mu} for these maps has almost the same form as the formula \eqref{eq:maps:T_k,m_cross_0mu} in the saddle case. The only  difference is that the coefficient governing the contraction/expansion in the central coordinate now contains trigonometric functions, which gives the possibility to obtain partially hyperbolic sets with different types of central dynamics by choosing appropriate $k$ and $m$. This fact is summarized in Lemma \ref{lem:conefields2}, a counterpart of Lemma \ref{lem:conefields} in the saddle case. Then, we prove in Proposition \ref{prop:sad-foc_blender} that, without destroying the original cycle, the simultaneous
existence of both  cs- and cu-blenders follows from Proposition \ref{prop:blender}, provided that $\theta$, $\frac{\omega}{2\pi}$ and 1 are rationally independent ($\theta$ is defined in \eqref{tht} and $\omega$ is the argument of the nonreal central multiplier).

Next, we study the local stabilization of saddle-focus cycles in Section \ref{sec:sta_sad-foc}. The key is to find a replacement of Lemma \ref{lem:type1_pert:includeO_1}, which tells us when the iterates  $F_{21}\circ F^m_2\circ F_{12} (W^u_{loc}(O_1))$ and $F^{-k}_1\circ F^{-1}_{12}(W^s_{loc}(O_2))$ enter the activating domain $\Pi'$ of the blenders. Thanks to the rotation, some of the iterates automatically cross $\Pi'$ properly, depending on which central multiplier is nonreal, i.e., where the rotation happens, in $F_1$ or $F_2$. The other manifold can also cross $\Pi'$ after unfolding the heterodimensional cycle, as in the saddle case (this is the only place where we unfold the original cycle). These facts are proven in Lemmas \ref{lem:sf:includeO_1} and  \ref{lem:sf_includeO_2}.

Finally, double-focus heterodimensional cycles are dealt with in Section \ref{sec:double-focus}. The main effort in this case is to obtain the necessary formula for the first-return maps. After that, one immediately finds that it has the same form as in the saddle-case. Consequently, a counterpart to Proposition \ref{prop:sad-foc_blender} stating the coexistence of blenders follows, see Proposition \ref{prop:df_blender}. The main difference is that the required arithmetic condition becomes the rational independence among $\theta$, $\frac{\omega_1}{2\pi}$, $\theta\frac{\omega_2}{2\pi}$ and 1, where $\omega_{1,2}$ are the arguments of the nonreal central multipliers. The stabilization part is also similar to the saddle-focus case but now both $W^u(O_1)$ and $W^s(O_2)$ enter $\Pi'$ automatically due to the simultaneous rotations near $O_1$ and $O_2$, see Lemmas \ref{lem:df_includeO_1} and \ref{lem:df_includeO_2}.

\subsection{First-return maps in the saddle-focus case}\label{sec:firstreturn_sf}
We assume that
$$\lambda_{1,1}=\lambda_{1,2}^*=\lambda e^{i\omega}, \;\omega\in(0,\pi), \quad\mbox{and}\quad \gamma:=\gamma_{2,1}\quad \mbox{is real},$$
where $\lambda>|\lambda_{1,3}|$ and $|\gamma_{2,1}|<|\gamma_{2,2} |$.  As mentioned in the introduction, the other case (where $\lambda_{1,1}$ is real while $\gamma_{2,1}$ is complex) can be reduced to this one by the time reversal. The main goal of this section is to obtain a formula for the first-return maps and prove its partial hyperbolicity in Lemma \ref{lem:conefields2}. As mentioned before, the way to find the first-return maps is the same as in the saddle case.

\subsubsection{Local maps}
We use the same coordinates near $O_2$ as in the saddle case so that the formulas for the local map $F_2$ remains the same (see \eqref{eq:maps:F_2^m}). Let us now introduce coordinates $(x_1,x_2,y,z)\in \mathbb{R}\times\mathbb{R}\times\mathbb{R}^{d_1}\times\mathbb{R}^{d-d_1-2}$ near $O_1$ such that the local map $F_1$ takes the form (see \citep[Lemmas 5 and 6]{GST08})
\begin{equation}\label{eq:complex:F_1}
\begin{aligned}
\bar{x}_1&=\lambda x_1 \cos k\omega+\lambda x_2 \sin k\omega + g_1(x_1,x_2,y,z),\\
\bar{x}_2&=-\lambda x_1 \sin k\omega+\lambda x_2 \cos k\omega + g_2(x_1,x_2,y,z),\\
\bar{y}&=P_1 y + g_3(x_1,x_2,y,z),\qquad  \bar{z} = P_2 z +g_4(x_1,x_2,y,z),
\end{aligned}
\end{equation}
where we do not indicate the dependence on parameters for simplicity. The eigenvalues of the matrices $P_1$ and $P_2$ are $\gamma_{1,1},\gamma_{1,2},\dots,\gamma_{1,d_1}$ and $\lambda_{1,3}\dots \lambda_{1,d-d_1}$, respectively. The functions $g$ vanish along their first derivatives at the origin, and satisfy
\begin{equation*}\label{eq:complex:F_1_nlnr}\!\!\!\!\!\!
g_{1,2,4}(0,0,y,0)=0,\quad\; g_3(x_1,x_2,0,z)=0, \quad\; g_{1,2}(x_1,x_2,0,z)=0,\quad\;
\dfrac{\partial g_{1,2,4}}{\partial (x_1,x_2)}(0,0,y,0)=0,
\end{equation*}
for all sufficiently small $x_1,x_2,y$ and $z$. Similar to the saddle case, in these coordinates the local manifolds $W^s_{loc}(O_1)$ and $W^u_{loc}(O_1)$ are given by $\{y=0\}$ and $\{x_1=0,x_2=0,z=0\}$; the leaves of the strong-stable foliation have the form $\{(x_1,x_2)=const,y=0\}$;
the restriction of the map to $W^s_{loc}(O_1)$ is linear in $x$. This is the same coordinate system as described in Section \ref{sec:intro2}.

By Lemma 7 of \citep{GST08}, for any point $(x_1,x_2,y,z)$ in $U_{01}$, we have $(\tilde x_1,\tilde x_2,\tilde y,\tilde z)=F^k_1(x_1,x_2,y,z)$ if and only if 
\begin{equation}\label{eq:complex:F^k_1}
\begin{aligned}
\tilde x_{1}&=\lambda^kx_1\cos k\omega+\lambda^k x_2\sin k\omega+ p_1(x_1,x_2,\tilde y,z),\\
\tilde x_{2}&=-\lambda^kx_1\sin k\omega+\lambda^k x_2\cos k\omega+ p_2(x_1,x_2,\tilde y,z),\\
y&=p_3(x_1,x_2,\tilde y,z),\qquad \tilde z= p_4(x_1,x_2,\tilde y,z),
\end{aligned}
\end{equation}
where
\begin{equation*}\label{eq:complex:F_1^k:derivatives}
\|p_{1,2,4}\|_{C^1}=o(\lambda^k), \qquad\qquad \|p_3\|_{C^1}=o(\hat\gamma^{-k}),
\end{equation*}
for some constant $\hat\gamma\in (1, |\gamma_{1,1}|)$, and these estimates are uniform for all systems $C^2$-close to $f$.

\subsubsection{Transition maps}
Let us now define the transition maps $F_{12}$ and $F_{21}$ from a neighborhood of $M^-_1=(0,0,y^-,0)$ to a neighborhood of 
$M^+_2=(0,v^+,0)$ and, respectively, from a neighborhood of $M^-_2=(u^-,0,w^-)$ to a neighborhood of
$M^+_1=(x_1^+,x_2^+,0,z^+)$. Arguing like in the saddle case, by the first part of condition C1 (that $F^{-1}_{12}(W_{loc}^{sE}(O_2))\pitchfork W^{u}_{loc}(O_1)$ at $M^-_1$), the transition map $F_{12}: (\tilde x_{1},\tilde x_{2},\tilde y,\tilde z)\mapsto (u,v,w)$ can be written as follows (the dots refer to the second and higher order terms in the Taylor expansion):
\begin{equation}\label{eq:complex:F_12}
\begin{aligned}
u &= \hat\mu+a_{11} \tilde x_{1}+a_{12}\tilde x_{2}+a_{13}\tilde z + a_{14}w  + \dots, \\
v -v^+&= a_{21}\tilde x_{1}+a_{22}\tilde x_{2}+a_{23}\tilde z + a_{24}w  + \dots, \\
\tilde y-y^- &=a_{31}\tilde x_{1}+a_{32}\tilde x_{2}+a_{33}\tilde z + a_{34}w  + \dots, 
\end{aligned}
\end{equation}
where the relation between $\hat\mu$ and $\mu$ is the same as in the saddle case and is given by \eqref{mmoh}. 

Recall the second part of C1 that $F_{12}(W_{loc}^{uE}(O_1))$ intersects $W^{s}_{loc}(O_2)$ transversely at $M^+_2$. When the dimension of the map $F_{12}$ is large than three (i.e., when the system has dimension higher than three if it is a diffeomorphism or higher than four if it is a flow), this condition means that the common directions shared by $\D F_{12}(\mathcal{T}_{M^-_1}W^{uE}_{loc}(O_1))\cap \mathcal{T}_{M^+_2}W^s_{loc}(O_2)$ is at most one-dimensional, which happens only if 
\begin{equation}\label{a11n0}
a_{11}^2+a_{12}^2\neq 0.
\end{equation}
When the dimension of $F_{12}$ equals three, the strong-stable coordinates $z$ are absent. In this case the above inequality holds automatically since $F_{12}$ is a diffeomorphism.

Similarly, condition C2 shows that the other transition map $F_{21}:(\tilde u,\tilde v,\tilde w)\mapsto(x_{1},x_{2},y,z)$ can be written as
\begin{equation}\label{eq:complex:F_21}
\begin{aligned}
{x}_1-x^+_1&= b_{11} (\tilde u-u^-) + b_{12}\tilde v + b_{13}{y} + \dots, \\
{x}_2-x^+_2&= b_{21}(\tilde u-u^-) + b_{22}\tilde v + b_{23}{y} + \dots, \\
\tilde {w}-w^-&=  b_{31}(\tilde u-u^-) + b_{32}\tilde v + b_{33}{y} + \dots,\\
{z}-z^+&= b_{41}(\tilde u-u^-) + b_{42}\tilde v + b_{43}{y} + \dots ,
\end{aligned}
\end{equation}
with $b^2_{11}+b^2_{21}\neq 0$.  By rotating the $x$-coordinates to a small angle, if necessary, we can assume 
\begin{equation}\label{b11n0}
b_{11}\neq 0.
\end{equation}
Thus, the above formula can be rewritten as
\begin{equation}\label{f21crc_complex}
\begin{aligned}
\tilde u - u^- &= b_{11}^{-1} ( x_1 -x^+_1 - b_{13} y) + O(\|\tilde v\|+ ( x_1-x_1^+)^2+ y^2), \\
 x_2 - x^+_2 &= b_{21}b_{11}^{-1}(x-x^+_1-b_{13}y)+b_{23}y+ O(\|\tilde v\|+(x_1 -x_1^+)^2+ y^2), \\
\tilde w - w^- &= O(| x_1 -x^+_1|+\|\tilde v\|+\| y\|), \\
 z - z^+ &= b_{41}b_{11}^{-1}(x-x^+_1-b_{13}y)+b_{43}y  + O(\|\tilde v\|+ ( x_1-x_1^+)^2+ y^2).
\end{aligned}
\end{equation}

\subsubsection{First-return maps and cone field lemma}
Now we can find a formula for the first-return map in the same way as in the saddle case. Namely, combining \eqref{eq:complex:F^k_1} and \eqref{eq:complex:F_12} yields the analogue of \eqref{eq:maps:112}:
\begin{equation}\label{122a}
\begin{aligned}
u &=\hat\mu + a_{11}\lambda^k x_1\cos k\omega +a_{11}\lambda^k x_2\sin k\omega 
-a_{12}\lambda^k x_1\sin k\omega +a_{12}\lambda^k x_2\cos k\omega + \tilde h_1(x_1,x_2,z,w)  \\
y &=\tilde h_2(x_1,x_2,z,w),\qquad v- v^+ =  \tilde h_3(x_1,x_2,z,w),
\end{aligned}
\end{equation}
where 
$$\tilde h_1(x_1,x_2,z,w)=O(\|w\|)+o(\lambda^k), \qquad \tilde h_2(x_1,x_2,z,w)=o(\hat\gamma^{-k}), \qquad  \tilde h_3(x_1,x_2,z,w)=O(\|w\|+|\lambda|^k).$$
Combining \eqref{eq:maps:F_2^m} and \eqref{f21crc_complex} yields the analogue of \eqref{eq:maps:221}
\begin{equation}\label{eq:221_complex}
\begin{aligned}
u &= \gamma^{-m}(u^- + b_{11}^{-1} ( x_1 -x_1^+ - b_{13}  y + \hat h_{01}( x_1- x_1^+, y)))  + \hat h_1( x_1,v, y), \\
 x_2 - x^+_2 &=b_{21}b_{11}^{-1}(x_1-x^+_1-b_{13}y)+b_{23}y +\hat h_{02}( x_1-x_1^+, y)+\hat h_2( x_1,v, y), \\
 z - z^+ &= b_{41}b_{11}^{-1}(x_1-x^+_1-b_{13}y)+b_{43}y+ \hat h_{03}( x_1-x_1^+, y)   + \hat h_3( x_1,v, y),\\
w &= \hat h_4( x_1, v, y),
\end{aligned}
\end{equation}
where
$$ \hat h_{0i}( x_1-x_1^+,\bar y)=O( (x_1-x_1^+)^2+ y^2)\;\; (i= 1,2,3), \qquad
\hat h_{1,4}( x_1, v, y)=o(\gamma^{-m}), \qquad \hat h_{2,3}( x_1, v, y)=o(\hat\lambda^m).
$$
Finally, combining \eqref{122a} and \eqref{eq:221_complex} with renaming $x_1,x_2,y,z$ by $\bar x_1,\bar x_2,\bar y,\bar z$ in \eqref{eq:221_complex}, we obtain the following form for the first-return map $T_{k,m}:=F_{21}\circ F_2^m \circ F_{12} \circ F^k_1:(x_1,x_2,y,z)\mapsto(\bar x_1,\bar x_2,\bar y,\bar z)$:
\begin{equation}\label{eq:complex:F_k,m_cross}
\begin{aligned}
\bar x_1 - x^+_1 =& b_{11} \gamma^m\hat\mu  + b_{11} \lambda^k\gamma^m\left((a_{11} x_1+a_{12}x_2)\cos k\omega + (a_{11} x_2-a_{12}x_1)\sin k\omega \right)- b_{11}u^-+ b_{13}\bar y \\
&- \hat h_{01}(\bar x_1-x^+_1,\bar y)+\gamma^m h_1(x_1,\bar x_1,x_2,\bar y,z),\\
\bar x_2 - x^+_2 = &b_{21}b_{11}^{-1}(\bar x_1-x^+_1-b_{13}\bar y)+b_{23}\bar y +\hat h_{02}( \bar x_1-x_1^+, y)+ h_2(x_1,\bar x_1,x_2,\bar y,z),\\
\bar z - z^+ =&  b_{41}b_{11}^{-1}(\bar x_1-x^+_1-b_{13}\bar y)+b_{43}\bar y+\hat h_{03}(\bar x_1-x^+_1,\bar y)+ h_3(x_1,\bar x_1,x_2,\bar y,z)\\
y =&h_4(x_1,\bar x_1,x_2,\bar y,z),
\end{aligned}
\end{equation}
where
\begin{equation*}
h_{1}=o(\lambda^k) + o(\gamma^{-m}),\quad
h_{2,3} = o(\hat\lambda^m),\quad
h_4 = o(\hat\gamma^{-k}).
\end{equation*}

We further do computations with the first-return map $T_{k,m}$ only at $\mu=0$, so we will omit the term $b \gamma^m\hat\mu$ in
formula \eqref{eq:complex:F_k,m_cross}. Make the coordinate transformation
\begin{equation}\label{eq:complex:coortrans}
\begin{array}{l}
X_1 = x_1 - x^+_1 -b_{13}y,\qquad Y=y,\\
X_2 =\delta^{-\frac{1}{2}}\; (x_2 - x^+_2 - b_{21}b_{11}^{-1}( x_1-x^+_1-b_{13}y)-b_{23}y  -\hat h_{02}( x_1-x_1^+, y)),\\
Z = z- z^+ - b_{41}b_{11}^{-1}(x-x^+_1-b_{13}y)-b_{43}y-\hat h_{03}(x_1-x^+_1, y).
\end{array}
\end{equation}
The first-return map acquires the form
\begin{equation}\label{eq:complex:T_k,m_cross}
\begin{aligned}
\bar X_1 &= A_{k,m} X_1 +  B_{k,m} + \delta^{\frac{1}{2}} O(\lambda^k\gamma^m)  X_2+\hat\phi_{01}(\bar X_1,\bar Y)  + \gamma^m\hat\phi_1(X_1,\bar X_1, X_2,\bar Y,Z),\\
\bar X_2 &= \delta^{-\frac{1}{2}}\; \hat\phi_2(X_1,\bar X_1, X_2,\bar Y,Z),\\
Y &=\hat\phi_3(X_1,\bar X_1, X_2,\bar Y,Z), \qquad
\bar Z = \hat\phi_4(X_1,\bar X_1, X_2,\bar Y,Z),
\end{aligned}
\end{equation}
where
\begin{equation}\label{eq:maps:complex:nonlinearterms}
\begin{aligned}
&\hat\phi_{01}=O(\bar X^2_1+\bar Y^2),\qquad
\|\hat\phi_1\|_{C^1}=o(\lambda^k)\!+o(\gamma^{-m}), \\
 &\|\hat\phi_{2,4}\|_{C^1}=o(\hat\lambda^m),\qquad
\|\hat\phi_{3}\|_{C^1}=o(\hat\gamma^{-k}),
\end{aligned}
\end{equation}
and 
\begin{equation}\label{akmbkm}
\begin{aligned}
A_{k,m}&=\lambda^k\gamma^m\left((a_{11} b_{21} -a_{12}b_{11})\sin k\omega+ (a_{11}b_{11}+a_{12}b_{21})\cos k\omega  \right)\\
&=\lambda^k\gamma^m  A \sin(k\omega + \eta_1),\\ 
B_{k,m}&= \lambda^k\gamma^m b_{11}(a_{11}x^+_2 -a_{12}x^+_1)\sin k\omega  +b_{11}(a_{11}x^+_1 +a_{12}x^+_2)\cos k\omega )-b_{11} u^-\\ 
&=\lambda^k\gamma^m B\sin(k\omega + \eta_2)- b_{11} u^-.
\end{aligned}
\end{equation}
Note the factor $\delta^{\frac{1}{2}}$ in front of $X_2$ in (\ref{eq:complex:T_k,m_cross}) which appears because we scale $x_2$ to
$\delta^{\frac{1}{2}}$ in \eqref{eq:complex:coortrans}. Note also that the coefficients
$$A =  \sqrt{(a_{11}^2 +a_{12}^2) (b_{11}^2+b_{21}^2)} \quad\mbox{and}\quad  B = |b_{11}| \|x^+\| \sqrt{a_{11}^2+a_{12}^2}$$
are non-zero by (\ref{a11n0}), (\ref{b11n0}) and because $x^+=(x^+_1,x^+_2)\neq 0$ due to the non-degeneracy condition C3.
The phases $\eta_{1,2}$ are defined by
\begin{align*}
A\sin\eta_1 = a_{11} b_{11} +a_{12}b_{21}, \qquad A\cos\eta_1 =a_{11} b_{21} - a_{12}b_{11}, \\
\frac{B}{b_{11}}\sin\eta_2 =  a_{11} x^+_1+a_{12}x^+_2,\qquad \frac{B}{b_{11}}\cos\eta_2 = a_{11} x^+_2-a_{12}x^+_1,
\end{align*}
which, together with $b_{11}/b_{21}\neq x^+_1/x^+_2$ due to the non-degeneracy condition C4.2 (see \eqref{eq:complex:F_21}), imply that 
\begin{equation}\label{tanet12}
\tan\eta_1\neq\tan\eta_2.
\end{equation}

Similarly to the saddle case, we consider pairs $(k,m)$ such that the maps $T_{k,m}$ take 
\begin{equation}\label{eq:complex:domain}
\Pi=[-\delta,\delta]\times[-\delta,\delta]\times[-\delta,\delta]^{d_1}\times [-\delta,\delta]^{d-d_1-1}
\end{equation}
into itself, which implies that 
\begin{equation}\label{eq:complex:lbdgm}
\lambda^k\gamma^m B\sin(k\omega + \eta_2) = b_{11} u^- + O(\delta).
\end{equation}
We will consider only such $k$ for which $\sin(k\omega + \eta_2)$ stays bounded away from $0$. 
This, along with \eqref{eq:complex:lbdgm}, implies that $\lambda^k\gamma^m$ is uniformly bounded. In particular, the term
$\gamma^m \hat\phi_1$ in \eqref{eq:complex:T_k,m_cross} tends to zero as $k,m\to\infty$. This allows to express $\bar X_1$ as a function
of $(X_1,X_2,Z,\bar Y)$ from the first equation of \eqref{eq:complex:T_k,m_cross}, and thus get rid of the dependence on $\bar X_1$ in the
right-hand side of \eqref{eq:complex:T_k,m_cross}. Thus, we can rewrite formula \eqref{eq:complex:T_k,m_cross} for $T_{k,m}$ as
\begin{equation}\label{eq:complex:T_k,m_cross_0mu}
\begin{aligned}
\bar X_1 &=  A_{k,m} X_1 + B_{k,m} +\phi_{1}(X_1, X_2,Z, \bar Y),\\
\bar X_2 &= \phi_2(X_1,X_2, Z, \bar Y),\quad \bar Z = \phi_4(X_1,X_2, Z, \bar Y),\\
Y &=\phi_3(X_1,X_2,Z,\bar Y),
\end{aligned}
\end{equation}
where
\begin{equation}\label{eq:complex:nonli_0mu}
\begin{array}{ll}
\phi_{1}=O(\delta^{\frac{3}{2}})+o(1)_{k,m\to\infty},\qquad & \dfrac{\partial\phi_{1}}{\partial(X_1,X_2,\bar Y,Z)}= O(\delta^{\frac{1}{2}})+o(1)_{k,m\to\infty},\\
  \|\phi_{2,4}\|_{C^1}=o(\hat\lambda^{m}),\qquad &\|\phi_3\|_{C^1} =o(\hat\gamma^{-k}).
\end{array}
\end{equation}
We have the above estimate for $\phi_2$ because in the $\bar X_2$-equation in \eqref{eq:complex:T_k,m_cross} the coefficient $\delta^{-\frac{1}{2}}$ can be absorbed by $\hat \phi_2$. Indeed, we always first take $\delta$ sufficiently small and then take $k,m$ sufficiently large.

Formula \eqref{eq:complex:T_k,m_cross_0mu} represents the first-return maps in the form used in Proposition \ref{prop:blender}
(where one should choose $(X_2,Z)$ as a new $Z$-variable). We will check in Section \ref{sec:sad-foc_blender} that these maps indeed
satisfy conditions of this proposition, thus establishing the existence of the blenders. Before doing that, we further restrict the choice of $(k,m)$
by the requirement that $\sin(k\omega + \eta_1)$ stays bounded away from zero. Then, the constant 
$A_{k,m}$ in \eqref{eq:complex:T_k,m_cross_0mu} stays bounded away from zero, and we notice that
the map \eqref{eq:complex:T_k,m_cross_0mu} assumes, upon setting $X=X_1$ and $Z^{new}=(X_2,Z)$,
the same form as \eqref{eq:maps:T_k,m_cross_0mu} (with the only difference being a slightly worse estimate for $\phi_1$, which does not affect
the further results). This gives us the following analogue of Lemma \ref{lem:conefields} for the saddle-focus case:
\begin{lem}\label{lem:conefields2} Let $\mu=0$.
Given any $K>0$, one can choose $\delta$ sufficiently small, such that for all sufficiently large $(k,m)$ such that
$\sin(k\omega+\eta_1)$ and $\sin(k\omega+\eta_2)$ stay bounded away from zero, the cone fields on $\Pi$
\begin{eqnarray}
&&\mathcal{C}^{cu} = \{(\Delta X_1,\Delta X_2,\Delta Y,\Delta Z): \|(\Delta X_2, \Delta Z)\| \leq K (|\Delta X_1|+\|\Delta Y\|)\},\label{eq:conefields2:u}\\
&&\mathcal{C}^{uu} = \{(\Delta X_1,\Delta X_2,\Delta Y,\Delta Z): \max\{|\Delta X_1|,\|(\Delta X_2, \Delta Z)\|\}\leq K  \|\Delta Y\|\},\label{eq:conefields2:uu}
\end{eqnarray}
are forward-invariant in the sense that if a point $M\in \Pi$ has its image $\bar M=T_{k,m}(M)$ in $\Pi$, then the cone at $M$ 
is mapped into the cone at $\bar M$ by $\D T_{k,m}$; the cone fields
\begin{eqnarray}
&&\mathcal{C}^{cs} = \{(\Delta X_1,\Delta X_2,\Delta Y,\Delta Z): \|\Delta Y\|\leq K (|\Delta X_1|+\|(\Delta X_2, \Delta Z)\|)\},\label{eq:conefields2:s}\\
&&\mathcal{C}^{ss}=\{(\Delta X_1,\Delta X_2,\Delta Y,\Delta Z):\max\{|\Delta X_1|,\|\Delta Y\|\}\leq K  \|(\Delta X_2, \Delta Z)\|\},\label{eq:conefields2:ss}
\end{eqnarray}
are backward-invariant in the sense that if a point $\bar{M}\in\Pi$ has its preimage $M=T^{-1}_{k,m}(\bar M)$ in $\Pi$, 
then the cone at $\bar{M}$ is mapped into the cone at $M$ by $\D T^{-1}_{k,m}$. Moreover, vectors in $\mathcal{C}^{uu}$ and, if $|A_{k,m}|>1$, also in $\mathcal{C}^{cu}$ are expanded by $\D T_{k,m}$; vectors in $\mathcal{C}^{ss}$ and, if $|A_{k,m}|<1$, also in $\mathcal{C}^{cs}$ are contracted by $\D T_{k,m}$.
\end{lem}

\subsection{Coexistence of mutually activating blenders}\label{sec:sad-foc_blender}
Here we prove the first part of Theorem \ref{thm:sad-foc_main} for the saddle-focus case. 

\begin{prop}\label{prop:sad-foc_blender}
If $\theta$, $\frac{\omega}{2\pi}$ and $1$ are rationally independent, then, at $\mu= 0$, there exist, arbitrarily close to the heterodimensional cycle $\Gamma$, a cs-blender 
$\Lambda^{cs}$ with an activating pair $(\Pi',\mathcal{C}^{ss})$ and a cu-blender $\Lambda^{cu}$ with an activating pair $(\Pi',\mathcal{C}^{uu})$ such that $W^s(\Lambda^{cs})\pitchfork W^u(\Lambda^{cu})\neq \emptyset$ and the two blenders activate each other, namely, $W^u(\Lambda^{cs})$ contains a piece crossing $\Pi'$ properly with respect to $\mathcal{C}^{uu}$ and $W^s(\Lambda^{cu})$ contains a piece crossing $\Pi'$ properly with respect to $\mathcal{C}^{ss}$. Here
\begin{equation}\label{eq:Pi'_complex}
\Pi'=[-q\delta,q\delta]\times[-\delta,\delta]\times[-\delta,\delta]^{d_1}\times [-\delta,\delta]^{d-d_1-2}
\end{equation}
for some fixed $q\in(0,1)$, and $\mathcal{C}^{ss},\mathcal{C}^{uu}$ are given by Lemma \ref{lem:conefields2} with some sufficiently small $K$.
\end{prop}
This proposition immediately leads to robust heterodimensional dynamics involving the two blenders, see the discussion after Definition \ref{defi:actcondition}.
\begin{proof}
The rational independence condition implies that the set
\begin{equation*}\label{eq:sf_bld:7}
\{(-k\theta+m,k\dfrac{\omega}{2\pi}-p)\}_{k,m\in \mathbb N, p\in \mathbb Z}
\end{equation*}
is dense in $\mathbb R^2$. So, given any $(s,t)\in \mathbb R^2$, one can find a sequence $\{(k_n,m_n,p_n)\}$ with $k_n,m_n\to\infty$  such that
$m_n$ are even and
\begin{equation}\label{knmnpn}
-k_n\theta + m_n\to t \quad\mbox{and}\quad k_n\dfrac{\omega}{2\pi}-p_n \to s.
\end{equation}
In fact, we fix a sufficiently large value of $t$ and let $s$ depend on $n$ so that
$$k_n\dfrac{\omega}{2\pi}-p_n - s_n\to 0.$$
We take the sequence $\{s_n\}$ dense in a sufficiently small closed interval $\Delta$ such that
the values
\begin{equation}\label{nov1}
|\gamma|^t B\sin(2\pi s_n+\eta_2) - b_{11}u^-
\end{equation}
are dense in a small interval around zero. In particular, this difference is $O(\delta)$ . Moreover, by \eqref{tanet12}, we can always choose $\Delta$ 
such that, for some constants $C_1$ and $C_2$,
\begin{equation}\label{eq:sf_bld:5s1}
0<C_1<\left|\dfrac{b_{11}u^- A \sin(2\pi s + \eta_1)}{B\sin(2\pi s+\eta_2)}\right|<C_2<1
\end{equation}
for all $s\in\Delta$, or such that, for some constants $C_3$ and $C_4$,
\begin{equation}\label{eq:sf_bld:5s2}
1<C_3<\left|\dfrac{b_{11}u^- A \sin(2\pi s + \eta_1)}{B\sin(2\pi s+\eta_2)} \right|<C_4<\infty
\end{equation}
for all $s\in\Delta$.

For the corresponding sequences $\{(k_n,m_n,p_n)\}$, the values of $\sin( k_n\omega + \eta_1)$ and $\sin( k_n \omega + \eta_2)$ get bounded away from zero for all sufficiently large $n$. Since $m_n$ are even, we have $\gamma^{m_n}=|\gamma|^{m_n}$ and, since $\lambda>0$, we have 
$\lambda^{k_n}\gamma^{m_n}=|\gamma|^{m_n-k_n\theta}\to |\gamma|^t$. It then follows from \eqref{nov1} that condition \eqref{eq:complex:lbdgm} is satisfied, and hence the first-return maps $T_n:=T_{k_n,m_n}$ can be represented in the form \eqref{eq:complex:T_k,m_cross_0mu} with estimates in \eqref{eq:complex:nonli_0mu}, and the cone lemma (Lemma \ref{lem:conefields2}) holds.

The coefficients $B_{k_n,m_n}$ given by \eqref{akmbkm} are dense in a small interval around zero, and, by passing to a subsequence if necessary, we have that $B_{k_n,m_n}$ lie in this interval for all $n$. Since by \eqref{akmbkm} one can write 
$$A_{k_n,m_n}=\dfrac{(B_{k_n,m_n}+b_{11}u^-)A\sin(k_n\omega+\eta_1)}{B\sin(k_n\omega+\eta_2)},$$
inequalities \eqref{eq:sf_bld:5s1} and \eqref{eq:sf_bld:5s2} imply that $C_1<|A_{k_n,m_n}|<C_2$ or $C_3<|A_{k_n,m_n}|<C_4$, depending on the choice of the interval $\Delta$.

Then, by setting $X=X_1$ and $Z^{new}=(X_2,Z)$, one immediately sees 
that conditions of Proposition \ref{prop:blender} are satisfied by the maps $T_n$, which gives us both
the cu- and cs- blender (by appropriate choices of the interval $\Delta$).

What remains is to show the homoclinic connection between the two blenders. Similarly to Remark \ref{rem:transintersection}, for every point $M_1\in \Lambda^{cs}$, we define
its local stable manifold as a connected piece of $W^s(\Lambda^{cs})\cap \Pi$ through $M_1$,
and the local unstable manifold as a connected piece of $W^u(\Lambda^{cs})\cap \Pi$ through $M_1$. By Lemma \ref{lem:conefields2}, the tangent space of $W^s_{loc}(M_1)$ at any point lies in the stable cone $\mathcal{C}^{cs}$, so $W^s_{loc}(M_1)$ is given by an equation
$Y=\xi^s_{M_1}(X_1,X_2,Z)$ where $\xi^s_{M_1}$ has its derivative small and is defined on $[-\delta,\delta]\times[-\delta,\delta]\times[-\delta,\delta]^{d-d_1-2}$. Similarly, the manifold $W^u_{loc}(M_1)$ is the graph of $(X_1,X_2,Z)=\xi^u_{M_1}(Y)$ defined for $Y\in[-\delta,\delta]^{d_1}$ and the tangent space of $W^u_{loc}(M_1)$ at any point lies in $\mathcal{C}^{uu}$. This means that $W^u_{loc}(M_1)$ crosses $\Pi'$ properly with respect to $\mathcal{C}^{uu}$, i.e., $\Lambda^{cs}$ activates $\Lambda^{cu}$.

For every point $M_2\in \Lambda^{cu}$, the manifold $W^u_{loc}(M_2)$ is the graph of $(X_2,Z)=\xi^u_{M_2}(X_1,Y)$ defined for $(X_1,Y)\in [-\delta,\delta]\times[-\delta,\delta]^{d}$ and its tangents lie in $\mathcal{C}^{cu}$. It immediately follows that it has a non-empty transverse intersection with $W^s_{loc}(M_1)$ for any point $M_1\in \Lambda^{cs}$, so we have $W^s(\Lambda^{cs})\pitchfork W^u(\Lambda^{cu})\neq \emptyset$. The manifold $W^s_{loc}(M_2)$ is the graph of
$(X_1,Y)=\xi^s_2(X_2,Z)$ defined for $(X_2,Z)\in[\delta,\delta]\times[-\delta,\delta]^{d-d_1-2}$, and its tangents lie in $\mathcal{C}^{ss}$. So, it crosses $\Pi'$ properly with respect to $\mathcal{C}^{ss}$, which means that $\Lambda^{cu}$ activates $\Lambda^{cs}$.
\end{proof}
\begin{rem}\label{remtrans} It follows from the arguments at the end of this proof that $W^u(\Lambda^{cs})$ intersects transversely $W^s_{loc}(O_1)$ (i.e., the manifold $\{Y=0\}$), while $W^s(\Lambda^{cs})$ intersects transversely any manifold
which crosses $\Pi$ properly with respect to $\mathcal{C}^{uu}$. Similarly,
$W^s(\Lambda^{cu})$ intersects transversely $F_{21}(W^u_{loc}(O_2))$ (i.e., the manifold $\{X_2=0,Z=0\}$ obtained by taking $m\to\infty$ in \eqref{eq:221_complex}), and
$W^u(\Lambda^{cu})$ intersects transversely any manifold
which crosses $\Pi$ properly with respect to $\mathcal{C}^{ss}$.
\end{rem}

\subsection{Local stabilization of  heterodimensional cycles in the saddle-focus case}\label{sec:sta_sad-foc}
In this section we prove the second part of Theorem \ref{thm:sad-foc_main} for the saddle-focus case. As in the saddle case, we investigate the iterates of the local invariant manifolds $W^u_{loc}(O_1)$ and $W^s_{loc}(O_2)$. 

\begin{lem}\label{lem:sf:includeO_1}
Let $q,\Pi',\mathcal{C}^{uu}$ be given by Proposition \ref{prop:sad-foc_blender}. Define the intervals
\begin{equation*}\label{eq:sf:includeO_1:0}
I^u_m=\left(\gamma^{-m} u^- -  \dfrac{1}{2}|b_{11}^{-1}\gamma^{-m}|q\delta,\gamma^{-m}u^-+ \dfrac{1}{2}|b_{11} ^{-1}\gamma^{-m}|q\delta\right).
\end{equation*}
Take sufficiently small $\delta$. For every sufficiently large $m$, if $\hat\mu\in I^u_m$, then the image $S^u_m:=F_{21}\circ F_2 \circ F_{12}(W^u_{loc}(O_1))$ is a disc of the form $(X_1,X_2,Z)=s(Y)$ for some smooth function $s$. The disc $S^u_m$ crosses the cube $\Pi'$ properly with respect to the cone field $\mathcal{C}^{uu}$.
\end{lem}

\begin{proof}
Since $W^u_{loc}(O_1)$ near $M_1^-$ has equation $(\tilde x,\tilde z)=0$,
it follows from formulas \eqref{eq:complex:F_12}, \eqref{eq:maps:F_2^m} and \eqref{eq:complex:F_21}, that the image $S^u_m$ is given by
\begin{equation*}\label{eq:sf:includeO_1:1}
\begin{aligned}
x_1 - x^+_1&= b_{11} \gamma^m\hat\mu-b_{11} u^- +b_{13} y+O((\gamma^m\hat\mu-u^-)^2+y^2)+o(1)_{m\to\infty},\\
x_2 - x^+_2&= b_{21}\gamma^m\hat\mu-b_{21}u^- +b_{23}y +O((\gamma^m\hat\mu-u^-)^2+y^2)+o(1)_{m\to\infty},\\
z - z^+&=  b_{41}\gamma^m\hat\mu-b_{41}u^- +b_{43}y +O((\gamma^m\hat\mu-u^-)^2+y^2)+o(1)_{m\to\infty},
\end{aligned}
\end{equation*}
which, after the transformation \eqref{eq:complex:coortrans}, recasts as
\begin{equation}\label{eq:sf:includeO_1:2}
\begin{aligned}
X_1&=b_{11} \gamma^m\hat\mu-b_{11} u^- +O((\gamma^m\hat\mu-u^-)^2+Y^2)+o(1)_{m\to\infty},\\
X_2&=\delta^{-\frac{1}{2}} O((\gamma^m\hat\mu-u^-)^2+Y^2)+\delta^{-\frac{1}{2}}o(1)_{m\to\infty},\\
Z&=O((\gamma^m\hat\mu-u^-)^2+Y^2)+o(1)_{m\to\infty}
\end{aligned}
\end{equation}

Since $\hat\mu\in I^u_m$, we have
$$|b_{11}\gamma^m\hat\mu-b_{11} u^-| < \dfrac{q\delta}{2} ,$$
which for $Y\in[-\delta,\delta]^{d_1}$ implies
$$|X_1|< \dfrac{q\delta}{2}+O(\delta^2)+o(1)_{m\to\infty}<q\delta
\quad\mbox{and}\quad
\|(X_2,Z)\|=O(\delta^{\frac{3}{2}})+\delta^{-\frac{1}{2}}  o(1)_{m\to\infty}< \delta,
$$
where we first take $\delta$ sufficiently small and then take $m$ sufficiently large (and we do the same whenever terms like $\delta^{-\frac{1}{2}}  o(1)_{m\to\infty}$ appear). 
This means that $S^u_m$ crosses $\Pi'$. One also finds from \eqref{eq:sf:includeO_1:2} that 
$$\dfrac{\partial (X_1,X_2,Z)}{\partial Y} =O(\delta^{\frac{1}{2}}) +\delta^{-\frac{1}{2}} o(1)_{m\to\infty},$$
which can be made sufficiently small so that the tangent spaces of $S^u_m\cap\Pi$ lie in $\mathcal{C}^{uu}$. So, the crossing is also proper with respect to $\mathcal{C}^{uu}$.
\end{proof}

\begin{lem}\label{lem:sf_includeO_2}
Let $\Pi'$ and $\mathcal{C}^{ss}$ be given by Proposition \ref{prop:sad-foc_blender}.
There exists a sequence $\{k_j\}\to\infty$ such that, at $\mu=0$, the preimage $S^s_j:=F^{-k_j}_{1}\circ F^{-1}_{12}(W^s_{loc}(O_2))$ is a disc of the form $(X_1,Y)=s(X_2,Z)$ for some smooth function $s$. The disc $S^s_j$ crosses $\Pi'$ properly with respect to the cone field $\mathcal{C}^{ss}$.
\end{lem}
\begin{proof}
Since $W^s_{loc}(O_2)$ near the point $M_2^+$ is given by $\{(u,w)=0\}$, the preimage $F^{-k}_{1}\circ F^{-1}_{12}(W^s_{loc}(O_2))$ can be found from  \eqref{122a} by setting $\hat\mu=0$ and $u=w=0$, which in coordinates \eqref{eq:complex:coortrans} is given by

\begin{equation}\label{eq:sf_includeO_2:1}
\begin{aligned}
X_1 &=-\dfrac{B\sin(k\omega+\eta_2)}{A\sin(k\omega+\eta_1)} + O(\delta^{\frac{3}{2}}) + o(1)_{k\to\infty}, \\ 
Y&=o(\hat\gamma^{-k}),
\end{aligned}
\end{equation}
where the right-hand sides are functions of $(X_2,Z)$.

The assumption that $\theta$, $\omega/2\pi$, and $1$ are rationally independent implies that $\omega/2\pi$ is irrational. As a result, one can find a sequence $\{k_j\}$ of positive integers with $k_j\to \infty$ such that
\begin{equation*}
\left|\dfrac{B\sin({k_j}\omega+\eta_2)}{A\sin({k_j}\omega+\eta_1)}\right|  <\dfrac{q\delta}{2}
\end{equation*}
for all sufficiently large $k_j$. It follows that for all sufficiently small $\delta$ and sufficiently large $k_j$, we have 
$$|X_1|<q\delta,\quad \|Y\|<\delta,\quad \dfrac{\partial (X_1,Y)}{\partial (X_2,Z)}=O(\delta^{\frac{3}{2}}) + o(1)_{k\to\infty}$$
in \eqref{eq:sf_includeO_2:1}, which completes the proof of the lemma.
\end{proof}

We can now finish the proof of Theorem \ref{thm:sad-foc_main} for the saddle-focus case. Let $\Lambda^{cs}$ and $\Lambda^{cu}$ be the cs- and cu-blenders of Proposition \ref{prop:sad-foc_blender}. By Remark \ref{remtrans}, Lemma \ref{lem:sf:includeO_1} implies that $\Lambda^{cs}$ is homoclinically related to $O_1$ when $\hat\mu \in I^u_m$; and Lemma \ref{lem:sf_includeO_2} implies that $\Lambda^{cu}$ is homoclinically related to $O_2$ at $\mu=0$, hence at all small
$\mu$, e.g. when $\hat\mu\in I^u_m$. \qed

\subsection{Double-focus case}\label{sec:double-focus}

We finish the paper by considering the case
$$\lambda_{1,1}=\lambda_{1,2}^*=\lambda e^{i\omega_1},\;\omega_1\in(0,\pi),\quad\mbox{and}\quad 
\gamma_{2,1}=\gamma_{2,2}^*=\gamma e^{ i\omega_2},\;\omega_2\in(0,\pi),$$
where $\lambda>|\lambda_{1,3}|$, $\gamma<|\gamma_{2,3}|$. We no longer need to split the heterodimensional cycle to get robust heterodimensional dynamics, so below we write formulas only for the unperturbed system (e.g. $\mu=0$).

\subsubsection{Local maps}
We use the same coordinates near $O_1$ as in the saddle-focus case, so the local map $F_1$ is given by \eqref{eq:complex:F_1} with replacing $\omega$ by $\omega_1$. Near $O_2$ we introduce coordinates $(u_1,u_2,v,w)\in \mathbb{R}\times\mathbb{R}\times\mathbb{R}^{d-d_1-1}\times\mathbb{R}^{d_1-1}$ such that the local map $F_2$ assumes the form (see \citep[Lemmas 5 and 6]{GST08})
\begin{equation*}\label{eq:df:F_2}
\begin{aligned}
\tilde{u}_1&=\gamma u_1 \cos k\omega_2+\gamma u_2 \sin k\omega_2 + \hat g_1(u_1,u_2,v,w),\\
\tilde{u}_2&=-\gamma u_1 \sin k\omega_2+\gamma u_2 \cos k\omega_2 + {\hat {g}}_2(u_1,u_2,v,w),\\
\tilde{v}&=Q_1 v + \hat g_3(u_1,u_2,v,w),\\
\tilde{w}&=Q_2 w + \hat g_4(u_1,u_2,v,w),
\end{aligned}
\end{equation*}
where the eigenvalues of the matrices $Q_1$ and $Q_2$ are $\lambda_{2,1},\dots,\lambda_{2,d-d_1-1}$ and $\gamma_{2,3}\dots \gamma_{1,d_1+1}$, respectively. Here the functions $\hat g$ satisfy
\begin{equation*}\label{eq:df:F_2_nlnr}
\begin{array}{l}
\hat g_{1,2,4}(0,0,v,0)=0,\quad\quad \hat g_3(u_1,u_2,0,w)=0, \quad\quad \hat g_{1,2}(u_1,u_2,0,w)=0,  \quad\quad  
\dfrac{\partial \hat g_{1,2,4}}{\partial (u_1,u_2)}(0,0,v,0)=0
\end{array}
\end{equation*}
for all sufficiently small $u_1,u_2,v$ and $w$. Similar to the saddle case, in these coordinates the local manifolds $W^u_{loc}(O_2)$ and $W^s_{loc}(O_2)$ are given by $\{v=0\}$ and $\{u_1=0,u_2=0,w=0\}$; the leaves of the strong-unstable foliation have the form $\{(u_1,u_2)=const,v=0\}$; and the restriction of the map to $W^u_{loc}(O_2)$ is linear in $u$. This is the same coordinate system as discussed in Section \ref{sec:intro2}.

By Lemma 7 of \citep{GST08}, the above coordinates can be chosen such that for any point $(u_{1},u_{2},v,w)\in U_{02} $, we have $(\tilde u_1,\tilde u_2,\tilde v,\tilde w)=F^m_2(u_{1},u_{2},v,w)$ if and only if 
\begin{equation}\label{eq:df:F_2^m}
\begin{aligned}
 u_{1} & =  \gamma^{-m} \tilde u_{1}\cos m\omega_2  + \gamma^{-m} \tilde u_{2}\sin m\omega_2+ q_1(\tilde u_{1},\tilde u_{2},v,\tilde w),\\
 u_{2} & =  -\gamma^{-m} \tilde u_{1}\sin m\omega_2  + \gamma^{-m} \tilde u_{2}\cos m\omega_2+ q_2(\tilde u_{1},\tilde u_{2},v,\tilde w),\\
\tilde v&= q_3(\tilde u_{1},\tilde u_{2},v,\tilde w),\\
w&= q_4(\tilde u_{1},\tilde u_{2},v,\tilde w),
\end{aligned}
\end{equation}
where
\begin{equation*}\label{eq:df:F_2^m:derivatives}
\|q_{1,2,4}\|_{C^1}=o(\gamma^{-m}), \qquad\qquad \|q_3\|_{C^1}=o(\hat\lambda^m),
\end{equation*}
for some constant $\hat\lambda\in (1, |\lambda_{2,1}|)$, and these estimates are uniform for all systems $C^2$-close to $f$.

\subsubsection{Transition maps}
We now take the heteroclinic points
$$
M^+_1=(x^+_1,x^+_2,0,z^+),\quad M^-_1=(0,0,y^-,0),\quad M^+_2=(0,0,v^+,0),\quad M^-_2=(u^-_1,u^-_2,0,w^-).
$$
By the non-degeneracy condition C3, we have $x^+\neq 0 $ and $ u^-\neq 0$. Up to a linear rotation of the coordinates $u$, we can always achieve
$$u_1^-\neq 0,$$
which will be our standing assumption.

The transition map $F_{12}:  (\tilde x_{1},\tilde x_{2},\tilde y,\tilde z)\mapsto (u_{1},u_{2},v,w)$ from a neighborhood of $M^-_1$ to a neighborhood of $M^+_2$ is given by 
\begin{equation*}
\begin{aligned}
u_1& = a'_{11} \tilde x_{1}+a'_{12}\tilde x_{2}+ a'_{13} (\tilde y-y^-)+ a'_{14}\tilde z+ \dots, \\
u_2 &= a'_{21} \tilde x_{1}+a'_{22}\tilde x_{2}+ a'_{23} (\tilde y-y^-)+ a'_{24}\tilde z+ \dots, \\
v -v^+&= a'_{31}\tilde x_{1}+a'_{32}\tilde x_{2}+a'_{33} (\tilde y-y^-)+a'_{34}\tilde z  + \dots, \\
w& =a'_{41}\tilde x_{1}+a'_{42}\tilde x_{2}+a'_{43} (\tilde y-y^-)+a'_{44}\tilde z  + \dots, 
\end{aligned}
\end{equation*}
where dots denote the second and higher order terms in the Taylor expansion. Note that here $\dim y=\dim w+\dim u_2$ and condition C1 means that $\det (a'_{23},a'_{43})\neq 0$. So, we can rewrite the formula as
\begin{equation}\label{eq:df:F_12_cross}
\begin{aligned}
u_1 &= a_{11} \tilde x_{1}+a_{12}\tilde x_{2}+a_{13}\tilde z + a_{14}u_2 +a_{15} w + \dots, \\
v -v^+&= a_{21}\tilde x_{1}+a_{22}\tilde x_{2}+a_{23}\tilde z + a_{24}u_2 +a_{25} w + \dots, \\
\tilde y-y^- &=a_{31}\tilde x_{1}+a_{32}\tilde x_{2}+a_{33}\tilde z + a_{34}u_2 +a_{35} w + \dots, 
\end{aligned}
\end{equation}
where $a^2_{11}+a^2_{12} \neq 0$ by the second part of C1. 
Similarly, it follows from condition C2 that the transition map $F_{21}: (\tilde u_{1},\tilde u_{2},\tilde v,\tilde w)\mapsto (x_1,x_2,y,z)$ from a neighborhood of $M^-_2$ to a neighborhood of $M^+_1$ is given by 
\begin{equation}\label{eq:df:F_21_cross}
\begin{aligned}
x_1-x^+_1&= b_{11} (\tilde u_{1}-u^-_1) + b_{12}\tilde v + b_{13}{y} + \dots, \\
x_2-x^+_2&= b_{21}(\tilde u_{1}-u^-_1) + b_{22}\tilde v + b_{23}{y} + \dots, \\
z-z^+&= b_{31}(\tilde u_{1}-u^-_1) + b_{32}\tilde v + b_{33}{y} + \dots ,\\
\tilde u_2-u^-_2&= b_{41}(\tilde u_{1}-u^-_1) + b_{42}\tilde v + b_{43}{y} + \dots, \\
\tilde w-w^-&=  b_{51}(\tilde u_{1}-u^-_1) + b_{52}\tilde v + b_{53}{y} + \dots,
\end{aligned}
\end{equation}
where $b^2_{21}+b^2_{12} \neq 0$. 

\subsubsection{First-return maps}
Slightly different from the previous cases, here we work with the first-return maps
\begin{equation*}\label{eq:df:frtrtn}
F_{k,m}:=F_2^m\circ F_{12}\circ F^k_1\circ F_{21}:(\tilde u_1,\tilde u_2,\tilde v,\tilde w)\mapsto(\bar{{u}}_1,\bar{{u}}_2,\bar{  v},\bar{  w})
\end{equation*}
defined in a small neighbourhood of $M^-_2$.

\begin{lem}\label{lem:df:F_km}
Suppose that $\cos m\omega_2$ and $(\cos m\omega_2+a_{14}\sin m\omega_2)$ are both bounded away from zero. Then, for a point $(\tilde u_1,\tilde u_2,\tilde v,\tilde w)$ in a small neighborhood of $M^-_2$, we have $(\bar u_1,\bar u_2,\bar v,\bar w)=F_{k,m}(\tilde u_1,\tilde u_2,\tilde v,\tilde w)$ if and only if
\begin{equation}\label{eq:df:F_km_cross}
\begin{aligned}
(\cos m\omega_2+a_{14}\sin m\omega_2)\bar{ u}_{1}
&=
\lambda^k\gamma^m (C_k(\tilde u_{1}-u^-_1)+ D_k)
+(a_{14}\cos m\omega_2-\sin m\omega_2)\bar{ u}_{2}\\
&\;\;\;\;+\lambda^k\gamma^m h_{01}(\tilde u_1-u^-_1,\tilde v)+\gamma^m h_1(\tilde u_1,\bar{ u}_2,\tilde v,\bar{ w}),\\
\tilde u_2-u^-_2 &=b_{41}(\tilde u_1-u^-_1)+b_{42}\tilde v+
h_{02}(\tilde u_1-u^-_1,\tilde v)+h_2(\tilde u_1,\bar{ u}_2,\tilde v,\bar{ w}),\\
\tilde w-w^-_2 &=b_{51}(\tilde u_1-u_1^-)+b_{52}\tilde v+
h_{03}(\tilde u_1-u^-_1,\tilde v)+h_3(\tilde u_1,\bar{ u}_2,\tilde v,\bar{ w}),\\
\bar{ v}&=h_4(\tilde u_1,\bar{ u}_2,\tilde v,\bar{ w}),
\end{aligned}
\end{equation}
where
\begin{equation}\label{ckdkeq}
\begin{array}{l}
C_k=(a_{11} b_{11} +a_{12}b_{21})\cos k\omega_1  + 
(a_{11} b_{21}-a_{12}b_{11})\sin k\omega_1,\\
D_k=(a_{11} x_1^+ +a_{12}x^+_2)\cos k\omega_1+
(a_{11} x^+_2-a_{12}x^+_1)\sin k\omega_1,
\end{array}
\end{equation}
and
$$h_{01}=O((\tilde u_1-u^-_1)^2 + \|\tilde v\|),\quad h_{0i}=O((\tilde u_1-u^-_1)^2+\tilde v^2)\quad (i=2,3),$$
$$h_1=o(\lambda^k)+o(\gamma^{-m}),
\quad h_{2,3}=o(\hat{\gamma}^{-k}),\quad 
h_4=o(\hat\lambda^m).
$$
\end{lem}
\begin{proof}

Consider first the composition $F^k_1\circ F_{21}$. Substituting the $y$-equation of \eqref{eq:complex:F^k_1} into the first three equations in \eqref{eq:df:F_21_cross}, one expresses $x_1,x_2,z$ as functions of $\tilde u_1,\tilde v,\tilde y$. Substituting these expressions together with the $y$-equation of \eqref{eq:complex:F^k_1} into the remaining equations in \eqref{eq:df:F_21_cross} leads to $\tilde u_2,\tilde w$ as functions of $\tilde u_1,\tilde v,\tilde y$. Then, combining \eqref{eq:df:F_21_cross} with the newly obtained equations for $x_1,x_2,z,\tilde u_2,\tilde w$ yields that for a point $(\tilde u_1,\tilde u_2,\tilde v,\tilde w)$ in a small neighborhood of $M^-_2$ we have $(\tilde x_1,\tilde x_2,\tilde y,\tilde z)=F^k_1\circ F_{21}(\tilde u_1,\tilde u_2,\tilde v,\tilde w)$ if and only if
\begin{equation}\label{eq:df:21}
\begin{aligned}
\tilde x_1 &= \lambda^k(x^+_1 + b_{11} (\tilde u_1-u^-_1)+b_{12}\tilde v)\cos k\omega_1
+\lambda^k(x^+_2 + b_{21} (\tilde u_1-u^-_1)+b_{22}\tilde v)\sin k\omega_1 \\
&\;\;\;\;+ \lambda^k O((\tilde u_1-u_1^-)^2 + \tilde v^2)+ \tilde h_1(\tilde u_1,\tilde v,\tilde y)
,\\
\tilde x_2 &= -\lambda^k(x^+_1 + b_{11} (\tilde u_1-u^-_1)+b_{12}\tilde v)\sin k\omega_1
+\lambda^k(x^+_2 + b_{21} (\tilde u_1-u^-_1)+b_{22}\tilde v)\cos k\omega_1 \\
&\;\;\;\;+ \lambda^kO((\tilde u_1-u_1^-)^2 + \tilde v^2)+ \tilde h_2(\tilde u_1,\tilde v,\tilde y),\\
\tilde u_2-u^-_2&=b_{41}(\tilde u_1-u_1^-)+b_{42}\tilde v+ O((\tilde u_1-u_1^-)^2 + \tilde v^2)+\tilde h_3(\tilde u_1,\tilde v,\tilde y),\\
\tilde w-w^-_2&=b_{51}(\tilde u_1-u_1^-)+b_{52}\tilde v+ O((\tilde u_1-u_1^-)^2 + \tilde v^2)+\tilde h_4(\tilde u_1,\tilde v,\tilde y),\\
\tilde z&=\tilde h_5(\tilde u_1,\tilde v,\tilde y),
\end{aligned}
\end{equation}
where $\|\tilde h_{1,2,5}\|_{C^1}=o(\lambda^k)$ and $\|\tilde h_{3,4}\|_{C^1}=o(\hat{\gamma}^{-k})$.

We proceed to find a formula for $F^m_2\circ F_{12}$. Substituting the $v$-equation in \eqref{eq:df:F_12_cross} into the $u_2$- and $w$-equations in \eqref{eq:df:F_2^m}, one obtains $u_2$ and $w$ as functions of $\tilde u_1,\tilde u_2,\tilde w,\tilde x_1,\tilde x_2,\tilde z$. Substituting these into the remaining equations, leads to $u_1$ and $\tilde v$ as functions of $\tilde u_1,\tilde u_2,\tilde w,\tilde x_1,\tilde x_2,\tilde z$. So we have the following relations:
\begin{equation*}
\begin{aligned}
 u_{1} & =  \gamma^{-m} \tilde u_{1}\cos m\omega_2  + \gamma^{-m} \tilde u_{2}\sin m\omega_2+ q'_1(\tilde u_1,\tilde u_2,\tilde w,\tilde x_1,\tilde x_2,\tilde z),\\
 u_{2} & =  -\gamma^{-m} \tilde u_{1}\sin m\omega_2  + \gamma^{-m} \tilde u_{2}\cos m\omega_2+ q'_2(\tilde u_1,\tilde u_2,\tilde w,\tilde x_1,\tilde x_2,\tilde z),\\
\tilde v&= q'_3(\tilde u_1,\tilde u_2,\tilde w,\tilde x_1,\tilde x_2,\tilde z),\\
w&= q'_4(\tilde u_1,\tilde u_2,\tilde w,\tilde x_1,\tilde x_2,\tilde z),
\end{aligned}
\end{equation*}
where $\|q'_{1,2,4}\|_{C^1}=o(\gamma^{-m})$ and $\|q'_3\|_{C^1}=o(\hat\lambda^m)$.

Now, since $\cos m\omega_2$ is bounded away from zero, one can express $\tilde u_1$ as a function of $u_1,\tilde  u_2,v,\tilde w$ from the $u_1$-equation. Substituting the result into the equations for $u_2,\tilde v$ and $w$, we obtain
\begin{equation}\label{eq:df:F_2^m_new}
\begin{aligned}
\gamma^{-m} \tilde u_{1}\cos m\omega_2& =  u_1   - \gamma^{-m} \tilde u_{2}\sin m\omega_2+ q''_1( u_{1},\tilde u_2,\tilde w,\tilde x_1,\tilde x_2,\tilde z),\\
 u_{2} & =  -\dfrac{ u_1\sin m\omega_2}{\cos m\omega_2} + 
 \dfrac{\gamma^{-m} \tilde u_2}{\cos m\omega_2}  + q''_2(u_{1},\tilde u_2,\tilde w,\tilde x_1,\tilde x_2,\tilde z),\\
\tilde v&= q''_3( u_{1},\tilde u_2,\tilde w,\tilde x_1,\tilde x_2,\tilde z),\\
w&= q''_4( u_{1},\tilde u_2,\tilde w,\tilde x_1,\tilde x_2,\tilde z),
\end{aligned}
\end{equation}
where $\|q''_{1,2,4}\|_{C^1}=o(\gamma^{-m})$ and $\|q''_3\|_{C^1}=o(\hat\lambda^m)$.

Recall that by assumption $\cos m\omega_2 + a_{14}\sin m\omega_2$ is also bounded away from zero. Then, substituting the above $u_2$-equation into the $u_1$-equation in \eqref{eq:df:F_12_cross} yields
\begin{equation*}\label{eq:df:u_1}
\begin{aligned}
\left(1+\dfrac{a_{14}\sin m\omega_2}{\cos m\omega_2}\right)u_1&=a_{11} \tilde x_{1}+a_{12}\tilde x_{2}+a_{13}\tilde z + 
 \dfrac{a_{14}\gamma^{-m} \tilde u_2}{\cos m\omega_2} +a_{15} w  \\
 &\;\;\;\;+ O(\tilde x_1^2+\tilde x_2^2+\tilde z^2 + \tilde u_2^2+w^2).
\end{aligned}
\end{equation*}
Combining this with \eqref{eq:df:F_2^m_new} and the $\tilde y$-equation in \eqref{eq:df:F_12_cross}, yields that for a point $(\tilde x_1,\tilde x_2,\tilde y,\tilde z)$ in a small neighborhood of $M^-_1$ we have $(\tilde u_1,\tilde u_2,\tilde v,\tilde w)=F^m_2\circ F_{12}(\tilde x_1,\tilde x_2,\tilde y,\tilde z)$ if and only if
\begin{equation}\label{eq:df:12}
\begin{aligned}
\gamma^{-m}(\cos m\omega_2 +a_{14}\sin m\omega_2)\tilde u_1
& =
\gamma^{-m}(a_{14}\cos m\omega_2 -\sin m\omega_2)\tilde u_2+a_{11}\tilde x_1+a_{12}\tilde x_2+a_{13}\tilde z  \\
&\;\;\;\;+O(\tilde x_1^2+\tilde x_2^2+\tilde z^2)
+\hat h_1(\tilde x_1,\tilde x_2,\tilde z,\tilde u_2,\tilde w),\\
\tilde y-y^- &=O(|\tilde x_1|+|\tilde x_2|+\|\tilde z\| )
+\hat h_2(\tilde x_1,\tilde x_2,\tilde z,\tilde u_2,\tilde w),\\
\tilde v&=\hat h_3(\tilde x_1,\tilde x_2,\tilde z,\tilde u_2,\tilde w),
\end{aligned}
\end{equation}
where
$
\|\hat h_1\|_{C^1}=o(\gamma^{-m}),\quad \|\hat h_2\|_{C^1}=O(\gamma^{-m}),\quad \|\hat h_3\|_{C^1}=o(\hat\lambda^m).$

Replacing $\tilde u_1,\tilde u_2,\tilde v,\tilde w$ by $\bar{{u}}_1,\bar{{u}}_2,\bar{  v},\bar{  w}$ in the above formula, and combining it with \eqref{eq:df:21} yields \eqref{eq:df:F_km_cross}.
\end{proof}

In formula \eqref{eq:df:F_km_cross}, write 
\begin{equation}\label{eq:df:T_k,m_cross_const0}
C_k=C\sin(k\omega_1 +\eta_1) \quad \mbox{and} \quad D_k=D\sin(k\omega_1+\eta_2)
\end{equation}
and
\begin{equation}\label{eq:df:T_k,m_cross_const1}
\begin{aligned}
\cos m\omega_2+a_{14}\sin m\omega_2&=\sqrt{1+a_{14}^2}\cos(m\omega_2+\eta_3),\\
 a_{14}\cos m\omega_2-\sin m\omega_2&= -\sqrt{1+a_{14}^2}\sin(m\omega_2+\eta_3).
\end{aligned}
\end{equation} 
Since $a^2_{11}+a^2_{12} \neq 0$, $b^2_{21}+b^2_{12} \neq 0$, and $x^+\neq 0$ by conditions C1 - C3, we have
$$C\neq 0 \quad\mbox{and}\quad D\neq 0.$$
Moreover, by condition C4.2 and \eqref{ckdkeq}, 
\begin{equation}\label{eqtan12}
\tan\eta_1\neq\tan\eta_2.
\end{equation}

Introduce the coordinate transformation
\begin{equation}\label{eq:df:coortrans}
\begin{array}{l}
U_1 = \tilde u_1-u^-_1,\quad \delta^{\frac{1}{2}}U_2 = \tilde u_2-u_2^--b_{41}(\tilde u_1-u^-_1)-b_{42} \tilde v,\\
\delta^{\frac{1}{2}}V=\tilde  v,\quad W=\tilde w-w^--b_{51}(\tilde u_1-u^-_1)-b_{52}\tilde v,
\end{array}
\end{equation}
and consider the restriction of $F_{k,m}$ to
\begin{equation}\label{eq:df:domain}
\Pi=[-\delta,\delta]\times[-\delta,\delta]\times[-\delta,\delta]^{d-d_1-1}\times [-\delta,\delta]^{d_1-1}
\end{equation}
in the new coordinates. Then, formula \eqref{eq:df:F_km_cross}, along with \eqref{eq:df:T_k,m_cross_const0} and \eqref{eq:df:T_k,m_cross_const1}, implies that for a point $(U_1,U_2,V,W)\in\Pi$ we have $(\bar U_1,\bar U_2,\bar V,\bar  W)=F_{k,m}(U_1,U_2,V,W)$ if and only if
\begin{equation}\label{eq:df:T_km_cross}
\begin{array}{l}
(1 + b_{41}\tan(m\omega_2+\eta_3)) \bar U_1 
= 
 \dfrac{\lambda^k\gamma^m C\sin(k\omega_1+\eta_1)U_1+\lambda^k\gamma^m D\sin(k\omega_1+\eta_2)}{\sqrt{1+a_{14}^2}\cos(m\omega_2+\eta_3)}-u^-_1 \\[10pt]

\qquad\qquad\qquad\qquad\qquad\qquad\qquad\quad 
-\tan(m\omega_2+\eta_3)
u^-_2+ \phi_1(U_1,\bar U_2,V,\bar W),\\[15pt]

U_2=\phi_2(U_1,\bar U_2,V,\bar W),\qquad
\bar{V}=\phi_3(U_1,\bar U_2,V,\bar W),\qquad
W=\phi_4(U_1,\bar U_2,V,\bar W),
\end{array}
\end{equation}
where
\begin{equation}\label{eq:df:nonlinear}
\begin{array}{l}
\phi_1= \lambda^k\gamma^m(O(\delta^{\frac{3}{2}})+ o(1)_{k,m\to\infty})+O(\delta^{\frac{3}{2}}) +o(1)_{k,m\to\infty},\\
\dfrac{\partial \phi_1}{\partial (U_1,\bar U_2,V,\bar W)}= \lambda^k\gamma^m(O(\delta^{\frac{1}{2}})+ o(1)_{k,m\to\infty}) +O(\delta^{\frac{1}{2}})+o(1)_{k,m\to\infty}\\
\phi_{2,4}=O(\delta^{\frac{3}{2}})+o(\hat\gamma^{-k}),\quad
\dfrac{\partial \phi_{2,4}}{\partial (U_1,\bar U_2,V,\bar W)}=O(\delta^{\frac{1}{2}})+o(\hat\gamma^{-k}),\quad
\|\phi_3\|_{C^1}=o(\hat\lambda^m).
\end{array}
\end{equation}

\subsubsection{Mutually activating blenders}
For the proof of the next proposition, we will consider only such $m$ that 
$\sin(m\omega_2+\eta_3)$ is close to zero. In this case, the condition of Lemma \ref{lem:df:F_km} that $\cos m\omega_2$ and
 $(\cos m\omega_2+a_{14}\sin m\omega_2)$ are both bounded away from zero  are automatically satisfied, as follows from \eqref{eq:df:T_k,m_cross_const1}.

\begin{prop}\label{prop:df_blender}
Let the system $f$ have a heterodimensional cycle $\Gamma$ of double-focus type, and $\theta,\frac{\omega_1}{2\pi}, \theta\frac{\omega_2}{2\pi}$ and 1 are rationally independent. There exist, arbitrarily close to $\Gamma$, a cs-blender with an activating pair $(\Pi',\mathcal{C}^{ss})$ and a cu-blender with an activating pair $(\Pi',\mathcal{C}^{uu})$ where
\begin{equation}\label{eq:Pi'_df}
\Pi'=[-q\delta,q\delta]\times[-\delta,\delta]\times[-\delta,\delta]^{d-d_1-1}\times [-\delta,\delta]^{d_1-1}
\end{equation}
for some fixed $q\in(0,1)$, and 
\begin{align*}
\mathcal{C}^{ss}&=\{(\Delta U_1,\Delta U_2,\Delta V,\Delta W):\max\{|\Delta U_1|,\|(\Delta U_2,\Delta W)\|\}\leq K  \|\Delta V\|\},\\
\mathcal{C}^{uu}&= \{(\Delta U_1,\Delta U_2,\Delta V,\Delta W): \max\{|\Delta U_1|,\|\Delta V\|\}\leq K  \|(\Delta U_2,\Delta W)\|\}
\end{align*}
for some sufficiently small $K$. The two blenders mutually activate each other (in the sense of Proposition \ref{prop:sad-foc_blender}).
\end{prop}
\begin{proof}
Denote 
\begin{equation}\label{akmdf}
\begin{array}{l}
A_{k,m}=\lambda^k\gamma^m \dfrac{C}{\sqrt{1+a_{14}^2}}\sin(k\omega_1+\eta_1),\\
B_{k,m}=\lambda^k\gamma^m \dfrac{D}{\sqrt{1+a_{14}^2}}\sin(k\omega_1+\eta_2)- 
u^-_1.
\end{array}
\end{equation}
We see that if we replace
\begin{equation}\label{eq:crspds}
U_1\to X_1,\qquad (U_2,W)\to Y,\qquad V\to(X_2,Z),
\end{equation}
then the first-return map \eqref{eq:df:T_km_cross} takes the same form as the first-return map \eqref{eq:complex:T_k,m_cross_0mu} in the saddle-focus case, provided 
we consider the values of $(k,m)$ such that
$\sin(m\omega_2+\eta_3)\to 0$ and $\lambda^k\gamma^m$ stays bounded.
Moreover, formulas \eqref{akmdf} coincide with formulas \eqref{akmbkm} for the saddle-focus case if we rename the constants $A=\frac{C}{\sqrt{1+a_{14}^2}}$,
$B=\frac{D}{\sqrt{1+a_{14}^2}}$, and replace  $u_1^-$
by $b_{11}u^-$, and $\omega_1$ by $\omega$.

Thus, we obtain Proposition \ref{prop:df_blender} in the same way as Proposition
\ref{prop:sad-foc_blender}, if we show that given any $(s,t)\in \mathbb R^2$, one can find a sequence $\{(k_n,m_n,p_n)\}$ with $k_n,m_n\to\infty$  such that
conditions \eqref{knmnpn} are satisfied, and
$$\sin(m_n\omega_2+\eta_3)\to 0$$
(we do not need the evenness of $m_n$ here, as $\gamma>0$ in the double-focus case).

These requirements are equivalent to the existence of a sequence of integers
$\{(k_n,m_n,p_n,l_n)\}$ such that
$$
k_n\theta -  m_n\to -t, \qquad k_n\dfrac{\omega_1}{2\pi}-p_n \to s,
\qquad k_n \frac{1}{\pi} \theta\omega_2 -  l_n\to 
- \frac{t\omega_2 + \eta_3}{\pi}.$$
This is guaranteed by our assumption that $\theta$, $\omega_1/2\pi$,
$\theta\omega_2/2\pi$, and $1$ are rationally independent.
\end{proof}

\subsubsection{Local stabilization of heterodimensional cycles in the double-focus case}
To complete the proof of Theorem \ref{thm:sad-foc_main} we need to show that the two periodic points $O_1$ and $O_2$ are homoclinically related to the cs-blender and, respectively, the cu-blender obtained in the above proposition.

\begin{lem}\label{lem:df_includeO_1}
Let $\Pi'$ and $\mathcal{C}^{uu}$ be given by Proposition \ref{prop:df_blender}.
There exists a sequence $\{m_j\}$ of positive integers with $m_j\to \infty$ such that the image $S^u_j:=F^{m_j}_2\circ F_{12}(W^u_{loc}(O_1))\cap \Pi$ is a disc of the form $(U_1,V)=s(U_2,W)$ for some smooth function $s$. The disc $S^u_j$ crosses $\Pi'$ properly with respect to the cone field $\mathcal{C}^{uu}$.
\end{lem}

\begin{proof}
By \eqref{eq:df:12} and \eqref{eq:df:coortrans}, the image $F^m_2\circ F_{12}(W^u_{loc}(O_1))$ satisfies
\begin{align*}
(\cos m\omega_2 +a_{14}\sin m\omega_2)(U_1+u^-_1)
&=(a_{14}\cos m\omega_2 -\sin m\omega_2)(\delta^\frac{1}{2}U_2+u^-_2+b_{41}U_1-b_{42}\delta^\frac{1}{2} V)+o(1)_{m\to\infty},\\
V&=\delta^{-\frac{1}{2}}o(\hat\lambda^m),
\end{align*}
which by \eqref{eq:df:T_k,m_cross_const1} reduces to
\begin{align*}
(1+b_{41}\tan(m\omega_2+\eta_3))U_1
&=-u^-_1 
- \tan(m\omega_2+\eta_3) u^-_2 
-\tan(m\omega_2+\eta_3)\delta^{\frac{1}{2}}U_2+o(1)_{m\to\infty},\\
V&=\delta^{-\frac{1}{2}}o(\hat\lambda^m).
\end{align*}

Now consider a sequence $\{m_j\}\to \infty$ such that $\tan(m_j\omega_2+\eta_3)\to -u^-_1/u^-_2$. Since $1/b_{41}\neq u^-_1/u^-_2$ by condition C4.2 and \eqref{eq:df:F_21_cross}, it follows that the above equations can be rewritten as
\begin{align*}
U_1
&=-\dfrac{\tan(m_j\omega_2+\eta_3)}{1+b_{41}\tan(m_j\omega_2+\eta_3)}\delta^{\frac{1}{2}}U_2+o(1)_{j\to\infty},\\
V&=\delta^{-\frac{1}{2}}o(\hat\lambda^{m_j}).
\end{align*}
This immediately shows that for all sufficiently small $\delta$ and all sufficiently large $j$, and for all $(U_2,W)\in[-\delta,\delta]\times[-\delta,\delta]^{d_1-2}$ we have
\begin{align*}
|U_1|<\dfrac{1}{2}\delta,\qquad \|V\|<\dfrac{1}{2}\delta,\qquad \dfrac{\partial(U_1,V)}{\partial(U_2,W)}=O(\delta^{\frac{1}{2}})+o(1)_{j\to 0}.
\end{align*}
The proper crossing of $S^u_j$ with respect to $\mathcal{C}^{uu}$ is thus guaranteed.
\end{proof}

\begin{lem}\label{lem:df_includeO_2} 
Let $\Pi'$ and $\mathcal{C}^{ss}$ be given by Proposition \ref{prop:df_blender}.
There exists a sequence $\{k_j\}$ of positive integers with $k_j\to \infty$ such that the preimage $S^s_j:=F^{-1}_{21}\circ F^{-k_j}_{1}\circ F^{-1}_{12}(W^s_{loc}(O_2))\cap\Pi$ is a disc of the form $(U_1,U_2,W)=s(V)$ for some smooth function $s$. The disc $S^s_j$ crosses $\Pi'$ properly with respect to the cone field $\mathcal{C}^{ss}$.
\end{lem}

\begin{proof}
By \eqref{eq:df:F_12_cross}, the preimage $F^{-1}_{12}(W^s_{loc}(O_2))$ satisfies
\begin{equation*}\label{eq:df:include1}
\begin{aligned}
0&= a_{11} \tilde x_{1}+a_{12}\tilde x_{2}+a_{13}\tilde z + O(\tilde x^2_1+\tilde x^2_2+\tilde z^2), \\
\tilde y-y^- &=a_{31}\tilde x_{1}+a_{32}\tilde x_{2}+a_{33}\tilde z + O(\tilde x^2_1+\tilde x^2_2+\tilde z^2). 
\end{aligned}
\end{equation*}
Substitute the equations for $\tilde x_1,\tilde x_2,\tilde z$ from 
formula \eqref{eq:df:21} for the map $F_1^k\circ F_{21}$ into the second equation above. This yields $\tilde y=y^-+O(\lambda^k)$ as a function of $(\tilde u_1,\tilde v)$. Combining the first equation above with \eqref{eq:df:21} and using the new expression for $\tilde y$, we obtain the following equation for $F^{-1}_{21}\circ F^{-k}_{1}\circ F^{-1}_{12}(W^s_{loc}(O_2))$ (see \eqref{ckdkeq} and \eqref{eq:df:T_k,m_cross_const0}):
\begin{equation*}
\begin{aligned}
\tilde u_1 -u^-_1&=-  \dfrac{ D\sin(k\omega_1+\eta_2)}{C\sin(k\omega_1+\eta_1)}+O((\tilde u_1-u^-)^2 + \|\tilde v\|)+o(1)_{k\to\infty}\\
\tilde u_2-u^-_2&=b_{41}(\tilde u_1-u_1^-)+b_{42}\tilde v+O((\tilde u_1-u^-)^2 + \tilde v^2)+o(\hat\gamma^{-k}),\\
\tilde w-w^-_2&=b_{51}(\tilde u_1-u_1^-)+b_{52}\tilde v+O((\tilde u_1-u^-)^2 + \tilde v^2)+o(\hat\gamma^{-k}),
\end{aligned}
\end{equation*}
where the small terms are functions of $\tilde u_1$ and $\tilde v$. Since $\omega_1$ is irrational, we can choose the sequence of values of $k\to\infty$ such that $\sin(k\omega_1+\eta_1)$ stays bounded away from zero.

Now, after solving $\tilde u_1$ from the first equation, we apply the coordinate transformation \eqref{eq:df:coortrans} and obtain the intersection $F^{-1}_{21}\circ F^{-k}_{1}\circ F^{-1}_{12}(W^s_{loc}(O_2))\cap\Pi$ as
\begin{equation*}
\begin{array}{l}
U_1=-  \dfrac{ D\sin(k\omega_1+\eta_2)}{C\sin(k\omega_1+\eta_1)}+O(\delta^{\frac{1}{2}}V) +o(1)_{k\to\infty},\\
U_2=O(\delta^{-\frac{1}{2}}U^2_1+\delta^\frac{1}{2} |V|)+ o(1)_{k\to\infty},\qquad
W=O(U^2_1+\delta |V|)+ o(1)_{k\to\infty}.
\end{array}
\end{equation*}
After comparing this with \eqref{eq:sf_includeO_2:1} using the correspondence \eqref{eq:crspds}, one just follows the proof of Lemma \ref{lem:sf_includeO_2}.
\end{proof}

It can be seen from \eqref{eq:df:coortrans} that $W^u_{loc}(O_2)=\{V=0\}$, and, after additionally taking $k\to\infty$ in \eqref{eq:df:21}, that the intersection $F^{-1}_{21}(W^s_{loc}(O_1))\cap\Pi$ is the graph of some smooth function $(U_2,W)=s(U_1,V)$ satisfying $s=O(U^2_1+\delta V^2)$ and defined on $[-\delta,\delta]\times[-\delta,\delta]^{d-d_1-1}$.

Therefore, since the above two lemmas are completely analogous to Lemmas \ref{lem:sf:includeO_1} and \ref{lem:sf_includeO_2},
one obtains the homoclinic relations between $O_1$ and the cs-blender, and between $O_2$ and the cu-blender in the same way as it is shown in the end of Section \ref{sec:sta_sad-foc}.

\section*{Acknowledgments}
\addcontentsline{toc}{section}{Acknowledgments}
The authors thank the referees for their advice. We are also grateful to Katsutoshi Shinohara and Pierre Berger for important discussions during the preparation of this paper. 
This work was supported by the Leverhulme Trust grant RPG 2021-072 and grants 19-11-00280 and 19-71-10048 of RSF. The research of D. Li was also supported by the ERC project 677793 StableChaoticPlanetM.

\section*{Appendix: Standard blenders}\label{sec:appen}
\addcontentsline{toc}{section}{Appendix: Standard blenders}

Here we give a detailed construction of standard blenders. Essentially, this is a summary (with a necessary generalization to the case of an arbitrary Markov partition) of the blender horseshoe construction \cite{BD12a}. We also show in Proposition \ref{prop:stanble} that the blenders found in this paper are standard. The main purpose is to make our blender-related results more accessible for further use. As an example, see the application in \citep{LLST22} where we show the $C^1$-robustness of homoclinic tangencies near homoclinic tangencies of effective dimension larger than one. 

Recall that a compact invariant set $\Lambda$ of a diffeomorphism $g$ is hyperbolic if it  possesses a continuous $Dg$-invariant splitting of its tangent bundle: $\mathcal{T}_{\Lambda}\mathcal{M}=E^{s}\oplus E^{u}$, such that $Dg$ is expanding in $E^u$ and contracting in $E^s$. This implies the existence of invariant continuous cone fields $\mathcal{C}^{s}$ and $\mathcal{C}^{u}$ in a small neighborhood $U$ of $\Lambda$: for any point $M\in U$, and any vectors
$v\in cl(\mathcal{C}^{s}_M)\subset \mathcal{T}_M U$ and $w\in cl(\mathcal{C}^{u}_M)\subset \mathcal{T}_M U$, one has $D g^{-1}(v)\in int(\mathcal{C}^{s}_{g^{-1}(M)})$ if $g^{-1}(M)\in U$, and $D g(v)\in int(\mathcal{C}^{u}_{g(M)})$ if $g(M)\in U$; moreover, $Dg$ is uniformly expanding in $\mathcal{C}^u$ and uniformly contracting in $\mathcal{C}^s$.

The hyperbolic set $\Lambda$ is cs-partially hyperbolic if  
the invariant subbundle $E^s$ admits a further invariant splitting $E^{s}=E^{ss}\oplus E^{cs}$ with $E^{cs}\neq 0$, where the contraction in $E^{ss}$ is uniformly stronger than in $E^{cs}$.
Similarly, it is cu-partially hyperbolic if  
the invariant subbundle $E^u$ admits an invariant splitting $E^{u}=E^{uu}\oplus E^{cu}$ with $E^{cu}\neq 0$, where the expansion in $E^{uu}$ is uniformly stronger than in $E^{cu}$.
The partial hyperbolicity implies the existence of  a forward-invariant cone field $\mathcal{C}^{uu}\subset \mathcal{C}^u$ that contains $E^{uu}$ or a backward-invariant cone field $\mathcal{C}^{ss}\subset \mathcal{C}^s$ that contains $E^{ss}$.

The compact hyperbolic set $\Lambda$ is called {\em basic} if it is transitive and locally-maximal, i.e.,
there is an orbit which is dense in $\Lambda$ and, when the neighborhood $U$ of $\Lambda$ is sufficiently small, any orbit which lies entirely in $U$ lies in $\Lambda$. Note that any point whose forward orbit never leaves $U$ belongs to the local stable manifold of some point in $\Lambda$, and
any point whose backward orbit never leaves $U$ belongs to the local unstable manifold of some point in $\Lambda$.

We specifically consider {\em zero-dimensional} basic sets. Such sets admit a description in terms of finite Markov chains. Namely, one can find a finite unidirectional graph $G$ with $k$ vertices
such that the set of orbits in $\Lambda$ is in 1-to-1 correspondence with the set of all infinite paths along the edges of $G$. The correspondence is established as follows. One finds a
finite collection of disjoint open subsets $U_1,\dots,U_n$ of $U$ (a Markov partition) such that
$\Lambda \subset U_1 \cup \dots \cup U_n$, and
$g(cl(U_s))\cap cl(U_j) \neq \emptyset$ if and only if there is an edge in the graph $G$ that goes from the vertex $s$ to the vertex $j$. Any orbit in $\Lambda$ defines a coding sequence $\{j_i\}_{i\in \mathbb{Z}} \in \{1,\dots,n\}^{\mathbb Z}$ such that the $i$-th point of the orbit
lies in $U_{j_i}$, for all $i$. Moreover, the partition $\{U_1, \dots, U_n\}$ can be chosen such that these sequences encode orbits in $\Lambda$ uniquely (i.e., different orbits produce different sequences $\{j_i\}$) and a sequence $\{j_i\}_{i\in \mathbb{Z}}$ corresponds to some orbit of $\Lambda$ if and only if the graph $G$ has an edge from the vertex $j_{i}$ to the vertex $j_{i+1}$  for all $i$.

One can choose the Markov partition such that $cl(U_j)$ is, for each $j$, a diffeomorphic image of a Cartesian product
$\mathbb{X}_j\times \mathbb{Y}_j\subset \mathbb{R}^{d^s}\times \mathbb{R}^{d^u}$ of two compact convex sets, where $d^s=\dim (E^s)$ and $d^u=\dim (E^u)$. Hence, we can use $x\in \mathbb{R}^{d^s}$ and $y\in \mathbb{R}^{d^u}$ as coordinates on $cl(U_j)$. Let us denote by $dx$ and $dy$ the vectors in the tangent spaces. The invariant cone fields in $U_1\cup \dots \cup U_n$ are then given by $\mathcal{C}^{u}= \{\|dx\| < K_j^u \|dy\|\}$ and
$\mathcal{C}^{s}= \{\|dy\| < K_j^s\|dx\|\}$ (where the numbers $K_j^{s,u}$ may depend on the point in $\mathbb{X}_j\times \mathbb{Y}_j$).
Moreover, given any {\em complete u-disc} $S\subset cl(U_j)$, i.e., a graph of a smooth function $x=h(y)$ defined {\em for all}
$y\in \mathbb{Y}_j$, taking values $x\in int(\mathbb{X}_j)$ for some $j$, and satisfying the Lipshitz property 
$(D h(dy), dy) \in \mathcal{C}^{s}$, if there is an edge of the graph $G$ going from $j$ to $s$, then the image $g(S\cap g^{-1}(cl(U_s)))$ is a complete u-disc in $cl(U_s)$,
i.e., it is a graph of a smooth function satisfying the the same Liprshitz property and defined for all
$y\in \mathbb{Y}_s$. Similarly, given any {\em complete s-disc} $S\subset cl(U_s)$, i.e., a graph of a smooth function $y=h(x)$ defined for all
$x\in \mathbb{X}_s$, taking values $y\in int(\mathbb{Y}_s)$, and satisfying the Lipshitz property $(dx, Dh(dx)) \in \mathcal{C}^{u}$, the preimage $g^{-1}(S\cap g(cl(U_j)))$ is a complete s-disc in $cl(U_j)$.

Since the tangents to the s-discs lie in $\mathcal{C}^{s}$ and the
tangents to the u-discs lie in $\mathcal{C}^{u}$, the map $g$ is uniformly expanding on the u-discs and uniformly contracting on the s-discs. The stable
lamination $W^s_{loc}(\Lambda)$ is comprised by complete s-discs and the unstable lamination $W^u_{loc}(\Lambda)$ by complete u-discs. Moreover,
$W^s_{loc}(\Lambda)=\bigcup_{J} W^s_{loc,J}$
where the union is taken over all forward paths $J=\{j_i\}_{i\geq 0}$  in the graph $G$ and  $W^s_{loc,J}$ is the uniquely defined complete s-disc that consists of all points $M\in cl(U_{j_0})$
such that $g^i(M)\in cl(U_{j_i})$ for all $i\geq 0$. Similarly,
$W^u_{loc}(\Lambda)=\bigcup_{J} W^u_{loc,J}$
where the union is taken over all backward paths $J=\{j_i\}_{i\leq 0}$ in the graph $G$ and  $W^u_{loc,J}$ is the uniquely defined completely u-disc that consists of all points $M\in cl(U_{j_0})$ such that $g^i(M)\in cl(U_{j_i})$ for all $i\leq 0$.

If the zero-dimensional basic set is also cs-partially hyperbolic, we assume that the Markov partition is chosen such that $\mathbb{X}_j = \mathbb{X}^{ss}_j\times \mathbb{X}^{cs}_j$, where
$\mathbb{X}^{ss}_j$ and $\mathbb{X}^{cs}_j$ are convex compact sets of dimension
$\dim(E^{ss})$ and $\dim(E^{cs})$, respectively. Thus, a point in $cl(U_j)$ has coordinates
$(x^{ss},x^{cs}, y)$ with $x^{ss}\in \mathbb{X}^{ss}_j$, $x^{cs}_j\in\mathbb{X}^{cs}_j$,
$y\in \mathbb{Y}_j$. A {\em complete ss-disc} of $U_j$ is the graph
$\ell:=(x^{cs},y)=p(x^{ss})$ of a smooth function $p$ defined {\em for all}
$x^{ss}\in \mathbb{X}^{ss}_j$  and such that the tangents to $\ell$ lie in
$\mathcal{C}^{ss}$ at every point of $\ell$. 

We say that the {\em covering property} is satisfied if for any complete ss-disc $\ell$ there exists $s$ such that the preimage $g^{-1}(\ell \cap g(cl(U_s))$ is a complete ss-disc again.

\begin{defia}[Standard blenders]\label{defi:blender}
A basic zero-dimensional cs-partially hyperbolic set $\Lambda$ is called a  standard cs-blender of coindex $d=\dim(E^{cs})$ if it satisfies the covering property.
\end{defia}

By the covering property, given any complete ss-disc, we have an infinite sequence of its  preimages by $g$ which lie in a small neighborhood $U$ of a cs-blender $\Lambda$. So, every such disc has
at least one point whose backward orbit never leaves $U$. Such points must belong to
$W_{loc}^u(\Lambda)$, which means that {\em every complete ss-disc has a non-empty intersection with the local unstable lamination of the standard cs-blender}. 

Since the sum of the dimensions of any leaf $W^u_{loc,J}$ of $W^u_{loc}(\Lambda)$ and any ss-disc is strictly less than $\dim(U)$, these intersections are non-transverse, so they get destroyed by arbitrarily small perturbations. However, one can easily show that the covering property is $C^1$-open: the basic set $\Lambda$ persists for all $C^1$-small perturbations of $g$, i.e., it depends continuously on $g$, the Markov partition and the cone families $\mathcal{C}^{ss}\subset \mathcal{C}^s$ and
$\mathcal{C}^u$ remain unchanged, and the covering property holds for the continuation of
$\Lambda$. It follows that, given any complete ss-disc, it has an intersection with some leaf of
$W^s_{loc}(\Lambda)$ for every map which is $C^1$-close to $g$. Thus, a standard cs-blender is a cs-blender in the usual terminology \citep{BD96,BDV,BCDW}, and, particularly, satisfies Definition \ref{defi:blender_old}.

The set $\Lambda$ is a standard cu-blender if it is a standard cs-blender of the system obtained from time inversion.  More specifically, when $\Lambda$ is cu-partially hyperbolic, one chooses the Markov partition such that $\mathbb{Y}_j = \mathbb{Y}^{uu}_j\times \mathbb{Y}^{cu}_j$, and defines {\em complete uu-discs} as graphs of functions of the form $(x,y^{cu})=p(y^{uu})$. The set $\Lambda$ is called a {\em standard cu-blender of coindex $d=\dim(E^{cu})$} if it satisfies the following covering property:
for any complete uu-disc $\ell$ there exists $s$ such that the preimage
$g(\ell \cap g^{-1}(cl(U_s))$ is a complete uu-disc again. Then {\em every complete uu-disc activates (i.e., has a non-empty intersection with the local stable lamination)  the standard cu-blender
$\Lambda$} and this property persists for every $C^1$-close map.

\begin{rema}\label{rem:a1}
Note that the same construction can be carried out in the case of continuous-time dynamical systems. Given a compact, transitive, locally-maximal, one-dimensional uniformly hyperbolic invariant set $\Lambda$ of a smooth flow, we can take a cross-section $U$ such that the flow near $\Lambda$ defines, on $U$, the return map $g$ -- a diffeomorphism acting from a small neighborhood of $\Lambda\cap U$ into $U$ such that $\Lambda\cap U$ is a zero-dimensional basic set of $g$. Then, $\Lambda$ is a standard blender of the flow if and only if $\Lambda\cap U$ is a standard blender of $g$.
\end{rema}


\begin{rema}\label{rem:a2}
Similarly to the case of a flow, given a diffeomorphism $f$ and a zero-dimensional basic set $\Lambda$
of $f$, take an open set $U$ such that $\partial U \cap \Lambda = \emptyset$. Then, $\Lambda \cap U$ is a basic set of an induced map $g$ defined on the union of finitely many disjoint open subsets $U_j$ of $U$, by the restriction of certain iterations of $f$ to $U_j$.  Iterating by $f$ the Markov partition $\{U_j\}$, one extends it to a small neighborhood of $\Lambda$, and proves that if $\Lambda\cap U$ is a standard blender for $g$, then $\Lambda$ is a standard blender of $f$.
\end{rema}

Finally, let us show that the blenders obtained in this paper are standard. By the symmetry of the problem, it suffices to consider only cs-blenders.  Recall that the finite set $\mathcal{J}_{\delta'}$ of pairs in Lemma \ref{lem:crossingsurfaces} corresponds to a hyperbolic basic set $\Lambda:=\Lambda_{\mathcal{J}_{\delta'}}$ by Lemma \ref{lem:hypset2}, and it is just the cs-blender in Proposition \ref{prop:cublendertype1}.
\begin{propa}\label{prop:stanble}
The cs-blender $\Lambda$ is  standard.
\end{propa}
\begin{proof}
Let $\Pi'$ be the cube of size $\delta'$, as in Lemma \ref{lem:crossingsurfaces}. For each pair $(k_j,m_j)\in \mathcal{J_{\delta'}} $, the intersection $T^{-1}_{k_j,m_j}(\Pi')\cap \Pi'$ is a ``horizontal strip'' of the form
\begin{equation*}\label{eq:strips}
U_j=\{(X,Y,Z)\mid X\in [-\delta',\delta'], Y\in \phi_{m_j,k_j}(X, [-\delta',\delta']^{d_1},Z),Z\in [-\delta',\delta']^{d-d_1-1}  \},
\end{equation*}
where $\phi_{m_j,k_j}$ is the function $\phi_2=o(\hat\gamma^{-k})$ from \eqref{eq:maps:T_k,m_cross_0mu} (take $(k,m)=(k_j,m_j)$). These strips are pairwise disjoint, and  their union contains $\Lambda\cap\Pi'$. If we define a map $g$ on $\bigcup U_j$ by
$$g(M)=T_{k_j,m_j}(M)\quad {if} \quad M\in int(\sigma_{k_j,m_j}),$$
then the set $\Lambda\cap\Pi'$ is a partially hyperbolic basic set of $g$. Moreover, it is cs-partially hyperbolic, with the cone fields  ${\mathcal{C}}^{cs}$, ${\mathcal{C}}^{ss}$ and ${\mathcal{C}}^{uu}$ described in Lemma \ref{lem:conefields}.

By Remarks \ref{rem:a1} and \ref{rem:a2}, it suffices to prove that $\Lambda\cap \Pi'$ is a standard cs-blender of $g$.
Obviously,  the sets $U_j$ form a Markov partition for $\Lambda\cap\Pi'$ with the property that for any $s$ and $j$ we have $g(cl(U_s))\cap cl(U_j)$ non-empty. Next, observe that every $cl(U_j)$ is the diffeomorphic image of
$$Q=\{(x^{cs},x^{ss},y)\mid x^{cs}\in [-\delta',\delta'], x^{ss}\in [-\delta',\delta']^{d-d_1-1} , Y\in [-\delta',\delta']^{d_1}
 \},$$
where the diffeomorphism is defined by
$$X = x^{cs},\quad Y = \phi_{m_j,k_j}(x^{cs},y,x^{ss}),\quad Z=x^{ss}. $$
The last thing to check is the covering property, but this is exactly Lemma \ref{lem:crossingsurfaces}. So $\Lambda\cap \Pi'$ is indeed a standard cs-blender of $g$.
\end{proof}

\end{document}